\documentclass[12pt]{amsart}

\usepackage[T1]{fontenc}
\usepackage[letterpaper,margin=1in]{geometry}
\usepackage{microtype}
\usepackage{graphicx}
\usepackage{subcaption}
\usepackage{float}
\usepackage{array,booktabs,tabularx}
\usepackage{comment}
\usepackage{xcolor}

\usepackage{amsmath,amssymb,amsfonts,amsthm}
\usepackage{mathrsfs}
\usepackage{mathtools}
\usepackage{bm}
\usepackage{thmtools,thm-restate}
\usepackage[shortlabels]{enumitem}

\declaretheorem[numberwithin=section]{theorem}
\declaretheorem[sibling=theorem]{corollary,lemma,proposition,question,conjecture}
\declaretheorem[style=definition,sibling=theorem]{definition,example}
\declaretheorem[style=remark,sibling=theorem]{remark,philosophy}
\numberwithin{equation}{section}

\DeclareMathOperator{\Fav}{Fav}
\DeclareMathOperator{\proj}{proj}
\DeclareMathOperator{\dist}{dist}

\DeclareMathOperator{\lcm}{lcm}

\DeclareMathOperator{\sign}{sign}

\allowdisplaybreaks

\newcommand{\cH}{\mathcal{H}}

\newcommand{\cM}{\mathcal{M}}

\newcommand{\RR}{\mathbb R}
\newcommand{\Q}{\mathbb Q}

\renewcommand{\j}[1]{\langle #1 \rangle}

\newcommand{\capexp}{\hbox{\sf EXP}}

\def\set4{\mathcal I}
\def\tup14{(1,2,3,4)}

\usepackage{graphicx}
\usepackage{mathrsfs}

\newcommand{\R}{\mathbb{R}}

\newcommand{\I}{\mathbb I}

\newcommand{\calC}{\mathcal{C}}

\usepackage[unicode=true,colorlinks=true,linkcolor=blue,citecolor=blue,urlcolor=blue,hypertexnames=false]{hyperref}

\begin{document}

\author[C.~Marshall]{Caleb Marshall}
\address{Department of Mathematics \\ The University of Toronto, St. George Campus, Toronto, ON, Canada}
\email{caleb.marshall@utoronto.ca}

\keywords{Favard length, rectifiability, orthogonal projections, self-similar sets, Riesz products, cyclotomic polynomials}
\subjclass[2020]{28A78, 28A80 (primary); 11C08, 42B05, 11B75 (secondary)}

\date{July 2026}
\title[Improved Power Laws for Favard Length in All Dimensions]{Improved Power Laws for the Favard Length Problem in All Dimensions}
\dedicatory{In loving memory of my Mother, J.L.L. Marshall.}

\begin{abstract}
The Favard length of a compact set in $\mathbb R^d$ is the average one-dimensional Hausdorff measure of its orthogonal projections onto lines.  We establish quantitative upper bounds for small neighbourhoods of a class of purely $1$-unrectifiable rational product Cantor sets in every dimension $d\geq2$.  The proof combines coordinatewise small-value estimates, multidimensional sets of large values, and an improved symbolic-cylinder combinatorial argument whose low-multiplicity estimate closes for every parameter $\rho>2$.  Applied to the classical four-corner Cantor set, the original Nazarov--Peres--Volberg Fourier input and the improved propagation yield the bound $N^{-1/5+u}$ for every $u>0$.  When the coordinate mask polynomials are zero-free on the unit circle, we obtain the general exponent $(2D+1)^{-1}-u$, where $D:=\sum_i\log_L\bigl(\#A_i/\min_{|z|=1}|A_i(z)|\bigr)$; in particular, explicit planar examples in bases $25$ and $36$ attain the exponents $1/3-u$ and $\delta_{36}-u$, respectively, where $\delta_{36}=1/(2\log_6(300/53)+1)=0.340720\ldots>1/3$.
\end{abstract}

\maketitle

\section{Introduction}
Assuming $d\geq 2$, we define for each $x \in \mathbb{R}^d$ and $\theta\in\mathbb{S}^{d-1}$, the dot product mapping
$
\pi_\theta(x):=x\cdot\theta,
$
which we identify with orthogonal projection onto the one-dimensional linear subspace spanned by $\theta$.  Letting $\sigma^{d-1}$ denote the normalized surface measure on $\mathbb{S}^{d-1}$, we then define the \emph{Favard length} of a compact set $E\subset\mathbb{R}^d$ as
\[
\Fav(E):=\int_{\mathbb{S}^{d-1}}\mathcal{H}^1\big(\pi_\theta(E)\big)\,d\sigma^{d-1}(\theta).
\]
This Favard length of $E$ is thus the average length of the orthogonal projections of $E$ onto one-dimensional linear subspaces.

The Besicovitch--Federer projection theorem (see \cite[Theorem 18.1]{mattila}) states that a compact set $E\subset\mathbb{R}^d$ satisfying $0<\mathcal{H}^1(E)<\infty$ is purely $1$-unrectifiable if and only if $\Fav(E)=0$.  If
\[
\mathcal{N}_\delta(E):=\{x\in\mathbb{R}^d:\dist(x,E)<\delta\},
\]
then compactness, monotonicity, and the dominated convergence theorem show that this qualitative statement is equivalent to the existence of a monotone function $\psi_E:(0,1)\to(0,\infty)$ such that
$
\Fav\big(\mathcal{N}_\delta(E)\big)\leq \psi_E(\delta)
$
and
$
\lim_{\delta\to0^+}\psi_E(\delta)=0.
$
The Favard length problem asks for quantitative information about the decay of this asymptotic upper bound $\psi_E$.

We establish upper bounds for two overlapping, countably infinite families of rational product Cantor sets.  The most general statement is Theorem \ref{thm:totalFavbound}.  Before introducing this result and its framework, we describe a canonical higher-dimensional example and state our corresponding upper bound.

\subsection{A higher-dimensional variant of the four-corner Cantor set}
Let $d\geq2$ and set $A:=\{0,2^d-1\}$.  For $N\geq1$, define
\begin{equation}\label{eq:2dmultidigit}
A^N:=2^{-d}A+2^{-2d}A+\cdots+2^{-Nd}A,
\end{equation}
and let
\begin{equation}\label{eq:2diterates}
\mathcal{K}_d^N:=\underbrace{(A^N\times\cdots\times A^N)}_{d\text{ factors}}+[0,2^{-dN}]^d.
\end{equation}
The set $\mathcal{K}_d^N$ is a union of $2^{dN}$ hypercubes of side-length $2^{-dN}$.  Equivalently, one begins with $[0,1]^d$, subdivides it into hypercubes of side-length $2^{-d}$, retains only the $2^d$ subcubes meeting a vertex, and repeats the same operation inside every retained cube.  The limiting set
\begin{equation}\label{eq:2dcantorset}
\mathcal{K}_d^\infty:=\bigcap_{N=1}^\infty \mathcal{K}_d^N
\end{equation}
is the \emph{$2^d$-corner Cantor set}.

The set $\mathcal{K}_d^\infty$ is compact, self-similar, and satisfies the strong separation condition.  Its similarity dimension is one, and it is purely $1$-unrectifiable; consequently,
\[
0<\mathcal{H}^1(\mathcal{K}_d^\infty)<\infty
\quad\text{and}\quad
\Fav(\mathcal{K}_d^\infty)=0.
\]
When $d=2$, this construction is the classical four-corner Cantor set; see \cite{Tay} for an introduction to its role in projection problems.

\begin{theorem}\label{thm:mainsimple}
Let $d\geq2$.  There exist $\epsilon_d\in(0,1)$ and $C_d\geq1$ such that, for every sufficiently large $N \in \mathbb{N}$,
\begin{equation}\label{eq:thmmainsimple}
C_d^{-1}N^{-1}
\leq
\Fav\big(\mathcal{N}_{2^{-dN}}(\mathcal{K}_d^\infty)\big)
\leq
C_dN^{-\epsilon_d}.
\end{equation}
\end{theorem}

The lower bound follows from Mattila's general estimate \cite[Theorem 4.1]{Mattila}; the upper bound is our new contribution.  For $d=2$, Nazarov, Peres, and Volberg proved the more precise family of bounds
\[
\Fav\big(\mathcal{N}_{4^{-N}}(\mathcal{K}_2^\infty)\big)\lesssim_u N^{-1/6+u},
\qquad u>0,
\]
in \cite{NPV}, with $\mathcal{K}_2^{\infty} \subset \mathbb{R}^2$ colloquially known as the \textit{four-corner Cantor set} or \textit{Garnett set}. As part of our method, we are able to concretely improve this exponent to the following bound.

\begin{theorem}
For every $u>0$, there are constants $C_u$ and $N_u$ such that, for every $N\geq N_u$,
\[
\Fav\big(\mathcal{N}_{4^{-N}}(\mathcal{K}_2^\infty)\big)\leq C_uN^{-1/5+u}.
\]
\end{theorem}

Our proof of this improvement relies on several general improvements of both the Fourier analytic and combinatorial methods originally posited by Nazarov, Peres and Volberg in \cite{NPV}, as well as careful book-keeping of all implicit parameters and implied constants. The improvement of the exponent for the four-corner Cantor set follows directly from our general upper bound machinery via a direct and simple calculation. To further highlight the strength of our method, we are also able to improve the exponent $1/5 - u$ to $1/3 - u$ for certain special classes of Cantor sets (with the specific result stated as Corollary \ref{cor:zerofreedigitsets}). The specific calculation of this exponent, as well as the general examples with improved exponents, will be given in Section \ref{sec:exponentexamples}.

If $E \subset \mathbb{R}^d$ is a compact and purely $1$-unrectifiable set, then an upper bound of the form
$$
\Fav \big(\mathcal{N}_{\delta} (E)\big) \lesssim \bigl(\log(\delta^{-1})\bigr)^{-\epsilon}, \textrm{ for all } \delta \ll 1,
$$
is colloquially known as a \textit{power law} for the Favard length of $\delta$-neighbourhoods of $E$. Theorem \ref{thm:mainsimple} gives such a law for a genuinely higher-dimensional deterministic self-similar set: unlike a planar example embedded into $\mathbb{R}^d$, the set $\mathcal{K}_d^\infty$ is not supported on any proper affine subspace. Hence, our Theorem \ref{thm:mainsimple} is the first genuinely non-trivial example of a power law for the Favard length problem in dimensions $d \geq 3$.

\subsection{Higher dimensional estimates for more general self-similar sets}
Theorem \ref{thm:mainsimple} is a special case of the more general estimate in Theorem \ref{thm:totalFavbound}.  The relevant sets are rational product Cantor sets whose coordinate digit sets have controlled accumulating cyclotomic factors.  Different coordinates may satisfy different hypotheses: small cardinality, a two-prime condition, or the existence of a fibered subset.

The planar method originating in the work of Nazarov, Peres, and Volberg \cite{NPV} was extended by Bond, {\L}aba, and Volberg \cite{BLV}, and was further enlarged by {\L}aba and the author through the analysis of vanishing sums of roots of unity \cite{LM}.  The present argument not only improves upon the Fourier analytic and combinatorial mechanisms developed in these previous works; it also extends these techniques to every ambient dimension $d\geq2$.  As a further extension, we also show that our method applies to an even larger class of self-similar sets satisfying the cardinality bounds developed by Kiss, {\L}aba, Somlai, and the author \cite{KLMS}; this part is new even in the plane.

The broader planar literature includes \cite{BV,BThesis,BLV,BV1,BV3,Bongers,CDOV,LZ,NPV,Tao,VV,Wilson,Zhang}; surveys and expository accounts include \cite{L1,PS1}.  Connections with analytic capacity and quantitative rectifiability appear in \cite{DV,To,Dab}, while curvilinear and nonlinear variants are studied in \cite{BV2,BT,CDT,LMT}.  The recent comparison theorem of {\L}aba, McDonald, and Taylor \cite{LMT} transfers quantitative classical Favard-length bounds to smooth nonlinear projection families; consequently, for the eligible planar self-similar sets, the exponent improvements proved here pass to their generalized Favard lengths without a loss in the power exponent.  We record the precise consequence in Section \ref{sec:exponentexamples}.  In a different, probabilistic direction, Chang, Shmerkin, and Suomala recently proved the sharp order $1/\log(1/r)$ for a broad class of planar one-dimensional random Cantor sets and exhibited other random models with an additional $\log\log(1/r)$ loss \cite{CSS}.  Their results provide useful context for the range of possible Favard decay, although the present article concerns deterministic rational product constructions and thus utilizes entirely distinct machinery.

\subsection{Main results summary}
The most general result regarding power laws in arbitrary dimension $d \geq 2$ is given as Theorem \ref{thm:totalFavbound}. It allows the coordinate hypotheses to be mixed and gives a fixed power law when every unit-circle zero is a root of unity, together with a quantitative sub-power estimate in the remaining branch.  The improved combinatorial estimate is Theorem \ref{thm:improved-L2}.  Its principal concrete consequences are the four-corner estimate in Theorem \ref{thm:fourcorneronefifth} and the zero-free-digit-set results in Corollary \ref{cor:zerofreedigitsets} and Theorem \ref{thm:explicitonethird}.

\subsection{Outline of the paper}\label{subsec:paperoutline}
Section \ref{sec:RPCSfavresults} introduces rational product Cantor sets, mask polynomials, and the associated algebraic structures necessary to establish our power laws.  Section \ref{sec:toolbox} develops the Fourier analytic tools needed for our main argument, which are: the counting function and its Fourier representation as a Riesz product; a projective-chart reduction of the sphere; and, a Fourier pigeonholing argument reducing Theorem \ref{thm:mainsimple} to an $L^2$ estimate on an almost unit interval.  Section \ref{sec:rieszproductbounds} gives an $L^2$ estimate for the Riesz-product, which is explicit in all underlying constants and parameters.   Sections \ref{sec:ssvtechnical} and \ref{sec:multiSLV} supply the small-value and large-value inputs, which are the main algebraic and combinatorial inputs for our argument (and rely on deep facts regarding divisibility by \textit{cyclotomic polynomials}).  Section \ref{sec:comblemmata} proves the improved combinatorial estimates, and Section \ref{sec:exponentexamples} performs the exponent calculation after presenting the principal examples.

\section*{Acknowledgments}
The author was supported by the Four-Year Doctoral Fellowship at The University of British Columbia, where the author was also a Vanier Fellow and Killam Scholar. The author also received the generous support of the SQuaREs program, hosted by the \textit{American Institute of Mathematics}, and this work was produced in conjunction with the ``Covering Fractals by Curves'' SQuaRE. The author expresses their gratitude to Paige Bright and Joshua Zahl for feedback on an earlier version of this manuscript, as well as Blair Davey, Krystal Taylor, and Malabika Pramanik for feedback and corrections on later versions. The author would also like to thank Alan Chang, Alex Iosevich and Ignacio Uriarte-Tuero for many helpful and encouraging discussions of ideas surrounding the Favard length problem and its variants.


\section*{Notation and Standing Assumptions}
\begin{itemize}
\item If $A$ is a finite set, we let $\# A$ denote the cardinality of $A$.

\smallskip

\item We let $\mathbb{N}_0:=\{0,1,2,\ldots\}$ and, if $M \in \mathbb{N}$, we let $[M] : = \{0,1,...,M-1\} \subset \mathbb{N}_0$.

\smallskip

\item If $X,Y \in \mathbb{R}$ satisfy $X \leq CY$ for some constant $C > 0$ which does not depend upon $X$ and $Y$, then we will write $X \lesssim Y$. If $X \lesssim Y$ and $Y \lesssim X$, then we write $X \approx Y$.  

\smallskip

\item A constant $C > 0$ depends upon a parameter $\upsilon$ if $C$ is a function of $\upsilon$. We sometimes write $C = C_{\upsilon}$ to denote this dependence.

\smallskip

\item If $X,Y \in \mathbb{R}$ satisfy $X \leq C_{\upsilon} Y$ for some $C = C_{\upsilon} > 0$, we write $X \lesssim_{\upsilon} Y$.

\smallskip

\item If $X,Y \in \mathbb{R}$ and $X$ is asymptotically dominated by $Y$, we write $X \ll Y$.

\smallskip

\item A ball $B(z,r) \subset \mathbb{R}^d$ is the open set
$
B(z,r) : = \{x \in \mathbb{R}^d : \vert x - z\vert < r\},
$
where $\lvert \cdot \rvert$ denotes the standard Euclidean distance on $\mathbb{R}^d$.
\smallskip

\item If $E \subset \mathbb{R}^d$, then the $\delta$-neighbourhood of $E$ is the open set $$
\mathcal{N}_{\delta} (E) : = \{x \in \mathbb{R}^d : \vert x - a \vert < \delta \textrm{ for some } a \in E\}.
$$

\smallskip

\item If $(X,d)$ is a metric space, $s \in \{0,...,d\}$ is an integer, and $E \subset X$ is Borel, then $\mathcal{H}^s (E)$ denotes the standard $s$-dimensional Hausdorff measure of $E$ induced by the underlying metric. 

\end{itemize}

\section{The Favard Length of Rational Product Cantor Sets}\label{sec:RPCSfavresults}

We prove Theorem \ref{thm:mainsimple} through a more general theory of rational product Cantor sets.  Subsection \ref{subsec:RPCS} fixes the geometric model and its finite-stage approximations.  Subsection \ref{subsec:digitsetpoly} introduces the digit-set-polynomial factorization that separates accumulating from non-accumulating zeros.  Subsection \ref{subsec:fiberedfav} formulates the fibered-subset extension.  The resulting estimates are collected in Theorem \ref{thm:totalFavbound}.

\subsection{Rational product Cantor sets}\label{subsec:RPCS}
We first begin with a discussion of the types of self-similar sets which we shall consider. We begin by introducing the basic structure which we require at a single-scale, which is accomplished in the following definition.

\medskip

\begin{definition}
Given $L \geq 3$ and non-empty sets $A_1,...,A_d \subset \{0,1,...,L-1\}$, we call $\mathcal{A} : = A_1 \times \cdots \times A_d$ a \textbf{(product) digit set} in \textbf{base $\bm{L}$}.
\end{definition}

The reasoning for introducing this definition is to construct a special kind of self-similar set, which is defined below, and is a major object of study throughout this article.

\medskip

\begin{definition}
Let $\mathcal{A} : = A_1 \times \cdots \times A_d$ be any product digit set in base $L \geq 3$. We then refer to the set
$$
\mathcal{S} = \mathcal{S}^{\infty}(L; \mathcal{A}) : = \big\{x \in \mathbb{R}^d : x = \sum\limits_{k=1}^{\infty} L^{-k} \alpha_k, \textrm{ for some } \bm{\alpha} : = (\alpha_k) \in \mathcal{A}^{\mathbb{N}} \big\},
$$
as the \textbf{rational product Cantor set} associated to $L$ and $A_1,...,A_d$. If $M \geq 1$ is any integer, we further let
$$
\mathcal{A}^M : = \big\{x \in \mathbb{R}^d : x = \sum\limits_{k=1}^{M} L^{-k} \alpha_k, \textrm{ for some } \bm{\alpha} : = (\alpha_k) \in \mathcal{A}^{M} \big\},
$$
and call the sets
$$
\mathcal{S}^M : = \mathcal{A}^M + [0,L^{-M}]^d
$$
the \textbf{Cantor iterates} at scales $M \geq 1$, where the sum above denotes the Minkowski sum of these two sets.
\end{definition}

The name \textit{rational product Cantor set} comes from the following geometric perspective. Let us begin with the unit hypercube $Q : = [0,1]^d \subset \mathbb{R}^d$, as well as the \textbf{base $\bm{L}$ lattice} over $Q$, which is given by
$$
\Lambda_L^1 : = \{(\frac{t_1}{L},\cdots,\frac{t_d}{L}) : t_j \in \{0,1,...,L-1\} \textrm{ for each } j =1,...,d \}.
$$
Now, if $A_1,...,A_d \subset \{0,1,...,L-1\}$ are given non-empty sets and $\mathcal{A} : = A_1 \times \cdots \times A_d$, then the set
\begin{equation}\label{ch2:eq:rationalproductdigits}
\mathcal{A}^1 := \{\big(\frac{t_1}{L},...,\frac{t_d}{L} \big) : t_j \in A_j \textrm{ for each } j = 1,...,d\}
\end{equation}
is, itself, a subset of the lattice $\Lambda_L^1$. The Cantor iterate $\mathcal{S}^1$, then, corresponds to placing a hypercube of side-length $L^{-1}$ whose minimal vertex (i.e. the vertex where each coordinate is simultaneously minimized) aligns with each point $a \in \mathcal{A}^1$. To obtain the second-scale Cantor iterate $\mathcal{S}^2$, we further refine the lattice $\Lambda_L^1 \subset \Lambda_L^2$, where
$$
\Lambda_L^2 : = \{(\frac{t_1}{L^2},\cdots,\frac{t_d}{L^2}) : t_j \in \{0,1,...,L^2-1\} \textrm{ for each } j =1,...,d \},
$$
and observe that the set
$$
\mathcal{A}^2 : = L^{-1} \mathcal{A} + L^{-2} \mathcal{A} \subset \Lambda_L^2.
$$
By placing hypercubes of side-length $L^{-2}$ whose minimal vertices align with each point $a \in \mathcal{A}^2$, we thus obtain $\mathcal{S}^2$, and a bit of thought shows that $\mathcal{S}^2 \subset \mathcal{S}^1$. Continuing inductively produces a chain of sets
$$
\mathcal{S}^N \subset \mathcal{S}^{N-1} \subset \cdots \subset \mathcal{S}^1 \subset Q
$$
whose limiting set $\mathcal{S}^{\infty} := \cap_{N=1}^{\infty} \mathcal{S}^N$
is the attractor of the iterated function system with similitudes
$$
S_m (x) : = L^{-1} x + \frac{\alpha_m}{L} \textrm{ for some } \alpha_m \in \mathcal{A}.
$$
In particular, from standard self-similar set theory, we obtain the following structural result for the rational product Cantor sets $\mathcal{S}^{\infty}$.

\medskip

\begin{proposition}
Let $L \geq 3$ and $A_1,...,A_d \subset \{0,1,...,L-1\}$ be given, with associated rational product digit set $\mathcal{A}$ and rational product Cantor set $\mathcal{S}^{\infty} \subset [0,1]^d \subset \mathbb{R}^d$.
\begin{enumerate}[(1)]
    \item The set $\mathcal{S}^{\infty}$ is a self-similar set which satisfies the open set condition.
    \smallskip
    \item As a consequence of the open set condition, one has that
    $$
    \dim \mathcal{S}^{\infty} = \dim_S \mathcal{S}^{\infty} = \frac{\log\big[(\# A_1)\cdots (\# A_d)\big]}{\log L}.
    $$
    \smallskip
    \item Moreover, if one assumes that $\# \mathcal{A} = L$ and there exist indices $1 \leq i \neq j \leq d$ such that $\#A_i, \#A_j \geq 2$, then $\mathcal{S}^{\infty}$ is a purely $1$-unrectifiable set.
\end{enumerate}
\end{proposition}

We remark that all of the information encoded in the rational product Cantor set $\mathcal{S}^{\infty}$ is essentially obtainable from the information encoded in the rational product digit set $\mathcal{A}$, which, itself, carries the same amount of information as the sets $A_1,...,A_d \subset \{0,1,...,L-1\}$. The only meaningful distinction between the two is that $\mathcal{S}^{\infty}$ carries information at multiple scales $\delta_N : = L^{-N}$, whereas $\mathcal{A}$ is a strictly discrete (i.e. single scale) object. This phenomenon of a discrete structure iterated across countably-infinite scales enables us to leverage ideas from combinatorics together with ideas from multiscale analysis in order to derive power laws for the Favard length problem. We will accomplish this by considering the structure of the sets $A_1,...,A_d$, which is the topic of the next two subsections.

\begin{remark}
Although the Cantor iterate $\mathcal{S}^N$ is not literally the neighbourhood $\mathcal{N}_{L^{-N}}(\mathcal{S}^\infty)$, each can be covered by a bounded number of translates of the other, with the bound depending only on $d$.  Consequently,
\begin{equation}\label{eq:squaresandballs}
\Fav\big(\mathcal{N}_{L^{-N}}(\mathcal{S}^\infty)\big)\approx_d \Fav(\mathcal{S}^N).
\end{equation}
Since the implicit constant in our power laws is allowed to depend upon the ambient dimension, this does no harm, and we therefore pass freely between these two formulations.
\end{remark}

Using our method, we are able to prove the following concrete exponent improvement over the original $1/6 - u$ exponent of Nazarov, Peres, and Volberg, for a certain special choice of digit set $\mathcal{A} : = A_1 \times A_2$.

\begin{theorem}\label{thm:explicitonethird}
Let $\mathcal S^\infty : = \mathcal{S}^{\infty} (25 ; A_1 \times A_2) \subset \mathbb{R}^2$ be the planar product Cantor set in base $25$ with digit sets
$$
A_1=A_2=\{0,1,2,6,9\}.
$$
Then, for every $u > 0$, there exists an $N_u \in \mathbb{N}$ such that
\[
 \Fav\bigl(\mathcal N_{25^{-N}}(\mathcal S^\infty)\bigr)
 \lesssim_u N^{-1/3+u}, \quad \forall N \geq N_u.
\]
\end{theorem}

The special structure of the digit set $A : = \{0,1,2,6,9\}$ which produces the improvement in the exponent is introduced in the next subsection, and further examined in Section \ref{sec:exponentexamples}. Nazarov, Peres, and Volberg observed that a bound for the four-corner Cantor set decaying faster than $O(N^{-1/4})$ would require ``new ideas'' \cite[p. 3]{NPV}.  The results here go concretely beyond that scale in the broader class of rational product Cantor sets: the base-$25$ example reaches the exponent $1/3-u$, while the base-$36$ example in Theorem \ref{thm:explicitbeyondthird} exceeds the one-third threshold.  These examples concern different digit sets and do not assert the corresponding improvement for the four-corner Cantor set.

\subsection{The cyclotomic structure of rational product Cantor sets}\label{subsec:digitsetpoly}

The main method of proof for establishing upper bounds for asymptotics of the Favard length of rational product Cantor sets is estimating (both from above and below) integrals of the form
\begin{equation}\label{eq:rieszproductfirsttime}
\int_{L^{-N}}^1 \prod_{j=0}^{N-1} \vert \phi_{(t_1,...,t_{d-1})} (L^j \xi) \vert^2 d \xi,
\end{equation}
where $(t_1,...,t_{d-1}) \in [0,1]^{d-1}$ and 
$$
\phi_{(t_1,...,t_{d-1})} (\xi) : = \frac{1}{L} \sum\limits_{(a_1,...,a_d) \in A_1 \times \cdots \times A_d} e^{2 \pi i \xi (a_1 + t_1 a_2 + \cdots + t_{d-1} a_d)}.
$$
Such estimates are the main focus of the upcoming Sections \ref{sec:toolbox} and \ref{sec:rieszproductbounds} of this article. The integrand in \eqref{eq:rieszproductfirsttime} is an example of a \textit{lacunary Riesz product}---a multiscale trigonometric polynomial written as linear combinations of multiscale frequencies that are integral powers of some initial frequency scale $L > 1$. These have a long history within Fourier analysis and one can read \cite[Chapter 13]{mattila2} or \cite[Chapter 5]{Zygmund} for more background on the general theory.

The product structure of the set $\mathcal{A} : = A_1 \times \cdots \times A_d$ allows us to write
\begin{equation}\label{eq:phit}
\phi_t (\xi) = \phi_{A_1} (\xi) \phi_{A_2} (t_1 \xi) \cdots \phi_{A_d} (t_{d-1} \xi),
\end{equation}
where, for each $i=1,...,d$, we have
\begin{equation}\label{eq:phifirsttime}
 \phi_{A_i} (\xi) : = \frac{1}{\# A_i} A_i (e^{2\pi i \xi}) , \textrm{ where } A_i (X) : = \sum_{a \in A_i} X^a,
\end{equation}
and we remind the reader that, by the dimensional assumption, we have
$$
L : = (\#A_1) \cdots (\#A_d).
$$

These polynomials $A_1 (X),...,A_d (X)$ are so ubiquitous in what follows that we define them formally below.

\medskip

\begin{definition}
Given a finite and non-empty set $A \subset \mathbb{N}_0$, we let
$$
A(X) : = \sum\limits_{a \in A} X^a \in \mathbb{Z}[X]
$$
denote the \textbf{mask polynomial} of $A$.
\end{definition}

\medskip

For our purpose, we are quite concerned when the function $\phi_{A_i} : [0,1] \rightarrow \mathbb{C}$ vanishes. To understand when this might happen, we need the following definition.

\begin{definition}
Given $N \in \mathbb{N}$, the \textbf{$\bm{N}$-th cyclotomic polynomial} $\Phi_N (X) \in \mathbb{Z}[X]$ is the (unique) monic and irreducible polynomial satisfying
$$
\Phi_N \big(e^{\frac{2 \pi i d}{N}} \big) = 0,
$$
for all integers $d$ such that $\textrm{gcd}(d,N) = 1$.
\end{definition}

We record only the portion of this theory required for the Favard-length argument, beginning with the factorization introduced in \cite{BLV}.

\medskip

\begin{definition}\label{def:cycfactorization}
Let $A(X)$ be a mask polynomial, and let $v_{s}(A)$ denote the multiplicity of the cyclotomic polynomial $\Phi_s$ in $A$.  We define
\[
A^{(1)}(X):=\prod_{\substack{s\geq1\\ \gcd(s,\#A)\neq1}}
\Phi_s(X)^{v_s(A)},
\qquad
A^{(2)}(X):=\prod_{\substack{s\geq1\\ \gcd(s,\#A)=1}}
\Phi_s(X)^{v_s(A)},
\]
where factors with zero multiplicity are omitted.  Let $A^{(3)}$ be the product, with multiplicity, of the remaining irreducible factors having a zero on the unit circle that is not a root of unity, and let $A^{(4)}$ be the remaining factor, which has no zero on the unit circle.  Thus
\[
A=A^{(1)}A^{(2)}A^{(3)}A^{(4)}.
\]
We put
\[
A':=A^{(1)}A^{(3)}A^{(4)},
\qquad
A'':=A^{(2)},
\]
and refer to this polynomial  $A'' : = A'' (X)$ as the \textbf{accumulating factor} of $A$.  The set
\[
S_A^{(2)}:=\{s:v_s(A)>0,\ \gcd(s,\#A)=1\}
\]
records the indices of its distinct accumulating cyclotomic divisors; the factor $A''$ retains their multiplicities.
For symmetry of notation, we also write
\[
S_A^{(1)}:=\{s:v_s(A)>0,\ \gcd(s,\#A)\neq1\}.
\]
\end{definition}

The four factors and their zeroes have different analytic roles when considering lower bounds for our Riesz products. The factors in $\Phi_s \mid A^{(1)}$ are cyclotomic polynomials which do not accumulate their zeroes infinitely under base-$L$ dilation, since $\gcd (s,L) > 1$. These factors will be handled by our \textit{Set of Small Values (SSV)} estimate, which is the main subject of Section \ref{sec:ssvtechnical}.  The cyclotomic factors $\Phi_s \mid A^{(2)}$, on the other hand, do accumulate their zeroes infinitely under $L$-adic dilations, and so their zeroes must be accounted for at all scales. This is handled by the \textit{Set of Large Values (SLV)} construction, which is the main subject of Section \ref{sec:multiSLV}.  The factor $A^{(3)}$ contains the unit-circle zeros which are not roots of unity; it is responsible for the weaker, so-called \textit{logarithmic Set of Small Values (log-SSV)}  branch, which was handled by Bond, {\L}aba and Volberg in \cite{BLV}.  Finally, $A^{(4)}$ is uniformly separated from zero on the unit circle and creates only a fixed pointwise loss in the lower bound. The purpose of the factorization is therefore not merely algebraic: it apportions and assigns each possible unit-circle behaviour to exactly the necessary analytic device capable of controlling it (i.e. the \textit{SSV property}, \textit{SLV property}, \textit{log-SSV property}, and pointwise control).

We shall show in Section \ref{sec:exponentexamples} that when certain factors are simplified, the corresponding final exponent can be concretely improved. For example, the digit set $A_{\text{corner}} : = \{0,3\}$, with associated mask polynomial $A_{\text{corner}}(X) : = 1 + X^3$, satisfies
$$
A_{\text{corner}} (X) : = (1 + X)(1 -X + X^2) = \Phi_2 (X) \Phi_6 (X),
$$
and hence $A_{\text{corner}} (X) : = A_{\text{corner}}^{(1)} (X)$. Hence, the four-corner Cantor set (which is generated by the digit set $\mathcal{A}_{\text{corner}} : = A_{\text{corner}} \times A_{\text{corner}}$) has an implicitly special cyclotomic structure. As another example, we verify that the mask polynomial associated to the digit set $A_{\text{zerofree}} : = \{0,1,2,6,9\}$ has no zeroes on the unit circle, and hence satisfies $A_{\text{zerofree}} (X) = A_{\text{zerofree}}^{(4)} (X)$, and so the digit set $\mathcal{A}_{\text{zerofree}} = A_{\text{zerofree}} \times A_{\text{zerofree}}$ also has a special (but distinct from $A_{\text{corner}}$!) cyclotomic structure. These special cyclotomic factorizations of the two mask polynomials $A_{\text{corner}}(X)$ and $A_{\text{zerofree}}(X)$ are \textit{precisely} the number theoretic inputs which allow us to obtain the concrete exponent values of $1/5 - u$ and $1/3 - u$, as opposed to a more nebulous value of exponent as in Theorem \ref{thm:mainsimple}. We discuss and state our most general result, under the hypothesis that no factors $A_i^{(j)}$ vanish, in the next subsection.

In the general setting where no factor $A_i^{(j)}$ is identically one, we are most concerned with analyzing each type of zero of the associated digit sets $A_1,...,A_d$ via a distinct technique. Moreover, there is additional complexity in combining the zero information from each polynomial $A_1(X),...,A_d(X)$ in order to produce a joint lower bound on their associated Riesz product. In some sense, then, understanding and utilizing the zero patterns of the cyclotomically-factored mask polynomials $A_i^{(j)}$ is the main number theoretic challenge which must be resolved to obtain our results.

The following is our first main power law for the Favard length problem, which we obtain even when all four factors $A_i^{(1)},...,A_i^{(4)}$ are non-trivial.

\begin{theorem}\label{thm:fav2prime}\label{thm:favardsmall}
    Let $L\geq4$, let $A_1,\ldots,A_d\subset\{0,1,\ldots,L-1\}$, and assume
    \[
    (\#A_1)\cdots(\#A_d)=L,
    \qquad
    \#A_i,\#A_j\geq2\quad\text{for some }i\neq j.
    \]
    Let $\mathcal S^\infty$ be the associated rational product Cantor set and define
    $
    s_{A_i} := \textrm{lcm} \big(S_{A_i}^{(2)}\big) \in \mathbb{N}.
    $
    If each one of the $A_i$ satisfies at least one of the following two conditions:
    \begin{enumerate}
    \item $\# A_i \leq 10$; or else,
    \smallskip
    \item $s_{A_i}$ has (at most) two distinct prime factors,
    \end{enumerate}
    then there exists an $\epsilon \in (0,1)$, which may depend upon $A_1,...,A_d$ and $d \geq 2$, such that, for all sufficiently large $N$,
    \begin{equation}\label{eq:fav2prime}
    \Fav (\mathcal{N}_{L^{-N}} (\mathcal{S}^{\infty})) \lesssim N^{-\upsilon(\epsilon,N)}
    \end{equation}
    where
    $$
    \upsilon (\epsilon,N) : =
\begin{cases}
    \epsilon, \textrm{ if } A_i^{(3)} \equiv 1 \textrm{ for each } i = 1,...,d \\[1ex]
    \frac{\epsilon}{\log \log N}, \textrm{ otherwise }
\end{cases}
    $$
\end{theorem}

The special case of $\max (\#A_1, \# A_2) \leq 6$ was proven in the plane by M. Bond, I. {\L}aba and A. Volberg in \cite{BLV}. {\L}aba and the author then improved this result in the plane to Theorem \ref{thm:favardsmall}(1) by examining the explicit classifications of vanishing sums of roots of unity given in \cite{CDK} and \cite{PR}. The hypothesis of Theorem \ref{thm:fav2prime}.(2) may appear incongruous when compared to the required size-bound of Theorem \ref{thm:fav2prime}.(1). However, the two-prime-factors assumption essentially forces the set $S_A^{(2)}$ associated to the accumulating factor $A_i'' (X)$ from Definition \ref{def:cycfactorization} to be less complicated (this is made precise in \cite[Section 5]{LM}). With this simplifying assumption on $S_A^{(2)}$, {\L}aba and the author demonstrate a size bound on the set $A_i$, which acts as a substitute for the assumption $\# A_i \leq 10$. One can see \cite[Theorem 1.2]{LM} for the precise statement of this size bound.

Our flagship Theorem \ref{thm:mainsimple} then follows directly from Theorem \ref{thm:fav2prime}, which we show in the following example.

\begin{example}[A Higher-Dimensional Four-Corner Set]\label{examp:highfourcorner}
    \rm{Let $d \geq 2$ and $\mathcal{K}_d^{\infty}$ be the $2^d$-corner Cantor set defined in equations \eqref{eq:2diterates} and \eqref{eq:2dcantorset}. Since $\#A_i = 2 \leq 10$, Theorem \ref{thm:fav2prime} already furnishes some $\epsilon \in (0,1)$ such that
    $$
    \Fav \big(\mathcal{N}_{2^{-dN}} (\mathcal{K}_d^{\infty})\big) \lesssim_d N^{-\frac{\epsilon}{\log \log N}}.
    $$
    However, we can improve upon this by noticing that the mask polynomials
    $$
    A_i (X) : = 1 + X^{2^d -1}
    $$
    together satisfy $A_i(\zeta) = 0$ if and only if
    $$
    \zeta = e^{\frac{i \pi}{2^d-1} + \frac{2 \pi ir}{2^d-1}}
    $$
    for some $r \in \mathbb{Z}$. In particular, each $A_i$ has only roots of unity on the complex unit circle (equivalently, $A_i^{(3)} (X) \equiv 1$ for all $i = 1,...,d$), and so Theorem \ref{thm:fav2prime} now furnishes an $\epsilon \in (0,1)$ such that
    $$
    \Fav \big(\mathcal{N}_{2^{-dN}} (\mathcal{K}_d^{\infty})\big) \lesssim_d N^{-\epsilon_d},
    $$
    as required.
    }
\end{example}

\subsection{Favard length estimates and fibered digit sets}\label{subsec:fiberedfav}
In this subsection, we state the extension of Theorem \ref{thm:fav2prime} for a special class of digit sets $A_1,...,A_d$ which are related to the so-called \textit{fibered sets} studied by G. Kiss, I. {\L}aba, G. Somlai and the author in \cite{KLMS}. The main advancement is that we are able to drop the assumption that $s_{A_i} : = lcm (S_{A_i})$ has at most two prime factors if, instead, we assume that our digit sets $A_i$ admit some subset $A' \subseteq A$ which is highly structured modulo $s$ for each $s \in S_{A_i}^{(2)}$. We introduce the basic notations and definitions needed to state our result, and remark that a thorough treatment is given in \cite{KLMS}.

Letting $A\subset \mathbb{N}_{0}$ be fixed, we set $S = S_A^{(2)}$ and $M=\lcm(S)$. Because $s|M$ we must have that $\Phi_s (X) \mid (X^M-1)$, and so 

$$
\Phi_s (X) \mid A (X)
\Leftrightarrow
\Phi_s (X)\mid A(X)\bmod X^{M} - 1.
$$
One can readily check (or see from \cite[Section 2]{KLMS}) that
\begin{equation}\label{eq:inducedpoly}
A(X) \equiv \sum\limits_{a \in \mathbb{Z}_M} w_A^M (a) X^a \mod X^M - 1
\end{equation}
where
\begin{equation}\label{eq:weightfunction}
w_A^M (a) : = \# \{a' \in A : a' \equiv a \bmod M \}.
\end{equation}
Hence we see that divisibility by $\Phi_s (X)$ in the polynomial ring $\mathbb{Z}[X]$ is equivalent to divisibility of the induced polynomial \eqref{eq:inducedpoly} in the polynomial quotient ring $\mathbb{Z}[X]/(X^M-1)$. This leads us to consider the following more general paradigm.

\begin{definition}
A \textbf{weight} is any non-trivial function $w : \mathbb{Z}_M \rightarrow \mathbb{Z}$. Given a weight, its associated \textbf{multiset} on $\mathbb{Z}_M$ is the collection of ordered pairs
$
B : = \{(y,w (y)) \in \mathbb{Z}_M \times \mathbb{Z}\}.
$
We let $\mathcal{M} (\mathbb{Z}_M)$ denote the collection of all multisets on $\mathbb{Z}_M$. We also let $\mathcal{M}^+ (\mathbb{Z}_M)$ denote the collection of all multisets whose weights satisfy $w \geq 0$ uniformly on $\mathbb{Z}_M$. The \textbf{generating polynomial} $B \bmod M$ is the unique polynomial $B(X) \in \mathbb{Z}[X]/(X^M-1)$ which satisfies
$$
B(X) \equiv \sum\limits_{y \in \mathbb{Z}_M} w(y) X^y \mod X^M - 1.
$$
\end{definition}

If $A \subset \mathbb{N}_{0}$ is a non-empty set and $M \geq 2$ is any integer, then our previous discussion shows that
$$
A \bmod M : = \{(y, w_A^M (y)) : y \in \mathbb{Z}_M \textrm{ and } w_A^M (y) \textrm{ is as in } \eqref{eq:weightfunction}\}.
$$
The associated generating polynomial of $A \bmod M$ then has equation
$$
A (X) \equiv \sum\limits_{y \in \mathbb{Z}_M} w_A^M (y) X^y \mod X^M - 1,
$$
and, if $A \subset \mathbb{N}_0$, then $(A \bmod M) \in \mathcal{M}^+ (\mathbb{Z}_M)$ for any $M \geq 2$. We also always have that
$$
\# A = A(1) = \sum\limits_{y \in \mathbb{Z}_M} w_A^M (y).
$$
One can see \cite[Section 4]{LM} or \cite[Section 2]{KLMS} for a more thorough introduction to multisets and their generating polynomials.

We work with sets $A \subset \mathbb{N}_0$ whose induced multisets $A \bmod s$ have a special structural property called \textit{fibering}, which was first studied in \cite{KLMS}. To this end, let 
$$
M : =\lcm (S) = \prod_{k=1}^{K} p_k^{n_k},
$$
where $p_1,\dots, p_K$ are distinct primes and $n_1,\dots,n_K \in\mathbb{N}$. Notice that, since we are assuming that $S_A = S_A^{(2)}$, the hypothesis of Theorem \ref{thm:fav2prime}.(2) assumes that $K \leq 2$. However, we now make no such assumption moving forward. 

If $s \mid M$ and $p_k \mid s$, we say that $F \in \mathcal{M}^+ (\mathbb{Z}_M)$ is an \textit{$s$-fiber in the $p_k$ direction} if the multiset $F \bmod s$ satisfies
$$
F (X) \equiv \underbrace{1 + X^{s/p_k} + X^{2s/p_k} + \cdots + X^{(p_k - 1)s/p_k}}_{F_k^s (X)} \mod X^s - 1.
$$
We will always use the notation $F_k^s \in \mathcal{M}^+ (\mathbb{Z}_M)$ to denote an $s$-fiber in the $p_k$ direction. When $s$ and $p_k$ are assumed from context, we simply call such sets \textit{fibers}.

While we will not use the following test for cyclotomic divisibility in this work, we state it to explain how fibers are related to cyclotomic divisibility. The following is known as the \textit{de Bruijn--R\'edei--Schoenberg theorem} on the structure of vanishing sums of roots of unity (see \cite{deB,LL,Mann,Re1,Re2,schoen}).

\begin{theorem}\label{thm:BRSvanishing}
Let $A \subset \mathbb{N}_{0}$ have mask polynomial $A(X) \in \mathbb{Z}[X]$ and suppose $s = \prod_{k=1}^K p_k^{\alpha_k}$ for some distinct prime numbers $p_1,...,p_K$ and exponents $\alpha_k \geq 1$. Then, $\Phi_s (X) \mid A(X)$ if and only if there exist coefficient polynomials $P_1(X),...,P_K (X) \in \mathbb{Z}[X]$ such that
\begin{equation}\label{eq:bruijnredeischoen}
    A(X) \equiv \sum\limits_{k=1}^{K} P_k (X) F_k^s (X) \mod X^s - 1,
\end{equation}
where, we remark, the coefficient polynomials $P_k (X)$ have integer, but not necessarily non-negative, coefficients.
\end{theorem}

We want to work with sets $A \subset \mathbb{N}_0$ which admit a subset $A' \subseteq A$ whose induced multisets $A' \bmod s$ all have a certain rigid structure for every $s \in S_A^{(2)}$. This structure is known as \textit{fibering} and was first introduced in \cite[Section 3]{KLMS} relative to the study of the Coven-Meyerowitz conjecture for integer tilings.

\begin{definition}\label{def:fiberedandassignmentfunct}
Let $A \subset \mathbb{N}_{0}$ and
$
S \subseteq S_A : = \{s \in \mathbb{N} : \Phi_s (X) \mid A(X) \}
$
be some set of cyclotomic divisors of $A$. Set $M = \lcm (S) = \prod_{k=1}^{K} p_k^{n_k}$. Fix some $ s \mid M$ (and $s \neq 1$) and suppose that $p_k \mid s$. We say that $A$ is \textbf{fibered in the $p_k$-direction at scale $s$} if there exists a polynomial $P_k (X) \in \mathbb{N}_0[X]$ such that
$$
A(X) \equiv P_k (X) F_k^s (X) \mod X^s - 1.
$$
Equivalently, the induced nonnegative multiset $A\bmod s$ is a nonnegative integer combination of translates of one $p_k$-fiber.  The nonnegativity here is part of the definition and is essential for the cardinality theorem below; it is stronger than the signed decomposition furnished by Theorem \ref{thm:BRSvanishing}.
An \textbf{assignment function} is any function $\sigma : S \rightarrow \{1,...,K\}$ such that $$\sigma (s) \in \{i : p_i \mid s \}.$$
Given an assignment function $\sigma$, we say that $A$ is $(S , \sigma)$-\textbf{fibered} if, for every $s \in S$, the associated multiset $A$ mod $s$ is fibered in the $p_{\sigma(s)}$ direction on the scale $s$. 
\end{definition}

We also introduce the following definition (which is new to the literature, but is based upon the work of G. Kiss, I. {\L}aba, G. Somlai and the author in \cite{KLMS}). This definition allows us to work with sets which may not, themselves, be fibered, but which at least contain a fibered subset along those scales $S_A^{(2)} \subset S_A$ which are associated to the accumulating zeros of $A(X)$.

\begin{definition}\label{def:persistentfiber}
Let $A \subset \mathbb{N}_0$.
We say that $A$ \textbf{admits a fibered subset} if there exists some $\emptyset \neq A' \subseteq A$ such that $A'$ is an $(S_A^{(2)}, \sigma)$-fibered set.
\end{definition}

We have the following new Favard length estimate for rational product Cantor sets whose digit sets $A \subset \mathbb{N}_0$ admit a fibered set. Note that this result is entirely new to the literature, even in the well-studied case of planar rational product Cantor sets.

\begin{theorem}\label{thm:fibered}
    Let $L\geq4$, let $A_1,\ldots,A_d\subset\{0,1,\ldots,L-1\}$, and assume
    \[
    (\#A_1)\cdots(\#A_d)=L,
    \qquad
    \#A_i,\#A_j\geq2\quad\text{for some }i\neq j.
    \]
    Let $\mathcal S^\infty$ be the associated rational product Cantor set and suppose that each $A_i$ admits a fibered subset (see Definitions \ref{def:fiberedandassignmentfunct} and \ref{def:persistentfiber}).
    Then there exists an $\epsilon$ (which may depend upon $A_1,...,A_d$ and $d \geq 2$) such that, for all sufficiently large $N$,
    $$
    \Fav (\mathcal{N}_{L^{-N}} (\mathcal{S}^{\infty})) \lesssim N^{-\upsilon(\epsilon,N)}
    $$
    where
    $$
    \upsilon (\epsilon,N) : =
\begin{cases}
    \epsilon, \textrm{ if } A_i^{(3)} \equiv 1 \textrm{ for each } i = 1,...,d \\[1ex]
    \frac{\epsilon}{\log \log N}, \textrm{ otherwise }
\end{cases}
    \, \, .$$
\end{theorem}

To understand why one might expect to have a Favard length estimate for such digit sets, consider the following lower bound on the cardinality of fibered multisets, which was recently proven by G. Kiss, I. {\L}aba, G. Somlai and the author in \cite{KLMS}.

\begin{proposition}\label{prop:multifibered}
Let $A \subset \mathbb{N}_{0}$, let $S$ be a finite set of integers greater
than one, set $M:=\lcm(S)=\prod_{k=1}^Kp_k^{n_k}$, and suppose that $A$ is
$(S,\sigma)$-fibered for an assignment function $\sigma$.  For each
$k = 1,...,K$, let
$$
\capexp_k(S, \sigma) : = \{ \alpha\in \mathbb{N} : \exists \, s \in S \textrm{ with } (s, p_k^{n_k}) = p_k^{\alpha} \textrm{ and } \sigma (s) = k \}.
$$
Let $E_k(S,\sigma) : = \# \capexp_k(S, \sigma )$ and define
$$
\mathsf{FIB}(S,\sigma) : = p_1^{E_1(S,\sigma)} \cdots p_K^{E_K(S,\sigma)}.
$$
Then $\mathsf{FIB}(S,\sigma)$ divides $\#A$.  In particular,
\begin{equation}\label{eq:fiberedlower bound}
\# A \geq \mathsf{FIB}(S,\sigma)
\geq \min_{\tau}\mathsf{FIB}(S,\tau),
\end{equation}
where the minimum is taken over all assignment functions $\tau:S\rightarrow\{1,\ldots,K\}$.
\end{proposition}

\begin{proof}
This is \cite[Proposition 3.2]{KLMS}, applied to the nonnegative induced
multiset $A\bmod M$.  The positivity in Definition
\ref{def:fiberedandassignmentfunct} is what permits that proposition to be
applied.  When $S=\varnothing$, we use the convention
$\mathsf{FIB}(\varnothing,\varnothing)=1$.
\end{proof}

Of course, we immediately obtain the following.

\begin{corollary}\label{cor:multifibered}
Suppose that $A \subset \mathbb{N}_0$ admits a fibered subset. Then,
$$
\# A \geq \min_{\sigma} \mathsf{FIB}(S_A^{(2)}, \sigma),
$$
where the minimum is taken over all possible assignment functions $\sigma$ on $S_A^{(2)}$.
\end{corollary}

Critically: the definition of \textit{admitting a fibered subset} requires that the subset $A'$ is $(S,\sigma)$-fibered with $S = S_A^{(2)}$,  thus allowing us to obtain a non-trivial lower bound on $\#A$ which depends upon the accumulating roots of $A$. This cardinality estimate is the key starting point for our constructions in Section \ref{subsec:SLVfiberconst}.

\begin{remark}
    \rm{
    It is reasonable to ask why one would even think to study digit sets which have this $(S,\sigma)$-fibering property. To help explain how such a structure arises, observe that whenever $\Phi_s (X) \mid A(X)$, we have to have that $A \bmod s$ is a linear combination of fibers (this is just the Bruijn-R{\'e}dei-Schoenberg theorem, stated as Theorem \ref{thm:BRSvanishing} in this work). Hence, the $(S,\sigma)$-fibered sets are exactly those which have the \textit{simplest allowable structure relative to these fibers} at all $s \in S$. 
    
    In principle, this observation suggests that fibered sets are an ideal test class when dealing with constructions which are sensitive to complicated configurations of roots of unity. Indeed, it is exactly the simple structure of fibered sets which enable the size bound of Proposition \ref{prop:multifibered} to hold in full generality. For more complicated sets, such estimates can fail (which is a nontrivial fact, see \cite{KLMS} for a thorough discussion of this phenomenon).
    }
\end{remark}

\subsubsection{Examples of digit sets that admit fibered subsets}
To explain why Theorem \ref{thm:fibered} represents genuine progress on the Favard length problem, we give an example of a digit set $A \subset \mathbb{N}_0$ which fails to satisfy the hypothesis of Theorem \ref{thm:fav2prime} but which admits a fibered subset. Because of Theorem \ref{thm:fibered}, we can still obtain a power law for the associated RPC sets $\mathcal{S}^{\infty}$ which are associated to the Cantor iterates
$$
\mathcal{S}_d^{N} : = \underbrace{(A^N \times \cdots \times A^N)}_{d \textrm{ times}} + [0, (\# A)^{-dN}]^d.
$$
We remark that one can find similar styles of examples in \cite[Section 8]{KLMS} and \cite[Section 7]{LM}, where a more geometric interpretation of these constructions is developed. However, the example which we present is new to the literature.

The first example is intentionally simple, providing the reader with an intuition for the type of structured digit sets which Theorem \ref{thm:fibered} addresses. We work in $\mathbb{Z}_{30}$ because $30 = 2 \times 3 \times 5$ is the smallest possible integer for which the two-prime hypothesis of Theorem \ref{thm:fav2prime} fails.

\begin{example}
\rm{
Consider the set $A \subset \mathbb{N}_0$ with mask polynomial
$$
A(X) = (1 + X +X^5)F_3^{30} (X) + X^2F_2^{30} (X) \mod X^{30} - 1,
$$
or, if one prefers,
$$
A = (\{0,10,20\} \cup \{1,11,21\} \cup  \{5,15,25\}) \cup \{2,17\}.
$$
Then $\# A = 11$ and, since $A(X) \bmod X^{30} - 1$ is the displayed linear combination of fibers, Theorem \ref{thm:BRSvanishing} implies that $\Phi_{30} (X) \mid A (X)$. In fact, one can also verify by direct calculation that
$$
A(X) = (1 + 2X)F_3^{6} (X) + X^2 F_2^6 (X) \mod X^6 - 1
$$
so that $\Phi_{6} (X) \mid A(X)$.
Moreover, we know that $\Phi_{s} (X) \nmid A(X)$ for any $s \not \in \{6, 30\}$, via standard arguments regarding fibers (see, for example, either \cite[Section 2]{KLMS} or \cite[Section 5]{LM} for examples). In particular, then, $S_A^{(2)} = \{6,30\}$ so that $s_A = \lcm (S_A^{(2)}) = 30$ has three distinct prime divisors. Nevertheless, we see that $A$ admits at least one (and, in fact, several) fibered subset, for example $A' : =\{0,10,20 \}$.
}
\end{example}

Our previous examples show that Theorem \ref{thm:fibered} is new and tangible progress towards power laws for rational product Cantor sets in $\mathbb{R}^d$ associated to arbitrary digit sets $A_1,...,A_d \subset \mathbb{N}_0$. However, the key ingredient in our proof of Theorem \ref{thm:fibered} is the lower bound \eqref{eq:fiberedlower bound}, which \textit{fails to hold} for arbitrary sets of integers $A \subset \mathbb{N}_{0}$. This is actually a complicated (and rather surprising) property of multisets, which is proven by direct construction in \cite[Theorem 1.3]{KLMS}. Hence, any progress beyond that of Theorem \ref{thm:fibered} would require significant new ideas in the study of vanishing sums of roots of unity and cyclotomic divisibility. The author is currently investigating this, with the hope of addressing this in a sequel to this work.

\subsection*{Some conventions for polynomials}

We use the conventions $\operatorname{lcm}(\varnothing)=1$ and that an empty product of polynomial factors is equal to $1$.  Cyclotomic divisors are recorded with multiplicity in the factorization itself, while the sets $S_A^{(j)}$ record only their indices.  Singleton digit sets and the cases of zero, one, or two accumulating prime divisors are interpreted using these conventions.

\subsection{Summary of general Favard length estimates for rational product Cantor sets}\label{subsec:summaryresults}

The preceding results combine into the following statement.

\begin{theorem}\label{thm:totalFavbound}
Let $L\geq4$, let $A_1,\ldots,A_d\subset\{0,1,\ldots,L-1\}$, and suppose that
\[
(\#A_1)\cdots(\#A_d)=L
\]
and that $\#A_i,\#A_j\geq2$ for some $i\neq j$.  Let $\mathcal{S}^\infty$ be the associated rational product Cantor set.  Assume, for every $i=1,\ldots,d$, that at least one of the following holds:
\begin{enumerate}[(1)]
\item $\#A_i\leq10$;
\item the integer
\[
s_{A_i}:=\lcm\{s\in\mathbb{N}:\Phi_s(X)\mid A_i(X),\ \gcd(s,\#A_i)=1\}
\]
has at most two distinct prime divisors;
\item $A_i$ admits a fibered subset in the sense of Definition \ref{def:persistentfiber}.
\end{enumerate}
Then there exists $\epsilon\in(0,1)$, depending only on the digit sets and on $d$, such that
\[
\Fav\big(\mathcal{N}_{L^{-N}}(\mathcal{S}^\infty)\big)
\lesssim_{A_1,\ldots,A_d,d} N^{-\upsilon(\epsilon,N)},
\]
where
\[
\upsilon(\epsilon,N):=
\begin{cases}
\epsilon, & A_i^{(3)}\equiv1\text{ for every }i,\\[1ex]
\epsilon/(\log\log N), & \text{otherwise}.
\end{cases}
\]
\end{theorem}

The hypotheses may be mixed from coordinate to coordinate.  Their common role is to compare the cardinality of each digit set with the complexity of its accumulating cyclotomic zeros. The exact argument is presented in Section \ref{sec:multiSLV}. The remainder of the article focuses on proving all aforementioned Favard length estimates.

\section{The Counting Function and Associated Riesz Products}\label{sec:toolbox}

This section introduces the geometric and Fourier framework for proving our power laws for the Favard length problem.

\subsection{A covering of the sphere by finitely-many projective charts}
Our first objective is to restrict our estimate to a single representative projective chart, taken from a finite collection of projective charts which cover the sphere $\mathbb{S}^{d-1}$. As a consequence, our later Fourier analysis will occur in a compact ``safe'' interval of frequencies adapted to the corresponding scale $\delta_N : = L^{-N}$ of our associated rational product Cantor set.

For $1\leq j\leq d$ and a given sign vector
$\chi : = (\chi_1,...,\chi_d) \in\{-1,1\}^d$, let $U_{j,\chi} \subset \mathbb{S}^{d-1}$ be those directions $\theta : = (\theta_1,...,\theta_d)$ satisfying
$|\theta_j|=\max_i|\theta_i|$ and the coordinate signs are prescribed by $\chi$. That is to say
$$
U_{j, \chi} : = \big\{ \theta \in \mathbb{S}^{d-1} : \vert \theta_j \vert = \max_{1 \leq i \leq d}\vert \theta_i\vert, \, \sign (\theta_i) = \chi_i \big\}.
$$
To assign coordinates on the chart $U_{j,\chi}$, we set
$$
t_i:=\frac{|\theta_i|}{|\theta_j|}, \qquad 1 \leq i \leq d, \, i \neq j.
$$
Since $|\theta_j|$ is maximal, we have that $t_i \in [0,1]^{d-1}$, and since $\theta \in \mathbb{S}^{d-1}$, one also has that
$$
1=|\theta_j|^2\bigg(1+\sum_{i\neq j}t_i^2\bigg) \Longrightarrow 
|\theta_j|=\bigg(1+\sum_{i\neq j}t_i^2\bigg)^{-1/2}.
$$
Up to $\sigma^{d-1}$-null sets, these charts form a finite cover of $\mathbb S^{d-1}$. Moreover, on each such chart $U_{j,\chi}$, we have the following smooth parametrization
\begin{equation}\label{eq:changeofvariablestheta}
 \theta=\Theta_{j,\chi}(t)
 :=\frac{\chi_j e_j+\sum_{i\ne j}\chi_i t_i e_i}
 {\bigl(1+\sum_{i\ne j}t_i^2\bigr)^{1/2}},
 \qquad t=(t_i)_{i\ne j}\in[0,1]^{d-1},
\end{equation}
where $e_1,...,e_d$ denotes the standard orthonormal basis in $\mathbb{R}^d$. We then observe that, for any $\theta \in U_{j, \chi}$, that the corresponding orthogonal projection $\pi_{\theta} : \mathbb{R}^d \rightarrow \mathbb{R}$ satisfies
\begin{equation}\label{eq:dotproductidentity}
 \pi_\theta(x)=\theta_j\left(x_j+\sum_{i\ne j}\sigma_i t_i x_i\right),
 \qquad \sigma_i:=\chi_i/\chi_j,
 \qquad d^{-1/2}\leq |\theta_j|\leq1.
\end{equation}
The main point of the calculations in this section is to demonstrate that we can replace the standard orthogonal projection map $\pi_{\theta}$ with the \textit{flattened orthogonal projection map}, which is the term appearing in the large braces in \eqref{eq:dotproductidentity}.

To this end, we include the Jacobian calculation, since it explains why no scale-dependent loss is introduced by this change of variables. Let
$$
r(t):=\bigg(1+\sum_{i\neq j}t_i^2\bigg)^{1/2}.
$$
After suppressing the fixed signs and permuting the coordinates, the parametrization is simply $(1,t)/r(t)$. The Gram matrix of its tangent vectors is
$$
G_{k\ell}(t)=\frac{\delta_{k\ell}}{r(t)^2}-\frac{t_kt_\ell}{r(t)^4}.
$$
The matrix determinant lemma gives $\det G(t)=r(t)^{-2d}$, and therefore the surface element on this chart is
\begin{equation}\label{eq:projectivechartjacobian}
d\sigma^{d-1}(\theta)=c_d\bigg(1+\sum_{i\neq j}t_i^2\bigg)^{-d/2}dt,
\end{equation}
where $c_d$ is the fixed normalization constant for surface measure. Since $t\in[0,1]^{d-1}$, the density in \eqref{eq:projectivechartjacobian} lies between two positive constants depending only upon $d$. Thus, once we have restricted any calculation to some local chart, integration with respect to $d\sigma^{d-1}(\theta)$ is equivalent to integration with respect to Lebesgue measure $dt$.

After permuting the coordinates so that the $j$th coordinate comes first, we thus define a \textit{projective orthogonal projection map} as
$$
 p_{j,\sigma,t}(x):=x_j+\sum_{i\ne j}\sigma_i t_i x_i, \quad \forall x \in \mathbb{R}^d.
$$
and observe that, if $\theta = \Theta_{j, \chi} (t)$, then one has the identity
\begin{equation}\label{eq:thetaprojectionproptprojection}
\lvert \pi_{\theta} (x) \rvert = \lvert \theta_j \rvert \lvert  p_{j,\sigma,t} (x) \rvert \approx_d \lvert  p_{j,\sigma,t} (x) \rvert, \quad \forall x \in \mathbb{R}^d.
\end{equation}
This relationship \eqref{eq:thetaprojectionproptprojection} shows that our original orthogonal projection map, and our new parametrized projection map $p_{j, \sigma, t}$, are proportional on each chart $U_{j, \chi}$. Moreover, our previous Jacobian calculation shows that, at least up to dimensional constants, there is no harm in substituting integration against $\sigma^{d-1}(\theta)$ with integration against $\mathcal{H}^{d-1}(t)$. Hence, up to issues with sign-changes and permutations of the digits $x_1,...,x_d$ (which are handled in due course of the argument), we may freely pass between $\pi_{\theta}$ and $p_{j, \sigma, t}$, once the underlying chart $U_{j, \chi}$ has been selected.

It is useful to make the preceding proportionality precise at the level of projected measures. If $\mu$ is a finite measure on $\mathbb{R}^d$, write
$$
\nu_{\theta}:=(\pi_{\theta})_{\#}\mu,
\qquad
\nu_{j,\sigma,t}:=(p_{j,\sigma,t})_{\#}\mu.
$$
Then $\nu_{\theta}$ is the pushforward of $\nu_{j,\sigma,t}$ by the one-dimensional dilation $y\mapsto\theta_jy$. In particular,
\begin{equation}\label{eq:chartmeasure-dilation}
 \nu_\theta=(D_{\theta_j})_\#\nu_{j,\sigma,t},
 \qquad
 \widehat{\nu_\theta}(\xi)
 =\widehat{\nu_{j,\sigma,t}}(\theta_j\xi),
 \qquad D_a(y):=ay.
\end{equation}
Hence, for every nonnegative integrand $F$,
\begin{equation}\label{eq:chartfrequency-dilation}
 \int_I F(\theta_j\xi)\,d\xi
 =|\theta_j|^{-1}\int_{\theta_j I}F(u)\,du,
 \qquad
 1\leq|\theta_j|^{-1}\leq\sqrt d.
\end{equation}
Thus only a dimensional constant is lost.  This is one common dilation; the
individual parameters $t_i$ remain inside their coordinate factors.

Finally, if $\sigma_i=-1$, then the $i$th digit set is reflected. Replacing $A_i$ by $A_i^*:=\{L-1-a:a\in A_i\}$ gives
$$
\phi_{A_i^*}(u)=e^{2\pi i(L-1)u}\overline{\phi_{A_i}(u)},
$$
and consequently $|\phi_{A_i^*}(u)|=|\phi_{A_i}(u)|$. The fixed chart signs therefore have no effect on any later estimate involving the absolute value of a Riesz product.

\subsection{Cylinders and the symbolic tree associated to a Cantor set}\label{subsec:cylindertree}

We now introduce the symbolic language which simultaneously describes the construction of the Cantor iterates and the combinatorial structure used in Section \ref{sec:comblemmata}. The purpose of this subsection is to distinguish carefully between a \textit{symbolic cylinder}, its \textit{geometric realization}, and its \textit{one-dimensional projection}. These objects are closely related, but they play different roles in the argument.

Recall that our alphabet is the finite product digit set
$$
\mathcal{A}:=A_1\times\cdots\times A_d,
$$
and that the standing dimension-one assumption gives $\#\mathcal{A}=L$. For each $n\geq1$, a \textbf{word of depth $\bm{n}$} is an ordered sequence
$$
\omega=(\alpha_1,\ldots,\alpha_n)\in\mathcal{A}^n.
$$
We write $|\omega|=n$ for its depth and use $\varnothing$ to denote the empty word, which has depth zero. The collection
$$
\mathcal{A}^{<\infty}:=\{\varnothing\}\cup\bigcup_{n\geq1}\mathcal{A}^n
$$
thus forms a rooted tree, whose root is $\varnothing$, with the children of a given word $\omega$ consisting of the $L$-many words $\omega\alpha$, where we take $\alpha\in\mathcal A$.

If $\tau\in\mathcal A^k$ and $\omega\in\mathcal A^n$ with $k\leq n$, we say that $\tau$ is a \textbf{prefix} of $\omega$, and write $\tau\preceq\omega$, when the first $k$ symbols of $\omega$ agree with $\tau$. Equivalently, there is a word $\eta$ such that $\omega=\tau\eta$. In this case, $\tau$ is an ancestor of $\omega$, while $\omega$ is a descendant of $\tau$. Every nonempty word has a unique parent, obtained by deleting its final symbol.

There is also a symbolic cylinder associated to each finite word. Namely, if $\mathcal A^{\mathbb N}$ denotes the space of infinite digit sequences, we put
$$
[\omega]:=\{(\beta_1,\beta_2,\ldots)\in\mathcal A^{\mathbb N}:(\beta_1,\ldots,\beta_n)=\omega\}.
$$
Thus $[\omega]$ consists of all infinite branches of the tree which pass through the vertex $\omega$.

Each symbolic cylinder is directly related to some subset of our rational product Cantor set $\mathcal{S}^{\infty}$ by the following coding map
$$
\Pi_{\mathcal A}(\beta_1,\beta_2,\ldots)
:=\sum_{k=1}^{\infty}L^{-k}\beta_k.
$$
For a word $\omega=(\alpha_1,\ldots,\alpha_n)$, define
\begin{equation}\label{eq:cylindercube}
z_\omega:=\sum_{k=1}^n L^{-k}\alpha_k,
\qquad
\mathcal S_\omega^\infty:=z_\omega+L^{-n}\mathcal S^\infty,
\qquad
Q_\omega:=z_\omega+[0,L^{-n}]^d.
\end{equation}
The set $\mathcal S_\omega^\infty=\Pi_{\mathcal A}([\omega]) \subseteq \mathcal{S}^{\infty}$ is thus the portion of the limiting rational product Cantor set determined by the prefix $\omega \in \mathcal{A}^n$, while $Q_\omega \supset \mathcal{S}_{\omega}^{\infty}$ is its associated enclosing constructive cube. We refer to $Q_\omega$ as the \textbf{cylinder cube} associated to $\omega$. At the root, we use the convention $Q_\varnothing=[0,1]^d$.

The \textbf{depth} of any cylinder $Q_\omega$ is defined as
$$
\operatorname{depth}(Q_\omega):=\lvert (\omega_1,...,\omega_n)\rvert : = n.
$$
Thus, the depth is an integer recording the generation of a given word in the total rooted tree tree. Critically: this value should not be confused with either the side length of the cube or the length of its projected interval.

The prefix relation of words in our rooted tree is exactly reflected by geometric containment of sub-cubes in our Cantor iterates. If $\tau\preceq\omega$, then all digit choices occurring after $\tau$ enter at scales strictly below $L^{-|\tau|}$, and therefore
\begin{equation}\label{eq:prefixcontainment}
\tau\preceq\omega
\quad\Longrightarrow\quad
\mathcal S_\omega^\infty\subset\mathcal S_\tau^\infty
\quad\text{and}\quad
Q_\omega\subset Q_\tau.
\end{equation}
Indeed, writing $\omega=\tau\eta$, with $|\tau|=k$ and $|\omega|=n$, gives
\[
 z_\omega-z_\tau
 =\sum_{r=k+1}^{n}L^{-r}\alpha_r
 \in[0,L^{-k}-L^{-n}]^d.
\]
Consequently,
\[
 z_\omega+[0,L^{-n}]^d
 \subset z_\tau+[0,L^{-k}]^d,
\]
which proves the cube containment in \eqref{eq:prefixcontainment}; the Cantor-set
containment follows by applying the same decomposition to every continuation
of $\omega$.
Consequently, every cylinder cube has a unique symbolic parent. Moreover, if $m\geq n$, then every depth-$n$ cylinder has exactly $L^{m-n}$ descendants at depth $m$.

\begin{figure}[H]
  \centering
  \includegraphics[width=\linewidth]{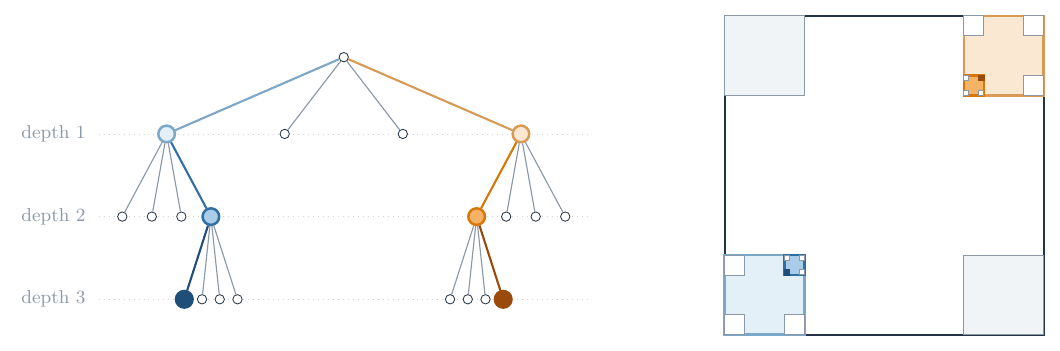}
  \caption{Two word branches and their corresponding nested squares in three
    generations of the rational product Cantor construction. The blue and
    orange chains distinguish the two branches, with each color darkening from
    ancestor to child to terminal child; the uncolored nodes and squares are
    the remaining sibling choices.}
  \label{fig:word-tree-cantor-correspondence}
\end{figure}

We further assign to each cylinder its \textbf{geometric length}
\begin{equation}\label{eq:symboliclength}
\ell(Q_\omega):=L^{-|\omega|}.
\end{equation}
This is both the side length of the cube, as well as its natural mass on the word tree; by which we mean, the mass of a cylinder is exactly the sum of the masses of its children:
$$
\sum_{\alpha\in\mathcal A}\ell(Q_{\omega\alpha})
=L\cdot L^{-(|\omega|+1)}
=\ell(Q_\omega).
$$
Iterating this identity, we obtain
\begin{equation}\label{eq:descendantmassconservation}
\sum_{\substack{R\in\mathscr Q_m\\R\subset Q_\omega}}\ell(R)
=\ell(Q_\omega),
\qquad m\geq|\omega|,
\end{equation}
where $\mathscr Q_m:=\{Q_\tau:\tau\in\mathcal A^m\}$.

We now arrive at our most critical definition for our word tree. A family of words, or of their corresponding cylinders, is called \textbf{prefix-free} if no member is a prefix of another. In the language of rooted trees, such a family is an antichain. Although the projected intervals of two prefix-free cylinders may overlap, the collections of their terminal descendants are disjoint as collections of words. This purely symbolic disjointness is the source of the following packing estimate, which we will use often in Section \ref{sec:comblemmata}.

\begin{lemma}[Prefix-free packing]\label{lma:prefixfreepacking}
Let $P$ be a cylinder cube and let $\mathscr F$ be a finite prefix-free family of descendants of $P$. Then
\begin{equation}\label{eq:prefixfreepacking}
\sum_{R\in\mathscr F}\ell(R)\leq\ell(P).
\end{equation}
\end{lemma}

\begin{proof}
Choose an integer $M$ which is at least the depth of every cylinder in $\mathscr F$. A cylinder $R\in\mathscr F$ has exactly $L^{M-\operatorname{depth}(R)}$ descendants at depth $M$. Because $\mathscr F$ is prefix-free, these collections of terminal descendants are pairwise disjoint. They are all contained among the $L^{M-\operatorname{depth}(P)}$ depth-$M$ descendants of $P$, and hence
$$
\sum_{R\in\mathscr F}L^{M-\operatorname{depth}(R)}
\leq L^{M-\operatorname{depth}(P)}.
$$
Dividing by $L^M$ gives \eqref{eq:prefixfreepacking}.
\end{proof}

The previous inequality is an example of a \textit{Kraft inequality} \cite{Kraft}; we include the short tree proof in the present notation because this is the form used repeatedly in Section \ref{sec:comblemmata}. In its essence, then, the Kraft inequality is a statement about the symbolic tree, rather than a geometric disjointness statement. Two prefix-free cylinder cubes may meet along their boundaries, and their one-dimensional projections may overlap substantially. Nevertheless, a terminal word cannot descend from two different members of a prefix-free family. It is precisely this distinction which allows the combinatorial argument to count descendants without assuming that their projections are disjoint.

As a general principle, we will often use the four-corner Cantor set as a model example. In the special case of the four-corner Cantor set, one has $L=4$ and $\mathcal A=\{0,3\}\times\{0,3\}$. The root then has four children, corresponding to the four first-generation corner squares, and every later vertex again has four children. A word of depth $n$ records the successive corners chosen during the first $n$ stages, while its prefixes record the chain of containing squares. This familiar example provides a concrete model for the much more general cylinder language used throughout the remainder of the paper.

\subsection{Projected measures and counting functions}
In this section, we define the \textit{counting function}, a key analytical tool for detecting when the orthogonal projections of our fixed rational product Cantor iterate $\mathcal{S}^N$ can be small. In our definitions, we will use the natural constructive cubes rather than enlarged Euclidean balls. This distinction is harmless for projection estimates, but important for the combinatorial arguments in Section \ref{sec:comblemmata}: constructive cubes form an exact rooted tree, whereas enlarged balls need not be nested.

For a word $\omega\in\mathcal A^n$, let
$$
I_{\omega,\theta}:=\pi_\theta(Q_\omega)
$$
be the corresponding projected cylinder interval; and, observe, the ambient scale parameter $n \geq 1$ is carried in the definition of $\omega : = (\omega_1,...,\omega_n)$. By our construction of rational product Cantor sets in Section \ref{subsec:cylindertree}, we have that
$$
\mathcal S^n=\bigcup_{\omega\in\mathcal A^n}Q_\omega.
$$
Moreover, because $Q_\omega$ is an axis-parallel cube of side length $L^{-n}$, we have
\begin{equation}\label{eq:projectedcylinderlength}
|I_{\omega,\theta}|=\|\theta\|_{\ell^1}L^{-n}\approx_d L^{-n} = L^{- \lvert \omega \rvert}
\end{equation}
Thus the geometric length of a projected cylinder is uniformly comparable to its geometric length $\ell(Q_\omega)=L^{-n}$, which is itself related to the standard length $\lvert \omega \rvert$ of its associated word $\omega \in \mathcal{A}^n$ . Replacing these intervals by concentric intervals whose lengths differ by a fixed dimensional factor changes our estimate by a harmless dimensional constant, as well as a fixed rescaling of multiplicity thresholds. We use that equivalent normalization later when taking Fourier transforms.

We can now define the main object of our analysis.

\begin{definition}
    Given $N \geq 1$ and $\theta \in \mathbb{S}^{d-1}$, the \textbf{counting function} is the mapping
    $$
    f_{N,\theta} (x) : = \sum\limits_{\omega \in \mathcal{A}^N} \mathbf{1}_{I_{\omega,\theta}} (x), \quad \forall x \in \mathbb{R},
    $$
    where $Q_\omega \subset \mathcal{S}^N$ is defined in \eqref{eq:cylindercube} and $I_{\omega,\theta} \subset \pi_{\theta} \big(\mathcal{S}^N\big)$ is its projected cylinder interval.
\end{definition}

For our calculations, it is much easier to work in parametrized projective coordinates. Hence, for notational simplicity, within a fixed chart $U_{j, \chi}$, we relabel the normalized coordinate as
the first coordinate and write its remaining parameters as $t_1,\ldots,t_{d-1}\in[0,1]$. The coordinate sign vector $\chi : = (\chi_1,...,\chi_d) \in \{1,-1\}^d$ then remains fixed throughout the argument. 

After reflecting the digit sets which correspond to negative chart signs, the associated discrete projected probability measure for the projective orthogonal projection map is given, for each scale $1 \leq n \leq N$, by the convolution identity
\begin{equation}\label{eq:bignu}
 \nu_{n,t}:=\mathop{\ast}_{k=1}^n\widetilde\nu_{k,t},
\end{equation}
where
\begin{equation}\label{eq:littlenu}
 \widetilde\nu_{k,t}
 :=L^{-1}\!\sum_{(a_1,\ldots,a_d)\in A_1\times\cdots\times A_d}
 \delta_{L^{-k}(a_1+t_1a_2+\cdots+t_{d-1}a_d)}.
\end{equation}
Each measure $\widetilde\nu_{k,t}$ records the contribution of the $k$th digit to the projection, and its total mass is one because $\#\mathcal A=L$. Their convolution $\nu_{n,t}$ therefore places mass $L^{-n}$ at the projected location associated to each word $\omega\in\mathcal A^n$, retaining the symbolic multiplicity when two different words have the same projected location.

We then set the following simplifying notation
\[
 \sigma_t:=1+t_1+\cdots+t_{d-1},
 \qquad
 b_{n,t}:=\frac{\sigma_t}{2}L^{-n},
 \qquad
 \nu_{n,t}^{\mathrm c}:=(D_{1,b_{n,t}})_\#\nu_{n,t},
\]
where $D_{1,b}(x)=x+b$.  After the fixed reflections used above, the
projective image of $Q_\omega=z_\omega+[0,L^{-n}]^d$ is the interval
\[
 p_t(z_\omega)+[0,\sigma_tL^{-n}]
 :=p_t(z_\omega)+b_{n,t}
  +[-b_{n,t},b_{n,t}].
\]
Since every atom of the measure $\nu_{n,t}^{\mathrm c}$ is assigned mass $L^{-n}$, the counting function therefore has the following exact convolution representation
\begin{equation}\label{eq:flattenedcountingfunction}
 f_{n,t}(x)
 :=L^n\bigl(\nu_{n,t}^{\mathrm c}
 *\mathbf 1_{[-b_{n,t},b_{n,t}]}\bigr)(x).
\end{equation}
Moreover,
\begin{equation}\label{eq:centeredmeasurephase}
 \widehat{\nu_{n,t}^{\mathrm c}}(\xi)
 =e^{-2\pi ib_{n,t}\xi}\widehat{\nu_{n,t}}(\xi),
 \qquad
 |\widehat{\nu_{n,t}^{\mathrm c}}(\xi)|
 =|\widehat{\nu_{n,t}}(\xi)|.
\end{equation}
Thus, and as we had previously claimed, re-centering the projected intervals introduces no loss in any Fourier
estimate.  As an added bonus, the bounds $1\leq\sigma_t\leq d$ permit the variable smoothing
interval in \eqref{eq:flattenedcountingfunction} to be replaced by fixed
dimensional inner or outer intervals whenever only one-sided estimates are
needed.

In the calculations which follow, we perform integration in $t$ chartwise, and we remind the reader that our previous calculation \eqref{eq:projectivechartjacobian} already gives the explicit numerical comparison
\begin{equation}\label{eq:chartmeasurecomparison}
 c_dd^{-d/2}|F|
 \leq \sigma^{d-1}(\Theta_{j,\chi}(F))
 \leq c_d|F|
\end{equation}
for every measurable $F\subset[0,1]^{d-1}$.  Hence, in our final Favard length upper bound, we need only sum over the (at most)
$d2^d$-many signed charts, and thus we modify only our implicit dimensional constants. The remainder of the calculation thus proceeds chart-by-chart.

\subsection{Exceptional directions and Fourier pigeonholing}\label{subsec:countingandlowmult}
We now use the cylinder tree from Section \ref{subsec:cylindertree} to compare the amount of projection overlap which occurs at different depths. To detect stacking at any depth up to $N$, we define the maximal counting function
$$
f_{N,\theta}^* (x) : = \max_{1 \leq n \leq N} f_{n,\theta} (x), \quad  x \in \mathbb{R},
$$
which encodes the largest stacking multiplicity observed along the first $N$ generations of the tree. Notice that the maximum is taken \textit{pointwise}: and thus, the depth at which the maximum occurs is allowed to depend upon $x$. This will become a key subtlety which must be handled in the combinatorial arguments in Section \ref{sec:comblemmata}.

The $L^2$ norm of a fixed-depth counting function has a useful geometric interpretation, as expanding the square and integrating provide that
\begin{equation}\label{eq:L2overlapidentity}
\|f_{n,\theta}\|_2^2
=\sum_{\omega,\tau\in\mathcal A^n}
\mathcal H^1(I_{\omega,\theta}\cap I_{\tau,\theta}).
\end{equation}
Thus, the quantity $\|f_{n,\theta}\|_2^2$ measures the total pairwise overlap of the projected depth-$n$ cylinders. A small $L^2$ norm indicates that the corresponding projection has relatively little \textbf{stacking} (i.e. a term we use to imply a failure of injectivity of the projection map of basic cubes $Q_{\omega}$ at fixed scale $n$ and direction $t \sim \theta$), while a large norm records many coincident or substantially overlapping projected intervals.

We use this maximal counting function to define a set of \textit{low multiplicity directions} $E_{N,K} \subset \mathbb{S}^{d-1}$. Intuitively, these are the directions for which the projection mapping rarely becomes $K$-to-one at any of the first $N$ depths. To make this notion rigorous, we define a \textit{stacking parameter} $K = K(N)$ which is a function of the maximal cutoff depth $N$, and whose value is determined by the cyclotomic factorization of the mask polynomials $A_1,...,A_d$.

\medskip

\begin{definition}\label{def:stackingparameter}
    Given some $\epsilon_0 \in (0,1)$ to be chosen later (but independent of $N$) define a parameter $K$, which depends upon $A_1,...,A_d$, $N$ and $\epsilon_0$ in the following manner:
    \begin{enumerate}
        \item If $A_i^{(3)} \equiv 1$ for each $i \in \{1,...,d\}$, we let $K = N^{\epsilon_0}$.
        \medskip
        \item Otherwise, we choose $K = N^{\epsilon_0/\log \log N}$.
    \end{enumerate}
\end{definition}

Using this choice of $K$, we now define the set of low multiplicity directions.

\medskip

\begin{definition}\label{ch3:def:lowmultiplicity}
Given $N \gg 1$ and its associated parameter $K$ (see Definition \ref{def:stackingparameter}), we let
$$
G = G_{N,K,\theta} := \{x \in \mathbb{R} : f_{N,\theta}^* (x) \geq K \}
$$
denote the level set of $f_{N,\theta}^*$ at height $K$. If $\rho \in (2, \infty)$ is some parameter to be chosen later, then we let
\begin{equation}\label{eq:lowmultfirsttime}
E_{N, K} : = \{\theta \in \mathbb{S}^{d-1} : \mathcal{H}^1 \big(G_{N,K,\theta}) \leq K^{-\rho} \}
\end{equation}
denote the \textbf{set of low multiplicity directions} for the counting function.
\end{definition}

In practice, we fix any $\rho>2$ and may subsequently let $\rho\downarrow2$.
The endpoint is not asserted (and, in fact, the author believes the current $L^2$ estimate to be sharp).  This is discussed further in Section \ref{sec:comblemmata}.

The terminology \textit{low multiplicity} refers to the size of the level set rather than to a pointwise upper bound. A direction may belong to $E_{N,K}$ even though $f_{N,\theta}^*$ exceeds $K$ at some points; membership means that the set of such points has length at most $K^{-\rho}$.

The following lemma makes precise the notion that the set of low multiplicity directions is where we do \textit{not} observe a large failure of injectivity in the stacking function $f_{N,\theta}^*$.

\medskip

In the following statement, we fix some geometric constant $C_{\mathrm{prop}}\geq1$, whose value is determined by the proof occurring in Section \ref{subsec:reverseholder}, and define an associated \textbf{propagation parameter} as
\begin{equation}\label{eq:propagationblocks}
J_{K,\rho}:=\left\lceil C_{\mathrm{prop}}K^{\rho-1}\log(2+K)\right\rceil.
\end{equation}

\begin{lemma}\label{lma:reverseholder}
If $N \gg 1$ and $K$ is its associated stacking parameter (see Definition \ref{def:stackingparameter}), then for each $\theta \in E_{N,K}^C$ one has
\begin{equation}\label{eq:propagationconclusion}
\mathcal{H}^1 \big(\pi_{\theta}(\mathcal{S}^{NJ_{K,\rho}}) \big) \lesssim \frac{1}{K}.
\end{equation}
\end{lemma}

We refer to Lemma \ref{lma:reverseholder} as the \textit{propagation inequality}.  The exceptional set is defined using scales at most $N$, whereas its complement gives a small projection after $J_{K,\rho}$ blocks of length $N$. By appealing to the propagation inequality, then, we thus obtain our first major reduction.

\medskip

\begin{proposition}\label{prop:favardreducedtoexceptional}
To establish our power laws for the Favard length problem as in Theorem \ref{thm:totalFavbound}, it is enough to show that there exists some $\beta \in (0,1)$ such that
\begin{equation}\label{eq:smallexceptionalset}
\sigma^{d-1} \big(E_{N,K}\big) \lesssim K^{-\beta}, \quad \forall N \gg 1,
\end{equation}
where $K = N^{\nu (\epsilon_0, N)}$ is the stacking parameter from Definition \ref{def:stackingparameter}.
\end{proposition}

\begin{proof}
This follows from the simple calculation
\begin{eqnarray*}
    \Fav (\mathcal{S}^{NJ_{K,\rho}}) & = & \bigg(\int_{\theta \in E_{N,K}} \mathcal{H}^1 \big(\pi_{\theta}( \mathcal{S}^{NJ_{K,\rho}}) \big) d \sigma^{d-1}(\theta) + \int_{\theta\in E_{N,K}^C} \mathcal{H}^1 \big(\pi_{\theta} (\mathcal{S}^{NJ_{K,\rho}}) \big) d \sigma^{d-1} (\theta) \bigg) \\[1ex]
    \quad & \lesssim & \bigg(K^{-1} + \sigma^{d-1} \big(E_{N,K}\big)\bigg) \\[1ex]
    \quad & \lesssim &  N^{- \beta \nu(\epsilon_0, N)}
\end{eqnarray*}
where in the final step we used the definition of the stacking parameter
$K=N^{\nu(\epsilon_0,N)}$.  The implicit constant may depend on the digit sets $A_1,...,A_d$, the ambient dimensional parameter $d \geq 2$, and the fixed value of $\rho$, but not on $N$.

It remains to pass from the propagated depths to all sufficiently large
construction depths.  Put
\[
 \mathfrak M(N):=NJ_{K(N),\rho}.
\]
In both branches of Definition \ref{def:stackingparameter}, the function
$\mathfrak M(N)$ is eventually increasing and
\[
 \frac{\mathfrak M(N+1)}{\mathfrak M(N)}\longrightarrow1.
\]
Indeed,
\[
 \frac{K(N+1)}{K(N)}\longrightarrow1,
 \qquad
 \frac{\log(2+K(N+1))}{\log(2+K(N))}\longrightarrow1,
\]
both for $K(N)=N^{\epsilon_0}$ and for
$K(N)=N^{\epsilon_0/\log\log N}$.  Since the ceiling in
\eqref{eq:propagationblocks} contributes a relative error $O(J_{K,\rho}^{-1})$,
these limits give
\[
 \frac{\mathfrak M(N+1)}{\mathfrak M(N)}
 =\frac{N+1}{N}
 \left(\frac{K(N+1)}{K(N)}\right)^{\rho-1}
 \frac{\log(2+K(N+1))}{\log(2+K(N))}(1+o(1))
 \longrightarrow1.
\]
Given a sufficiently large integer $M$, choose the largest $N$ such that
$\mathfrak M(N)\leq M$.  The Cantor iterates are nested, so
\[
 \Fav(\mathcal S^M)\leq\Fav(\mathcal S^{\mathfrak M(N)}).
\]
Maximality also gives
\[
 1\leq\frac{M}{\mathfrak M(N)}
 <\frac{\mathfrak M(N+1)}{\mathfrak M(N)}=1+o(1).
\]
Thus replacing $\mathfrak M(N)$ by $M$ does not change a power exponent, or
the leading coefficient in the $1/\log\log M$ branch.  The exact conversion
is recorded in Section \ref{sec:exponentexamples}.
\end{proof}

While we phrase the lower multiplicity estimate \eqref{eq:smallexceptionalset} for some arbitrary $\beta \in (0,1)$, in the course of our argumentation, we will show that $\beta$ may be taken arbitrarily close to $1$.  The precise effect of this idea, as well as the general input of all auxiliary parameters, on the final exponent is recorded in Section \ref{sec:exponentexamples}.

\begin{remark}\label{rmk:lossy}
\rm{
Recognize that the parameter $J_{K,\rho}$, defined in \eqref{eq:propagationblocks} and utilized in the statement of Lemma \ref{lma:reverseholder}, introduces critical loss in our final choice of power law exponent. Indeed, if one chooses $K=N^{\nu(\epsilon_0,N)}$, then
\[
NJ_{K,\rho}\asymp N^{1+(\rho-1)\nu(\epsilon_0,N)}\log(2+K).
\]
Thus propagation weakens the power obtained at scale $N$ by the exact factor of $(1+(\rho-1)\nu(\epsilon_0,N))^{-1}$, up to an arbitrarily small loss used to absorb the logarithm. This suggests that, if one could take $\rho \rightarrow 1^+$, then the above propagation would become lossless. However, the author is currently unaware of how to accomplish this task. We recall that Section \ref{sec:exponentexamples} records all necessary bookkeeping to support this observation.
}
\end{remark}

In addition to our propagation inequality, one of our main analytic tools is the following $L^2$ estimate, which holds for directions in our low multiplicity set $E_{N,K}$. While this result was already known in the plane (first in \cite{NPV} and then with an alternative proof given \cite[Section 5]{BThesis}), we will prove the analogous version for all $d \geq 2$ in Section \ref{sec:comblemmata} (specifically, see Section \ref{subsec:L2exceptional}) while also significantly improving the threshold of parameters utilized in the proof (and, thus, obtain an improved exponent over \cite{NPV} and subsequent results).

\medskip

\begin{lemma}\label{lma:L2exceptionalset}
    If $\theta \in E_{N,K}$, then
    \begin{equation}\label{eq:L2exceptionalset}
    \max_{1 \leq n \leq N} \vert \vert f_{n, \theta}\vert\vert_{L^2 (\mathbb{R})}^2 \leq C K
    \end{equation}
    where $C=C_{L,d,\rho}>0$ is independent of $N$, $K$, and $\theta$.
\end{lemma}

In the next section, we convert the previous $L^2$ upper bound on the counting function into an upper bound on a Riesz product of the form \eqref{eq:rieszproductfirsttime}. This is the main benefit in having an $L^2$ inequality, in that we can leverage Plancherel's theorem to transfer the above geometric estimate to a Fourier analytic estimate without any loss.

\subsection{High- and low-frequency factors of Riesz products}
\label{subsec:countingandtrig}

We now translate the geometric overlap estimate into a Fourier estimate. We use the convention
\[
\widehat{g}(\xi):=\int_{\mathbb R}g(x)e^{-2\pi ix\xi}\,dx.
\]

Fix one signed projective chart. As above, we absorb its fixed coordinate reflections into the digit sets and then relabel the reflected sets as $A_1,\ldots,A_d$. For each single-scale digit
\[
a=(a_1,\ldots,a_d)\in\mathcal A
:=A_1\times\cdots\times A_d
\]
and to each flattened direction $t=(t_1,\ldots,t_{d-1})\in[0,1]^{d-1}$, define the projected one-step digit
\[
X_t(a):=a_1+t_1a_2+\cdots+t_{d-1}a_d.
\]
Recall that the associated normalized trigonometric polynomial is
\begin{equation}\label{eq:phit-probabilistic}
\phi_t(u)
:=\frac1L\sum_{a\in\mathcal A}e^{2\pi iuX_t(a)}
=\phi_{A_1}(u)\phi_{A_2}(t_1u)\cdots
 \phi_{A_d}(t_{d-1}u).
\end{equation}
Thus $\phi_t$ is the characteristic function of the uniformly distributed projected digit, with the positive-exponent convention used in Section \ref{subsec:digitsetpoly}.

For notational simplicity, continue to write $E_{N,K}$ for the pullback of the exceptional set to the fixed parameter cube. The comparison \eqref{eq:chartmeasurecomparison} shows that this changes measure estimates by at most a dimensional constant. Lemma \ref{lma:L2exceptionalset}, the convolution representation \eqref{eq:flattenedcountingfunction}, and Plancherel's theorem give
\begin{equation}\label{eq:exactplancherelcounting}
CK\geq\|f_{N,t}\|_2^2
=\int_{\mathbb R}
 \left|\widehat{\nu_{N,t}}(\xi)\right|^2
 \left|
 L^N\widehat{\mathbf 1}_{[-b_{N,t},b_{N,t}]}(\xi)
 \right|^2\,d\xi.
\end{equation}
Here we used
\[
\left|\widehat{\nu_{N,t}^{\mathrm c}}(\xi)\right|
=\left|\widehat{\nu_{N,t}}(\xi)\right|,
\]
since centering the projected intervals introduces only a unimodular phase. It is important for us to retain the exact smoothing radius
\[
b_{N,t}=\frac{\sigma_t}{2}L^{-N},
\]
for, indeed, one can check that
\begin{equation}\label{eq:exactsmoothingmultiplier}
L^N\widehat{\mathbf 1}_{[-b_{N,t},b_{N,t}]}(\xi)
=\sigma_t
 \frac{\sin(\pi\sigma_tL^{-N}\xi)}
 {\pi\sigma_tL^{-N}\xi},
\end{equation}
where the quotient is interpreted as one at the origin. Since
$1\leq\sigma_t\leq d$,
for $|\xi|\leq L^{N/2}$ and all sufficiently large $N$ we have
\[
|\pi\sigma_tL^{-N}\xi|
\leq\pi dL^{-N/2}\leq\frac{\pi}{2}.
\]
A straightforward calculation thus shows that
\begin{equation}\label{eq:smoothingmultiplierlower}
\left|
L^N\widehat{\mathbf 1}_{[-b_{N,t},b_{N,t}]}(\xi)
\right|
\geq\frac{2\sigma_t}{\pi}
\geq\frac2\pi,
\qquad |\xi|\leq L^{N/2}.
\end{equation}

We thus restrict our nonnegative integral in
\eqref{eq:exactplancherelcounting} to this frequency range, applying
\eqref{eq:smoothingmultiplierlower}, and averaging over the chartwise exceptional set, to obtain that 
\begin{equation}\label{eq:prepigeonholeK}
\frac{1}{\mathcal H^{d-1}(E_{N,K})}
\int_{t\in E_{N,K}}
\left[
 \sum_{r=1}^{\lfloor N/2\rfloor}
 \int_{L^{r-1}}^{L^r}
 \left|\widehat{\nu_{N,t}}(\xi)\right|^2\,d\xi
\right]dt
\leq CK.
\end{equation}

We next remove the low-frequency portion of the region of integration by pigeonholing relative to the following parameter.

\medskip

\begin{definition}\label{def:littlem}
Let $s_0>0$ be independent of $N$, to be chosen later. If
$A_i^{(3)}\equiv1$ for every $i=1,\ldots,d$, define
\[
m:=\left\lceil s_0\log_LN\right\rceil.
\]
Otherwise, define
\[
m:=
\left\lceil
s_0\frac{\log_LN}{\log_L\log_LN}
\right\rceil.
\]
\end{definition}

The parameter $s_0$ is the analytic scale-selection parameter. It is distinguished from the lowercase letter $s$, which is reserved throughout the paper for a cyclotomic index. The admissible range of $s_0$ is determined by the constants associated to the Set of Small Values and Set of Large Values constructions, and will be optimized in Section \ref{sec:exponentexamples}.

We now produce a frequency pigeonholing relative to this chosen parameter
\[
m = m(N; A_1,...,A_d).
\]
First, select a floating parameter $r \in \mathbb{N}$ which satisfies the inequality
\[
\left\lceil\frac N4\right\rceil
\leq r\leq
\left\lfloor\frac N2\right\rfloor,
\]
and consider the (by definition, overlapping) frequency windows
\[
W_r : =[L^{r-m},L^r], \textrm{ for all admissible } r.
\]
A simple covering argument shows that each $\xi > 0$ can belong to, at the absolute most, ($m+1$)-many windows $W_r$ simultaneously. Hence, one obtains that 
\[
\sum_{r=\lceil N/4\rceil}^{\lfloor N/2\rfloor}
\int_{W_r}
\left|\widehat{\nu_{N,t}}(\xi)\right|^2\,d\xi
\leq
(m+1)\int_1^{L^{N/2}}
\left|\widehat{\nu_{N,t}}(\xi)\right|^2\,d\xi.
\]
After averaging in $t$, we may apply the inequality \eqref{eq:prepigeonholeK} to obtain that the right-hand side quantity is bounded above by, at most, $C(m+1)K$. If $N$ is chosen large enough, then there are (say) at least $N/5$ admissible values of $r$; and thus, there exists at least one admissible integer
\[
\frac N4\leq n\leq\frac N2
\]
such that
\begin{equation}\label{eq:averagingoverbadpigeonhole}
\frac{1}{\mathcal H^{d-1}(E_{N,K})}
\int_{E_{N,K}}
\left[
 \int_{L^{n-m}}^{L^n}
 \left|\widehat{\nu_{N,t}}(\xi)\right|^2\,d\xi
\right]dt
\leq\frac{CKm}{N}.
\end{equation}
We fix one such integer $n \in [N/4, N/2]$ and thus proceed with our argument localized to the frequency window $W_n : = [L^{n - m}, L^n]$.

Having localized our argument to a slightly smaller band of frequencies, we next justify removing all oscillatory factors whose frequencies occur at scales $k > n$. There are two different scale orderings in this calculation, and it is useful to distinguish them. In the original product
\[
\widehat{\nu_{N,t}}(\xi)
=\prod_{k=1}^N\phi_t(-L^{-k}\xi),
\]
the indices $k>n$ correspond to the finest spatial generations. On any frequency interval of the form
\[
0\leq\xi\leq c_*L^n,
\]
these factors are thus evaluated at very small arguments. They may therefore be removed at a fixed multiplicative cost. To make the small-argument estimate precise, define
\begin{equation}\label{eq:projected-mean-variance}
\mu_t
:=\frac1L\sum_{a\in\mathcal A}X_t(a),
\qquad
\operatorname{Var}_t(X)
:=\frac1L\sum_{a\in\mathcal A}
\bigl(X_t(a)-\mu_t\bigr)^2.
\end{equation}
Taking a Taylor expansion of the function $\phi_t (u)$ at the origin produces the following natural asymptotic identity
\[
\phi_t(u)
=
1+2\pi i\mu_tu
-2\pi^2
 \left(\frac1L\sum_{a\in\mathcal A}X_t(a)^2\right)u^2
+O(|u|^3),
\]
and so, as a consequence, we see that
\begin{equation}\label{eq:phit-modulus-expansion}
|\phi_t(u)|^2
=
1-4\pi^2\operatorname{Var}_t(X)u^2
+O(|u|^3).
\end{equation}
As the underlying digit set $\mathcal{A} := A_1 \times \cdots \times A_d$ is a fixed finite set and the associated parameter cube $t \in [0,1]^{d-1}$ is compact, the $X_t$ and the implicit remainder in
\eqref{eq:phit-modulus-expansion} are necessarily uniformly bounded. It follows that there exist constants $r_0,C_0>0$, which are independent of $t$ and of the signed chart, such that
\begin{equation}\label{eq:digitsetnearzero}
|\phi_t(u)|\geq\exp(-C_0|u|^2),
\qquad |u|\leq r_0.
\end{equation}

We thus fix a constant
\begin{equation}\label{eq:safefrequencycutoff}
0<c_*<\min\{1,r_0L\},
\end{equation}
and further restrict the region of integration
\eqref{eq:averagingoverbadpigeonhole} to the interval
\[
[L^{n-m},c_*L^n]
\]
and thus define 
\[
E_{N,K}^*
:=
\left\{
t\in E_{N,K}:
\int_{L^{n-m}}^{c_*L^n}
\left|\widehat{\nu_{N,t}}(\xi)\right|^2\,d\xi
\leq\frac{2CKm}{N}
\right\} \subseteq E_{N,K}.
\]
We then let
\[
F(t)
:=
\int_{L^{n-m}}^{c_*L^n}
\left|\widehat{\nu_{N,t}}(\xi)\right|^2\,d\xi
\]
and put
\[
A:=\frac{CKm}{N}.
\]
Then \eqref{eq:averagingoverbadpigeonhole} gives
\begin{equation}\label{eq:average-F-bound}
\int_{E_{N,K}}F(t)\,d\mathcal H^{d-1}(t)
\leq
A\,\mathcal H^{d-1}(E_{N,K}),
\end{equation}
so that, in this new notation, we have that
\[
E_{N,K}^*
:=
\left\{
t\in E_{N,K}:F(t)\leq2A
\right\}.
\]
We claim that $E_{N,K}^*$ contains at least half the measure of
$E_{N,K}$. Indeed, this follows by Markov's inequality together with \eqref{eq:average-F-bound}, since we then must have that
\begin{eqnarray*}
\mathcal H^{d-1}
\bigl(E_{N,K}\setminus E_{N,K}^*\bigr)
&=&
\mathcal H^{d-1}
\left(
\{t\in E_{N,K}:F(t)>2A\}
\right)\\[1ex]
&\leq&
\frac{1}{2A}
\int_{E_{N,K}}F(t)\,d\mathcal H^{d-1}(t)\\[1ex]
&\leq&
\frac{1}{2A}
\left(
A\,\mathcal H^{d-1}(E_{N,K})
\right)\\[1ex]
&=&
\frac12\mathcal H^{d-1}(E_{N,K}).
\end{eqnarray*}
Consequently,
\begin{equation}\label{eq:good-exceptional-subset}
\mathcal H^{d-1}(E_{N,K}^*)
\geq
\frac12\mathcal H^{d-1}(E_{N,K}).
\end{equation}
Thus the restriction to $E_{N,K}^*$ discards at most half of the
original exceptional parameter set, while providing the pointwise
Fourier upper bound
\[
\int_{L^{n-m}}^{c_*L^n}
\left|\widehat{\nu_{N,t}}(\xi)\right|^2\,d\xi
\leq
\frac{2CKm}{N}
\qquad
\text{for every }t\in E_{N,K}^*.
\]

As a final reduction, we show that any oscillatory factor at scale $n < k \leq N$ can be suppressed. To this end, fix some
$
\xi\leq c_*L^n
$
and write $k=n+r$, where $r\geq1$. Then
\begin{equation}\label{eq:tailargumentsmall}
|L^{-k}\xi|
\leq c_*L^{-r}
\leq c_*L^{-1}
<r_0,
\end{equation}
where the final inequality follows by our choice of constant $c_*$. It then follows from \eqref{eq:digitsetnearzero} that
\begin{align*}
\prod_{k=n+1}^N|\phi_t(L^{-k}\xi)|
&\geq
\exp\left(
 -C\sum_{r=1}^{N-n}|L^{-(n+r)}\xi|^2
\right)\\
&\geq
\exp\left(
 -Cc_*^2\sum_{r=1}^{\infty}L^{-2r}
\right)
=:c_{\mathrm{tail}}>0.
\end{align*}
Hence, the constant $c_{\mathrm{tail}}$ depends only on the fixed digit data and the chosen cutoff.
Therefore, for each $\xi \in [L^{n-m},c_*L^n]$, we have shown that
\begin{equation}\label{eq:tailcomparisonpointwise}
\left|\widehat{\nu_{N,t}}(\xi)\right|^2
\approx_{c_*, A_1,...,A_d}
\left|
 \prod_{k=1}^n\phi_t(L^{-k}\xi)
\right|^2.
\end{equation}
Thus the generations $n+1,\ldots,N$ may be removed without any scale-dependent loss.
The role of the cutoff $c_*$ is now explicit. It ensures that every omitted factor is evaluated inside the fixed neighbourhood of the origin on which \eqref{eq:digitsetnearzero} holds.

We conclude our Fourier analytic set-up by performing a natural change of variables. Letting
$
\xi=L^nu
$
we then have that
\[
\prod_{k=1}^n\phi_t(L^{-k}\xi)
=
\prod_{j=0}^{n-1}\phi_t(L^ju),
\qquad
d\xi=L^n\,du.
\]
Combining this identity with
\eqref{eq:tailcomparisonpointwise} and the definition of
$E_{N,K}^*$ gives, for every $t\in E_{N,K}^*$,
\begin{equation}\label{eq:modifiedexcepsetineq}
\int_{L^{-m}}^{c_*}
\left|
 \prod_{j=0}^{n-1}\phi_t(L^ju)
\right|^2\,du
\lesssim
\frac{Km}{L^nN}.
\end{equation}
Hence, all of our subsequent lower-bound arguments are carried out on the same frequency interval
\[
[L^{-m},c_*],
\]
which we will refer to as our \textbf{safe frequency interval}.  Thus
\eqref{eq:modifiedexcepsetineq} isolates precisely the portion of the Riesz product upon which our more delicate analysis (and, in particular, our Set of Small Values and Set of Large Values constructions) must be conducted.

\section{Lower Bounds on Riesz Products}\label{sec:rieszproductbounds}

The preceding section converted the assumption of low projection multiplicity into the upper bound \eqref{eq:modifiedexcepsetineq}. To obtain a contradiction, we must now prove that the same rescaled Riesz product has substantially more $L^2$ mass on the identified safe frequency interval $[L^{-m}, c_*]$. Our main difficulty in proving a lower bound is that the remaining high and low frequencies of the Riesz product can become small for different arithmetic reasons; and, even more, depending upon the cyclotomic structure of our digit sets $A_i$ and their corresponding mask polynomials $A_i (X)$, we may further encounter several distinct types of zeroes and small-values pointwise behaviour in our Riesz products. To handle these issues, we separate both the frequency scales and the algebraic factors of our Riesz product before constructing the lower bound.

Write
\[
 P_{1,t}(\xi):=\prod_{k=m}^{n-1}\phi_t(L^k\xi),
 \qquad
 P_{2,t}(\xi):=\prod_{k=0}^{m-1}\phi_t(L^k\xi).
\]
The first factor contains the high frequencies and the second contains the low
frequencies, all adapted to our pigeonholed scale parameter $n \in [N/4, N/2]$.  The cyclotomic factorization (see Definition \ref{def:cycfactorization}) further gives corresponding decompositions
\[
 P_{j,t}=P_{j,t}'P_{j,t}'',
 \qquad j=1,2,
\]
which corresponds the splitting of the mask polynomials $A(X) : = A' (X) A'' (X)$.

The two decompositions serve different purposes. The split $P_{1,t}P_{2,t}$ separates the generations of the cylinder construction: $P_{2,t}$ contains the first $m$-many frequency scales and will be bounded below pointwise; whereas, the factor $P_{1,t}$ contains the remaining $(n-m)$-many scales and supplies the averaged $L^2$ mass. The split into primed and double-primed factors separates the zero sets: the primed zeros are removed by a small exceptional set of frequencies, whereas the accumulating double-primed zeros are avoided by placing the argument on the difference set of a suitably constructed large-value set. These two splits correspond to the Set of Small Values technique (for the primed factor) and the Set of Large Values technique (for the double-primed factor). In this section, we will assume that both the SSV and SLV calculations have been accomplished, and prove the lower bound. The remaining work towards these two intermediate results is then accomplished in Sections \ref{sec:ssvtechnical} and \ref{sec:multiSLV}, respectively.

\subsection{Poisson localization}\label{subsec:poisson}
Before examining the structure of our Riesz product on the safe frequency interval, we first need to bound exactly how much $L^2$ mass has been lost by restricting to this safe interval. We accomplish this via the following upper bound, which is pointwise in the flattened direction parameter $t \in E_{N,K}^*$.

\begin{proposition}\label{prop:shortintervalbound}
Suppose that $t \in E_{N,K}^*$. Then,
\begin{equation}\label{eq:shortintervalbound}
    \int_0^{L^{-m}} \lvert P_{1,t} (\xi) \rvert^2 d \xi  = \int_0^{L^{-m}} \bigg\lvert \prod_{k = m}^{n - 1} \phi_t (L^k \xi )\bigg\vert^2 d \xi  \leq C_0 K L^{-n}
\end{equation}
for some constant $C_0=C_{L,d,\rho}>0$ independent of $N$, $K$, and $t$.
\end{proposition}

The proof we present requires the following Lemma, which is proven in detail in \cite[Section 5.4]{BThesis} and first appeared in \cite{BV3}. The lemma is a version of the Fourier uncertainty principle adapted to a discrete family of intervals. 

\medskip

\begin{lemma}\label{lma:poissonlocalization}
Let $\delta>0$, let $\mathcal J$ be a finite index set, and let
$(\alpha_j)_{j\in\mathcal J}$ be an indexed family of real numbers; repetitions
are allowed.  Suppose
$$
\int_{\mathbb{R}} \bigg(\sum\limits_{j\in\mathcal J} \mathbf{1}_{[\alpha_j - \delta, \alpha_j + \delta]} (x) \bigg)^2 dx \leq S
$$
for some $S > 0$ which may depend upon $\delta > 0$. Then there exists an absolute constant $C > 0$ such that for any $\xi_0 \in \mathbb{R}$ one has
$$
\delta^2 \int_{\xi_0}^{\xi_0 + \delta^{-1}} \bigg\lvert \sum\limits_{j\in\mathcal J} c_j e^{2\pi i \alpha_j \xi} \bigg\rvert^2 d \xi \leq C S,
$$
where $(c_j)_{j\in\mathcal J}\subset\mathbb C$ is any family satisfying
$|c_j|=1$.  Thus the same symbolic multiplicities occur in the spatial and
Fourier sums.
\end{lemma}

\begin{proof}[Proof of Proposition \ref{prop:shortintervalbound}]
Let $\mathcal W_{n-m}:=\mathcal A^{n-m}$ be the set of symbolic words of
length $n-m$; and, for each $\omega\in\mathcal W_{n-m}$, let $\alpha_\omega$ be
the corresponding flattened projected vertex; so that
$$
\alpha_\omega = z_{1,\omega} + t_1 z_{2,\omega} + \cdots + t_{d-1} z_{d,\omega}
$$
where $z_\omega \in Q_{\omega}$ is the corresponding minimal vertex of the associated cylinder cube.
Using this notation convention, we thus have that
\begin{equation}\label{eq:poissonloc1}
\prod_{k = 1}^{n - m}\phi_t (L^{-k} \xi) = \widehat{\nu_{n - m, t}} (\xi)
 =L^{-(n-m)} \sum\limits_{\omega\in\mathcal W_{n-m}} e^{- 2 \pi i \alpha_\omega \xi},
\end{equation}
where, we remind ourselves, the rescaled convolution $\nu_{n-m,t}$ is defined in the identities \eqref{eq:bignu} and
\eqref{eq:littlenu}.
If $t\in E_{N,K}^*\subset E_{N,K}$, Lemma \ref{lma:L2exceptionalset} and
\eqref{eq:flattenedcountingfunction} imply together that
$$
\int_{\mathbb R}\left(
 \sum_{\omega\in\mathcal W_{n-m}}
 \mathbf1_{[\alpha_\omega-L^{-(n-m)},\alpha_\omega+L^{-(n-m)}]}(x)
 \right)^2dx\lesssim K.
$$
Consequently, Cauchy--Schwarz and translation invariance give
\[
 \left\|\sum_\omega
 \mathbf1_{[\alpha_\omega-\delta,\alpha_\omega+\delta]}
 \right\|_2^2
 \leq R_d\sum_{r=1}^{R_d}\|f_{n-m,t}(\cdot-h_r)\|_2^2
 \lesssim_d K,
\]
where $R_d : = O_d(1)$ essentially counts the number of projected intervals of side-length $[-\sigma_t \delta/2, \sigma_t \delta /2]$ are needed to cover the interval $[-\delta, \delta]$.  Making the change of variables
$\xi\mapsto L^n\xi$ in \eqref{eq:shortintervalbound}, we obtain
\begin{eqnarray*}
\int_0^{L^{-m}} \bigg\lvert \prod_{k = m}^{n - 1} \phi_t (L^k \xi )\bigg\vert^2 d \xi & = & L^{-n} \int_{0}^{L^{n - m}} \bigg\lvert \prod_{k = 1}^{n - m} \phi_t (L^{-k} \xi) \bigg\rvert^2 d \xi \\[1ex]
\quad & = & L^{-n} \bigg( L^{-2(n-m)} \int_0^{L^{n-m}} \bigg\lvert \sum\limits_{\omega\in\mathcal W_{n-m}} e^{- 2 \pi i \alpha_\omega \xi} \bigg\rvert^2 d \xi \bigg) \\[1ex]
\quad & \leq & C_0 L^{-n} K
\end{eqnarray*}
where the final inequality comes by applying Lemma \ref{lma:poissonlocalization} with $\delta = L^{-(n-m)}$ and $S = K$.
\end{proof}

The normalization in this calculation is worth recording. The unnormalized exponential sum contains $L^{n-m}$ terms, while the probability measure contributes the factor $L^{-(n-m)}$. The factor $\delta^2=L^{-2(n-m)}$ in Lemma \ref{lma:poissonlocalization} exactly matches the square of this normalization. Consequently, Poisson localization produces the scale $KL^{-n}$ required by the later Riesz-product comparison, with no additional power of $L^m$. Thus the discarded interval contributes precisely the error term $C_0KL^{-n}$ which the later Set of Large Values gain must dominate.

\subsection{Applying the Set of Small Values and Set of Large Values properties}

We are now in a position to prove all necessary bounds for our Riesz products, under the assumption that the set of low multiplicity directions is large, and the condition that the SSV and SLV properties are in place for the function $\phi_t$ and for all directions $t \in E_{N,K}^*$. As explained in Proposition \ref{prop:favardreducedtoexceptional}, this suffices to prove our Favard length estimate.

We phrase the hypotheses of the following proposition as a list of conditions, which will then be supplied by the Set of Small Values and Set of Large Values properties in the following two sections.

\begin{proposition}[Riesz-product lower bound]
\label{prop:keyproposition}
Fix $t\in E_{N,K}^*$.  Suppose that there are constants
$\eta\in(0,1]$, $C_1\geq0$, $c_{1,1},\ldots,c_{1,d}\geq0$, and
$c_\Gamma>0$ such that the following conditions are simultaneously satisfied.
\begin{enumerate}
\item There exists a measurable set $\Gamma\subset[0,c_*/2]$ whose measure satisfies
$$
\mathcal{H}^1 (\Gamma) \geq c_\Gamma L^{-(1-\eta)m}
$$
and such that the low-frequency accumulated zeroes of our Riesz product satisfy the pointwise estimate
\[
|P_{2,t}''(\xi)|\geq L^{-C_1m},
\qquad
\forall\xi\in\Gamma-\Gamma.
\]
\item There exists a symmetric and measurable set $\operatorname{SSV}(t)\subset[-c_*,c_*]$, as well as two associated sets
\[
\mathcal B_t
:=
[-L^{-m},L^{-m}]
\cup
\operatorname{SSV}(t),
\qquad
\Omega_t
:=
(\Gamma-\Gamma)\setminus\mathcal B_t,
\]
as well as a fixed constant $C_{\mathrm{err}}>0$, such that the high-frequency portion of our Riesz product satisfies
\[
\int_{(\Gamma-\Gamma)\cap (\operatorname{SSV}(t) \cup [-L^{-m}, L^{-m}])}
|P_{1,t}(\xi)|^2\,d\xi
\leq
C_{\mathrm{err}}KL^{-n},
\]
and such that the non-accumulating part of the low-frequency product satisfies
\[
|P_{2,t}'(\xi)|
\geq
L^{-m\sum_{i=1}^dc_{1,i}},
\qquad
\forall\xi\in \Gamma - \Gamma \setminus \big(\operatorname{SSV}(t) \cup [-L^{-m}, L^{-m}]\big)
\]
\item Finally, assume that the associated parameters satisfy the following inequality
\[
c_\Gamma L^{\eta m}
\geq
2(C_{\mathrm{err}}+1)K.
\]
\end{enumerate}

\noindent
Then, by letting
\[
D:=C_1+\sum_{i=1}^dc_{1,i},
\]
we then have the following lower bound for our Riesz product
\begin{equation}
\label{eq:parameterneutralRieszlower}
\int_{L^{-m}}^{c_*}
|P_{1,t}(\xi)P_{2,t}(\xi)|^2\,d\xi
\gtrsim
KL^{-n}L^{-2Dm}.
\end{equation}
\end{proposition}

\begin{proof}
We first use Item $(1)$ to construct a function $h :[-c_*, c_*] \rightarrow [0,1]$ with certain special properties. For us, the correct choice is to take
$$
h:=\frac{\mathbf1_\Gamma*\mathbf1_{-\Gamma}}{\mathcal{H}^1(\Gamma)}.
$$
The function $h$ is thus supported in the topological closure of the difference set $\Gamma - \Gamma \subseteq [-c_*,c_*]$, and also has the following key properties:
$$
 0\leq h(\xi) \leq 1, \qquad \forall \xi \in [-c_*,c_*]
$$
$$
 \widehat h(\xi)= \frac{|\widehat{\mathbf1_\Gamma}(\xi)|^2}{\mathcal{H}^1(\Gamma)} \geq 0, \qquad \forall\xi \in \mathbb R
$$
$$
\widehat h(0)=\int_{\mathbb{R}} h (\xi) d \xi =\mathcal{H}^1 (\Gamma).
$$
Using $\lambda_\omega$ to denote the frequencies (with multiplicity) associated to each word$\omega\in\mathcal A^{n-m}$, we then have that
\[
 P_{1,t}(\xi)=L^{-(n-m)}\sum_{\omega\in\mathcal A^{n-m}}
 e^{2\pi i\lambda_\omega\xi}.
\]
Hence, and using that our Riesz product is a long trigonometric polynomial written in frequencies $\lambda_{\omega}$ for $\omega \in \mathcal{A}^{n-m}$, we have that
\begin{align*}
 \int_{\mathbb R}|P_{1,t}(\xi)|^2h(\xi)\,d\xi
 &=L^{-2(n-m)}\sum_{\omega,\tau\in\mathcal A^{n-m}}
 \widehat h(\lambda_\tau-\lambda_\omega)\\
 &\geq L^{-2(n-m)}\sum_{\omega\in\mathcal A^{n-m}}\widehat h(0)
 =L^{-{n-m}}\mathcal{H}^1 (\Gamma),
\end{align*}
and here, we only use that $\widehat h\geq0$ and that any repetitions among the frequencies
$\lambda_\omega$ can only increase the value of the sum.

Since, by construction, we have that $0 \leq h\leq1$ and $h$ is supported on the closure of $\Gamma - \Gamma$, the size hypothesis on $\Gamma$ yields
\begin{equation}\label{eq:salemparameterneutral}
 \int_{\Gamma-\Gamma}|P_{1,t}(\xi)|^2\,d\xi
 \geq L^{-(n-m)} \mathcal{H}^1 (\Gamma)
 \geq c_{\Gamma}L^{-n}L^{\eta m}.
\end{equation}

By Item $(3)$, we are now guaranteed that
$c_{\Gamma}L^{\eta m}\geq2(C_{\mathrm{err}}+1)K$. In particular, by the upper bound guaranteed in Item $(2)$, we now know that  
\[
 \int_{(\Gamma-\Gamma)\setminus
  ([-L^{-m},L^{-m}]\cup\operatorname{SSV}(t))}
 |P_{1,t}(\xi)|^2\,d\xi
 \gtrsim KL^{-n}.
\]
By applying the pointwise bounds on $P_{2,t}'$ and $P_{2,t}''$ from Items $(2)$ and $(1)$, respectively, we thus obtain that
\begin{eqnarray*}
 \int_{\Gamma - \Gamma \setminus (\operatorname{SSV}(t) \cup [-L^{-m}, L^{-m}])}\lvert P_{1,t} (\xi)\rvert^2 |P_{2,t}(\xi)|^2
 & \gtrsim &  K L^{-n} L^{-2C_1m}L^{-2m\sum_i c_{1,i}} \\[1ex]
 & = & K L^{-n} L^{-2Dm}.
\end{eqnarray*}
To pass from this modified set to the desired interval $[L^{-m}, c_*]$, observe that both the difference set $\Gamma-\Gamma$ as well as the set $\operatorname{SSV}(t)$ are origin-symmetric subsets of $[-c_*, c_*]$. Since the integrand is an even function, and our underlying region of integration is origin symmetric, we must then have that
one half of the remaining mass lies exactly in the interval $[L^{-m},c_*]$; and, thus, restricting our estimate to this interval proves the lower bound
\eqref{eq:parameterneutralRieszlower}.
\end{proof}

\subsection{The final lower bound calculation}

Combining \eqref{eq:modifiedexcepsetineq} and
\eqref{eq:parameterneutralRieszlower} gives
\[
 KL^{-n}L^{-2Dm}
 \lesssim\frac{Km}{L^nN},
\]
or equivalently
\begin{equation}\label{eq:finalparametercontradiction}
 N L^{-2Dm}\lesssim m.
\end{equation}
In the power-law branch, take
\[
 m=\lceil s_0\log_LN\rceil.
\]
If $2Ds_0<1$, the left-hand side of
\eqref{eq:finalparametercontradiction} cannot be bounded by its right-hand
side.  In fact, the ceiling in $m$ gives
\[
 L^{-2Dm}\asymp N^{-2Ds_0},
 \qquad m\asymp\log N,
\]
so
\[
 \frac{NL^{-2Dm}}m
 \gtrsim\frac{N^{1-2Ds_0}}{\log N}\longrightarrow\infty.
\]
If $K=N^\kappa$, the simultaneous SLV requirement is
\[
 \frac{L^{\eta m}}K
 \asymp N^{\eta s_0-\kappa}\longrightarrow\infty,
\]
equivalently $\kappa<\eta s_0$.  Section \ref{sec:exponentexamples}
optimizes the two strict inequalities $\kappa<\eta s_0$ and $2Ds_0<1$.

Thus two competing requirements determine the eventual exponent. Taking $m$ larger makes the SLV set strong enough to dominate the error term $K$, since its gain is $L^{\eta m}$. However, increasing $m$ also worsens the pointwise loss $L^{-2Dm}$. The exponent calculation in Section \ref{sec:exponentexamples} chooses $m$ at the smallest logarithmic scale which balances these two effects.

\section{Small Values and Quantitative SSV Theory}\label{sec:ssvtechnical}

The purpose of this section is to isolate, and then quantify, the exact sub-level set bound which is needed to produce our Riesz-product lower bound from the previous Section \ref{sec:rieszproductbounds}.  
The property which will furnish this sub-level set estimate is called the \textit{Set of Small Values (SSV)} property, which is (in essence) a quantitative covering result for the set of approximate zeroes of the low frequency portion of our Riesz product.

\subsection{The Set of Small Values property and associated parameters}\label{subsec:SSVSLV}
We begin by stating the Set of Small Values property in the abstract, as a covering theorem for the sub-level set of an arbitrary Riesz product.

\begin{definition}\label{def:setofsmallvalues}
Let $L>1$ and let $\upsilon:\mathbb R\to\mathbb C$ be any $1$-periodic function.  We say that
$\upsilon$ has the \textbf{Set of Small Values property} with threshold function $\psi(m)$ if there are
constants $c_2,c_3>0$ and $C\geq1$, independent of $m$, such that, for every
sufficiently large integer $m$, the set
\[
 B_m(\upsilon,\psi)
 :=\left\{\xi\in[0,1]:
 \left|\prod_{k=0}^{m-1}\upsilon(L^k\xi)\right|\leq\psi(m)
 \right\}
\]
is covered by at most $CL^{c_2m}$ intervals of length at most
$CL^{-c_3m}$.  The constants $c_2,c_3$ and $C$ may all depend upon the choice of $\upsilon$ and $L$, but may not depend upon $m$ or on any later construction scale.
\end{definition}

We emphasize the order of quantifiers in Definition \ref{def:setofsmallvalues}.  The function $\upsilon$, the base $L$, and the constants $c_2,c_3,C$ are chosen first.  With these quantities fixed, the same covering estimate must hold for every sufficiently large $m$.

For us, our goal is to show that, for any choice of flattened direction $t \in [0,1]^{d-1}$ and appropriately chosen SSV parameters, we may take $\upsilon : = \phi_t'$. However, in order for the SSV property to be verified, this means that absolutely none of the SSV constants, nor the threshold function $\psi$, may be prescribed in such a way that they might depend upon the pigeonholed depth parameters $n$ or $m$, nor the terminal depth $N$, nor the direction parameter $t$. Having this kind of uniformity in both direction and scale, then, is what then permits us to utilize the SSV property \textit{simultaneously} on every signed projective chart, and at all intermediate scales $n,N$ and $m$. 

The remainder of this section is thus devoted to establishing an appropriately strong version of the SSV property for the function $\phi_t'$, while also conducting careful optimization and book-keeping of the involved parameters $\psi$, $c_2, c_3$ and $C$. The beginning of our argument is the following single digit-set result, which supplies the initial SSV input used throughout the remainder of the paper.

\begin{theorem}[The Set of Small Values Property for Individual Digit Sets]\label{thm:onedimensionalSSV}
Let $A$ be a finite digit set such that $\#A$ divides $L$, and let 
$$
\phi_A'(\xi) : = \frac{1}{A'(1)} A'\big(e^{2 \pi i \xi}\big) \quad \forall\xi \in [0,1].
$$
be the normalized trigonometric
polynomial associated to the polynomial $A'=A^{(1)}A^{(3)}A^{(4)}$ (see Definition \ref{def:cycfactorization}). We then have the following SSV properties for $\phi_{A}'$.

\begin{enumerate}
\item{
The function $\phi_A'$ has the
log-SSV property with a threshold of the form
\[
 \psi(m)=L^{-c_1m\log_Lm},
\]
where $c_1 \geq 1/2$ is a constant whose value will depend explicitly upon the underlying SSV constants $c_2, c_3$.}

\medskip

\item{If $A^{(3)}\equiv1$, so that every unit-circle zero of $A$ is a root of unity,
then $\phi_A'$ has the ordinary SSV property with a threshold of the form
\[
 \psi(m)=L^{-c_1m}.
\]
where $c_1 \geq 1/2$ is a constant whose value will depend explicitly upon the underlying SSV constants $c_2, c_3$.
}
\medskip

\item{In either case, given any fixed proportion constant $\Lambda>1$, the threshold constant $c_1$ may
be enlarged so that the covering constants satisfy $c_3/c_2 \geq\Lambda$.
}
\end{enumerate}
\noindent
In all stated SSV properties, the underlying constants depend only upon the fixed digit set $A$ and upon the integer base $L$, and never upon the scale parameter $m$.
\end{theorem}

\begin{proof}
We first reconcile our notation with the normalization used in \cite{BLV}. Suppose that $Q(1)\neq 0$ and put
\[
 u(\xi):=Q(e^{2\pi i\xi}),
 \qquad
 \widetilde u(\xi):=\frac{Q(e^{2\pi i\xi})}{Q(1)}.
\]
For every $m\geq1$, one has the exact identity
\[
 \prod_{k=0}^{m-1}\widetilde u(L^k\xi)
 =Q(1)^{-m}\prod_{k=0}^{m-1}u(L^k\xi).
\]
Suppose first that the unnormalized product has the ordinary SSV property at threshold $L^{-am}$. Then
\[
 \left|\prod_{k=0}^{m-1}\widetilde u(L^k\xi)\right|
 \leq L^{-(a+\log_L|Q(1)|)m}
 \quad\Longleftrightarrow\quad
 \left|\prod_{k=0}^{m-1}u(L^k\xi)\right|
 \leq L^{-am}.
\]
Thus normalization increases the ordinary threshold exponent by precisely the fixed quantity $\log_L|Q(1)|$. If the unnormalized threshold is instead $L^{-am\log_Lm}$, then, for every fixed $\varepsilon>0$ and all sufficiently large $m$,
\[
 |Q(1)|^mL^{-(a+\varepsilon)m\log_Lm}
 \leq L^{-am\log_Lm}.
\]
Hence an arbitrarily small fixed enlargement of the logarithmic threshold exponent is sufficient. Moreover, the product indexed by $k=0,\ldots,m-1$ is obtained from the product indexed by $k=1,\ldots,m$ by the change of variables $\xi\mapsto \xi/L$. This change multiplies covering lengths by the fixed factor $L$ and does not alter either SSV exponent.

We next compare the cyclotomic conventions. Since $\#A$ divides $L$, one has
\[
 \gcd(s,\#A)\neq1
 \quad\Longrightarrow\quad
 \gcd(s,L)\neq1.
\]
Consequently, \cite[Proposition 4.2]{BLV} applies to every cyclotomic factor in $A^{(1)}$, while \cite[Proposition 4.3]{BLV} applies to every irreducible factor in $A^{(3)}$. Since $A^{(4)}$ has no zero on the unit circle, compactness gives
\[
 \min_{|z|=1}|A^{(4)}(z)|>0,
\]
and hence this factor contributes only a fixed exponential pointwise loss. The product of finitely many such factors retains the same SSV type after increasing $c_1$ and summing the corresponding covering counts. This is the product argument recorded in \cite[Proposition 4.1]{BLV}; the freedom to impose $c_3/c_2\geq\Lambda$ at the expense of increasing $c_1$ is part of \cite[Definition 1.3]{BLV}.

Finally, an origin-symmetric reflection does not change the sublevel sets. Indeed, by the identity preceding Section \ref{subsec:cylindertree},
\[
 \prod_{k=0}^{m-1}\phi_{A^*}'(L^ku)
 =e^{2\pi i\gamma_m(u)}
 \overline{\prod_{k=0}^{m-1}\phi_A'(L^ku)}
\]
for a real phase function $\gamma_m(u)$. Hence the two products have the same absolute value.
\end{proof}

The proportionality hypothesis of $c_3$ to $c_2$ will become critical in our explicit calculation of exponent values (and, indeed, this property was leveraged imprecisely in \cite{BLV} and also the author's previous work with I. {\L}aba in \cite{LM}). The exact optimal choice of $c_1$ to produce a sharp proportion $c_3/c_2$ is recorded in \eqref{eq:exactSSVmargin}, while the corresponding pointwise cost is incorporated into the total analytic loss parameter $D$, which was introduced in Proposition \ref{prop:keyproposition}.

\subsection{Sets of Small Values on coordinate projections}

We now transfer the preceding one-dimensional statement to the product $\phi_t'$ associated to a direction in $\mathbb R^d$.  On a fixed projective chart, one coordinate is normalized and the remaining $d-1$ coordinates are described by slope parameters.  The $i$th one-dimensional Riesz product is consequently evaluated at $q_i\xi$, where $q_i$ is either the normalized value $1$ or one of the slope parameters.  The multidimensional bad set must therefore be formed by pulling back each coordinate SSV set under the map $\xi\mapsto q_i\xi$, and then taking their union.

For every factor $Q$ with $Q(1)\neq0$, define its normalized trigonometric
polynomial by
\[
 \phi_Q(u):=\frac{Q(e^{2\pi iu})}{Q(1)}.
\]
In particular, and so as to clarify our notation, we have
\[
 \phi_{A_i}'(u):= \phi_{A_i'} (u) :=\frac{A_i'(e^{2\pi iu})}{A_i'(1)},
 \qquad
 \phi_{A_i}''(u):= \phi_{A_i''} (u) :=\frac{A_i'' (e^{2\pi iu})}{A_i''(1)},
\]
so that
\[
 \phi_{A_i}=\phi_{A_i}'\phi_{A_i}'',
 \qquad
 \phi_{A_i}'(0)=\phi_{A_i}''(0)=1.
\]
On a fixed signed projective chart (where we have relabeled, so ordering is lexicographic), we then set
\[
 q_1:=1,
 \qquad
 q_i:=t_{i-1}\quad(2\leq i\leq d),
\]
where the fixed chart signs have already been absorbed by reflection of the associated digit sets. We then define the single-digit Riesz product adapted to this chart as
\[
 R_{i,m}(u):=\prod_{k=0}^{m-1}\phi_{A_i}'(L^ku),
\]
so that, along this chart, our larger Riesz product satisfies:
$$
\bigg\lvert \prod_{k=0}^{m-1} \phi_t' (L^k u) \bigg\rvert : = \bigg\lvert \prod_{i=1}^{d} R_{i,m} (q_i u)\bigg\rvert.
$$
Associated to each coordinate $i = 1,...,d$, we then choose a coordinate threshold $\psi_i(m)$, whose optimal value will be determined at a later stage.  The corresponding sub-level set associated to the $i$-th coordinate, which needs to be covered appropriately, is then
\[
 B_i:=\{u\in[0,1]:|R_{i,m}(u)|\leq\psi_i(m)\}.
\]
We then define a candidate Set of Small Values associated to our total product function $\phi_t'$ in direction $t = (t_1,...,t_{d-1})$ as the following union
\begin{equation}\label{eq:productSSVunion}
 \operatorname{SSV}(t)
 :=\bigcup_{i:q_i>0}\big\{\xi\in[-c_*,c_*]:q_i|\xi|\in B_i\big\}
 =\bigcup_{i:q_i>0}\big\{\xi\in[-c_*,c_*]:|R_{i,m}(q_i\xi)|\leq\psi_i(m)\big\}.
\end{equation}

Let us examine what lower bounds hold for $\phi_t'$ over this candidate Set of Small Values chosen in \eqref{eq:productSSVunion}. First observe that, should one have that $q_i=0$, then by construction $R_{i,m}(q_i\xi)=R_{i,m}(0)=1$, and so, consequently, the associated subset of this union would be empty, and hence it does no harm to assume that $q_i > 0$. Outside of the union in \eqref{eq:productSSVunion}, every single-digit product $R_{i,m}$ is
larger than its own threshold function $\psi_i$.  Consequently, one has the uniform lower bound
\begin{equation}\label{eq:productSSVlowerbound}
 \left|\prod_{i=1}^dR_{i,m}(q_i\xi)\right|
 \geq\prod_{i=1}^d\psi_i(m),
 \qquad \xi\in [-c_*, c_*] \setminus \operatorname{SSV}(t),
\end{equation}
But, observe, that this is exactly the product principle required by Proposition \ref{prop:keyproposition}.$(2)$: the collection of one-dimensional thresholds gives a lower bound for the entire low-frequency product associated to $\phi_t'$. In particular, by choosing (say) $\psi_i (m) : = L^{-c_{1,i} m}$, the bound \eqref{eq:productSSVlowerbound} then implies that
$$
\left|\prod_{i=1}^dR_{i,m}(q_i\xi)\right|^2 : = \big\lvert P_{2,t}' (\xi)\big\rvert^2 \geq L^{-(2 \sum_{i} c_{1,i})m}
$$
as required by Proposition \ref{prop:keyproposition}.

We remark that the union \eqref{eq:productSSVunion} is essential; since, outside of this union, every active coordinate satisfies its own pointwise lower bound. In particular, if one were to attempt to construct a single threshold function, independent of each coordinate's threshold functions, then several factors could each be moderately small, causing the full product to fall below the desired bound even though no single factor crosses that common threshold. Thus, the set constructed via union in \eqref{eq:productSSVunion} is, in fact, the optimal Set of Small Values for the function $P_{2,t}'$.

\subsection{Quantitative SSV estimates and special cases}

We now quantify the amount of high-frequency mass which can lie above one coordinate SSV set, with the goal of determining optimal choices for the SSV constants $c_{1,i}, c_{2,i}$ and $c_{3,i}$ associated to the single-digit factors $\phi_{A_i}'$.  The calculation has two inputs: the SSV covering of the $i$-th coordinate, and orthogonality for the high-frequency products in complementary coordinates.

We record the result in the weakest allowable formulation which still allows for the necessary deletion estimate in Proposition \ref{prop:keyproposition}.
Let
\[
 b_{-i}:=\log_L\prod_{j\neq i}\#A_j,
\]
and assume that $ 0 < b_{-i} < 1 $ (which is allowed, since we can assume we are working with a $d$-dimensional purely $1$-unrectifiable rational product Cantor set, which does not concentrate on any lower dimensional sub-flat). The exponent $b_{-i}$ thus records the total symbolic growth contributed by all coordinates other than $i$. Indeed, if $a_j:=\log_L\#A_j$, then
\[
 \sum_{j=1}^d a_j=1,
 \qquad
 b_{-i}=1-a_i.
\]
We then use this notation in the following averaging argument for the SSV property.

\begin{proposition}[Coordinate SSV averaging]\label{prop:coordinateSSVaverage}
Fix some coordinate $i \in \{1,...,d\}$, and suppose that $c_{3,i}>1$ and that $B_i$ is
covered by $O(L^{c_{2,i}m})$ intervals of length $L^{-c_{3,i}m}$.  Choose
$E^* \subseteq [0,1]^{d-1}$ to be any measurable subset of chart parameters.  Then, if we let
\[
 Q_{j,t}(\xi):=\prod_{k=m}^{n-1}\phi_{A_j}(L^kq_j\xi),
\]
and define the following average value 
\[
 \mathcal I_i:=\frac1{|E^*|}
 \int_{[0,1]^{d-1}}
 \int_{[L^{-m},c_*]\cap\{q_i\xi\in B_i\}}
 \prod_{j=1}^d|Q_{j,t}(\xi)|^2\,d\xi\,dt,
\]
then we have the following coordinate-wise upper bound
\begin{equation}\label{eq:correctedSSVaverage}
 \mathcal I_i
 \lesssim
 m|E^*|^{-1}L^{-n}
 L^{[c_{2,i}-b_{-i}(c_{3,i}-1)]m}.
\end{equation}
Moreover, if $q_i = 1$ is the normalized coordinate, then one may omit the factor of $m$ in the upper bound. Regardless, in every case the implicit constant is uniform over the finite
signed chart family.
\end{proposition}

\begin{proof}
We divide the proof into three steps.  First, we compute the mean square of a single high-frequency coordinate product on one period.  Second, we use the supplied covering intervals (which are of small length) to improve the corresponding estimate in the identified bad coordinate.  Finally, we average in the remaining slope parameters.  As we shall see, the normalized coordinate and a slope coordinate require slightly different changes of variables, which accounts for the possible logarithmic loss factor $m \approx \log N$ in \eqref{eq:correctedSSVaverage}.

Write $a_j:=\log_L\#A_j$, so that $\sum_ja_j=1$ and
$b_{-i}=1-a_i$.  The one-coordinate high product
\[
 \widetilde Q_j(u):=\prod_{k=m}^{n-1}\phi_{A_j}(L^ku)
\]
is exactly $L^{-m}$-periodic.  Writing $q=n-m$ and expanding in base $L$ thus gives the simple estimate
\[
 \widetilde Q_j(u)
 = (\#A_j)^{-q}
 \sum_{\omega\in A_j^q}e^{2\pi iL^m\lambda_\omega u},
\]
where the integers $\lambda_\omega$ are, by definition, distinct, and we thus have that
\[
 \int_0^{L^{-m}}|\widetilde Q_j(u)|^2\,du
 =L^{-m}(\#A_j)^{-2q}\#A_j^q
 =L^{-m}(\#A_j)^{-q}.
\]
Since $\widetilde Q_j$ is exactly $L^{-m}$-periodic, summing this identity over complete periods, together with at most two incomplete endpoint periods, gives the following estimate for every interval $J$ of length at least $L^{-m}$:
\begin{equation}\label{eq:periodiccoordinateaverage}
 \int_J|\widetilde Q_j(u)|^2\,du
 \lesssim |J|(\#A_j)^{-(n-m)}.
\end{equation}

We next exploit the fact that an SSV covering interval has a fixed length factor which is much shorter than $L^{-m}$, to obtain an upper bound on the factor $Q_{i}$ on each SSV interval which is proportionally weaker, but still valid.

Let $I$ be one such interval supplied by the SSV property and put
$\ell=\lfloor c_{3,i}m\rfloor+O(1)$, with the harmless $O(1)$ chosen so that
$|I|\lesssim L^{-\ell}$.  We then have that
\begin{equation}\label{eq:badcoordinatesplit}
 \int_I|\widetilde Q_i(u)|^2\,du
   \leq \int_{I} \big\lvert \prod_{k = \ell}^{n-1} \phi_{A_i} (L^k u) \big\rvert^2 \lesssim L^{-\ell}(\#A_i)^{-(n-\ell)},
\end{equation}

We first suppose that $q_i = 1$ (so that our $i$-th coordinate is the normalized coordinate).  Then, for every coordinate $j \neq i$, we use the change-of-variables $u_j=t_j\xi$ together with the upper bound \eqref{eq:periodiccoordinateaverage} to obtain
\[
 \int_0^1|\widetilde Q_j(t_j\xi)|^2\,dt_j
 =\frac1\xi\int_0^\xi|\widetilde Q_j(u_j)|^2\,du_j
 \lesssim(\#A_j)^{-(n-m)}, \qquad \forall \xi \in [L^{-m}, c_*],
\]
where we have used that $\xi \geq L^{-m}$, together with the periodicity argument implicit in \eqref{eq:periodiccoordinateaverage} to cancel the additional Jacobian factor $\xi^{-1}$.  We may now apply \eqref{eq:badcoordinatesplit} on each of the $O(L^{c_{2,i}m})$ covering intervals, and thus obtain
\begin{align*}
 \mathcal I_i
 &\lesssim |E^*|^{-1}L^{c_{2,i}m}L^{-\ell}
 (\#A_i)^{-(n-\ell)}
 \prod_{j\ne i}(\#A_j)^{-(n-m)}\\
 &\lesssim |E^*|^{-1}L^{-n}
 L^{[c_{2,i}-b_{-i}(c_{3,i}-1)]m}.
\end{align*}
This resolves the case when the $i$-th coordinate is normalized.

Suppose that, instead, we have that $q_i \neq 1$.  In this case, we then set $u=t_i\xi$ and reverse the order of integration in the variables $u$ and $\xi$. Before doing so, however, we observe that the bad set remains uniformly separated from the origin: there is a fixed constant $c>0$ for which
\begin{equation}\label{eq:ssv-away-from-zero}
 B_i\cap[0,cL^{-m}]=\varnothing.
\end{equation}
To verify this, observe that the same style of Taylor expansion and quadratic loss estimate as in \eqref{eq:digitsetnearzero}
then guarantees that, for every  $0\leq u\leq cL^{-m}$, that one has the lower bound
\[
 |R_{i,m}(u)|
 \geq\exp\left(-Cu^2\sum_{k=0}^{m-1}L^{2k}\right)
 \geq\exp(-C'c^2) = O(1).
\]
Hence, and regardless of the threshold constant $c > 0$, so long as $m$ is chosen large enough, the SSV threshold function $\psi$ will be dominated by this lower bound, thus proving \eqref{eq:ssv-away-from-zero}.

Hence, after averaging the remaining slope coordinates and using the general upper bound \eqref{eq:periodiccoordinateaverage}, the only coupled variables remaining are $t_i$ and $\xi$.  The change-of-variables $u=t_i\xi$, for which $dt_i=du/\xi$ and $0\leq u\leq\xi$, then bounds their contribution by some fixed multiple of the quantity
\begin{align*}
\prod_{\substack{j\neq i\\j\text{ a slope coordinate}}}
 (\#A_j)^{-(n-m)}                  \times
 \int_{B_i\cap[0,c_*]}|\widetilde Q_i(u)|^2
 \left(\int_{\max\{L^{-m},u\}}^{c_*}
 |\widetilde Q_1(\xi)|^2\frac{d\xi}{\xi}\right)du.
 \end{align*}
Hence, to finish the calculation, simply partition the interval $[cL^{-m},c_*]$ into the $O(m)$-many annular neighbourhoods
\[
 [L^{-r-1},L^{-r}]\cap[cL^{-m},c_*],
 \qquad 0\leq r\leq m+O(1).
\]
On each annulus, $\xi^{-1}\lesssim L^r$, while
\eqref{eq:periodiccoordinateaverage} gives
\[
 \int_{L^{-r-1}}^{L^{-r}}|\widetilde Q_1(\xi)|^2\,d\xi
 \lesssim L^{-r}(\#A_1)^{-(n-m)}.
\]
As in the previous case, and on each such annular region, the Jacobian factor $\xi^{-1}$ is  canceled (up to constant factors) by the annulus length.  Summing over the $O(m)$ annuli gives, uniformly in the variable $u$, the estimate
\begin{equation}\label{eq:logarithmic-coordinate-average}
 \int_{\max\{L^{-m},u\}}^{c_*}
 |\widetilde Q_1(\xi)|^2\frac{d\xi}{\xi}
 \lesssim m(\#A_1)^{-(n-m)}.
\end{equation}

To conclude the calculation, we then sum
\eqref{eq:badcoordinatesplit} over the $O(L^{c_{2,i}m})$-many covering intervals, which then gives
\begin{align*}
 \mathcal I_i
 &\lesssim m|E^*|^{-1}L^{c_{2,i}m}L^{-\ell}
 (\#A_i)^{-(n-\ell)}
 \prod_{j\neq i}(\#A_j)^{-(n-m)}\\
 &\lesssim m|E^*|^{-1}L^{-n}
 L^{[c_{2,i}-b_{-i}(c_{3,i}-1)]m},
\end{align*}
thus proving \eqref{eq:correctedSSVaverage}.
\end{proof}

We now compare the preceding averaged error with the natural scale $KL^{-n}$ appearing in Proposition \ref{prop:keyproposition}.  Suppose $|E^*|\gtrsim K^{-\beta}$ and put
\[
 \varkappa:=\frac{\log_LK}{m}.
\]
Put
\[
 \Delta_i:=b_{-i}(c_{3,i}-1)-c_{2,i}.
\]
Equation \eqref{eq:correctedSSVaverage} gives
\begin{equation}\label{eq:ssv-margin-calculation}
 \frac{\mathcal I_i}{KL^{-n}}
 \lesssim mK^{\beta-1}L^{-\Delta_im}
 =mL^{-[\Delta_i+(1-\beta)\varkappa]m}.
\end{equation}
The right-hand side tends to zero provided
\begin{equation}\label{eq:exactSSVmargin}
 b_{-i}(c_{3,i}-1)-c_{2,i}
 >-(1-\beta)\varkappa.
\end{equation}
Thus the strict inequality in \eqref{eq:exactSSVmargin} produces an exponentially decaying factor in $m$, which in particular absorbs the harmless prefactor $m$.  No additional power of $K$ is lost at this stage; and, in fact, one can descend to the much simpler sufficient condition
\begin{equation}\label{eq:coordinatewiserestrictionSSVconstants}
 b_{-i}(c_{3,i}-1)-c_{2,i}>0,
\end{equation}
in order to obtain the necessary contradiction.

We then finish by showing that our previous calculation, which was uniform in the parametrized chart parameter $t$, can be imported into a proportional subset of our set $E_{N,K}^*$ of adapted Fourier directions from Section \ref{sec:toolbox}.

\begin{corollary}[Selection of one good direction]\label{cor:ssvgooddirection}
Assume that the hypotheses of Proposition
\ref{prop:coordinateSSVaverage} hold for every coordinate, that
$|E^*|\gtrsim K^{-\beta}$, and that the strict margin
\eqref{eq:exactSSVmargin} holds for every coordinate.  Then a subset
$E^{**}\subset E^*$ with $|E^{**}|\geq|E^*|/2$ satisfies
\[
 \int_{\operatorname{SSV}(t)\cap[L^{-m},c_*]}
 |P_{1,t}(\xi)|^2\,d\xi
 \lesssim KL^{-n},
 \qquad \forall t\in E^{**}.
\]
\end{corollary}

\begin{proof}
This is another application of Markov's inequality. For $t\in E^*$, set
\[
 F_i(t):=\int_{[L^{-m},c_*]\cap\{q_i\xi\in B_i\}}
 |P_{1,t}(\xi)|^2\,d\xi,
 \quad F(t):=\sum_{i=1}^dF_i(t).
\]
The union bound \eqref{eq:productSSVunion} gives
\[
 \int_{\operatorname{SSV}(t)\cap[L^{-m},c_*]}
 |P_{1,t}(\xi)|^2\,d\xi\leq F(t).
\]
Since $E^*$ lies in the chart parameter cube, Proposition
\ref{prop:coordinateSSVaverage} and the estimate
\eqref{eq:ssv-margin-calculation} yield, for suitably large $m$, the estimate
\[
 \frac1{|E^*|}\int_{E^*}F(t)\,dt
 \leq\sum_{i=1}^d\mathcal I_i
 \leq C_0KL^{-n}.
\]

We thus let 
\[
 \mathcal B:=\{t\in E^*:F(t)>2C_0KL^{-n}\}.
\]
Since $F$ is a non-negative function, Markov's inequality gives
\[
 |\mathcal B|
 \leq \frac{1}{2C_0KL^{-n}}\int_{E^*}F(t)\,dt
 \leq\frac{|E^*|}{2}.
\]
Set $E^{**}:=E^*\setminus\mathcal B$. Then
\[
 |E^{**}|=|E^*|-|\mathcal B|\geq\frac{|E^*|}{2},
\]
and, for every $t\in E^{**}$,
\[
 F(t)\leq2C_0KL^{-n}.
\]
Combining this inequality with the union bound in \eqref{eq:productSSVunion} proves the stated estimate.
\end{proof}

We finish this section by demonstrating how our exact book-keeping of the Set of Small Values parameters admits immediate and computationally simple improvements to the final power law exponent for our Favard length estimate.

\subsubsection{Zero-free mask polynomials on the complex unit circle}
Suppose each of our digits sets $A_1,...,A_d$ have mask polynomials satisfying
\[
 \mu_{A_i}:=\min_{|z|=1}|A_i(z)|>0, \textrm{ for each } i \in \{1,...,d\}.
\]
Then, as a simple consequence of our normalization of the function $\phi_{A_i}$, we must have that
\[
 |\phi_{A_i}(u)|\geq\frac{\mu_{A_i}}{\#A_i}, \qquad \forall u \in [0,1].
\]
Consequently, we have the natural lower bound
\[
 \prod_{k=0}^{m-1}|\phi_{A_i}(L^ku)|
 \geq L^{-c_im}, \textrm{ where one takes }
 c_i:=\log_L\frac{\#A_i}{\mu_{A_i}}.
\]
In particular, then, if we take
\begin{equation}\label{eq:zerofreeD}
 D=\sum_{i=1}^d\log_L\frac{\#A_i}{\mu_{A_i}},
\end{equation}
then the previous simple estimate shows that, for any $t \in [0,1]^{d-1}$, we have that
$$
\operatorname{SSV}_{\varepsilon_0}(t) : = \bigg\{\xi \in [-c_*,c_*] : \bigg\lvert \prod_{k=0}^{m-1} \phi_t' (L^k\xi)\bigg\rvert \leq L^{-(D+\varepsilon_0)m}  \bigg\} = \emptyset
$$
for any $\varepsilon_0 > 0$. Hence, one does not even need to appeal to the Set of Small Values argument to handle any sort of approximate zeroes at this threshold. The computational consequences of this observation for allowable final exponents are discussed in further detail in Section \ref{sec:exponentexamples}.

\section{Multidimensional Sets of Large Values}\label{sec:multiSLV}

We now introduce the Set of Large Values (SLV) property, which is designed to handle the terms in our Riesz product associated to the accumulating zeroes. To this end, and for each $m \gg 1$, we construct a measurable set $\Gamma_m$ whose difference set avoids the relevant zeros simultaneously at each of the first $m$ scales.  This $\Gamma_m$ is our \textit{Set of Large Values}.

The appearance of the difference set $\Gamma_m-\Gamma_m$ is not optional.  In Proposition \ref{prop:keyproposition}, observe that Salem's positive kernel function $h : [-c_*,c_*] \rightarrow (0,1]$ is supported exactly on this difference set; and, moreover, it was exactly the non-negativity of the Fourier transform $\hat{h}$ which preserves the diagonal lower bound for the high-frequency exponential sum.  We must therefore balance two competing requirements:  the set $\Gamma_m$ must be sufficiently large to produce a useful diagonal contribution; while, simultaneously, the difference set $\Gamma_m-\Gamma_m$ must be sufficiently structured to avoid every approximate zero of the accumulating products $\phi_t'' (L^k u)$ for $k = 0,1,...,m-1$.

\subsection{The Set of Large Values property and associated parameters}\label{subsec:SLVset}

We begin, as with the Set of Small Values property, by defining the Set of Large Values property relative to some underlying auxiliary parameters.

\begin{definition}\label{def:SLVset}
Let $\upsilon:\mathbb R\to\mathbb C$ be any $1$-periodic function.  We say that the function $\upsilon$ is \textbf{Set of Large Values structured} if there are constants $\eta\in(0,1]$, $C_1\geq0$, and $c_\Gamma>0$, each of which is independent of $m$, and such that for every sufficiently large integer $m$ there is a Borel set $\Gamma_m\subset[0,c_*/2]$ satisfying
\[
\Gamma_m-\Gamma_m\subset
\left\{\xi:\left|\prod_{k=0}^{m-1}\upsilon(L^k\xi)\right|\geq L^{-C_1m}\right\}
\]
while, simultaneously, one also has the measure lower bound
\[
|\Gamma_m|\geq c_\Gamma L^{-(1-\eta)m}.
\]
For the chart-dependent products below, the constants are uniform in the chart and slope parameters, even though the underlying choice of set $\Gamma_m : = \Gamma_{m,t}$ may, itself, depend upon on those parameters.
\end{definition}

As we have already seen in the Set of Small Values (see Definition \ref{def:setofsmallvalues} and its ensuing discussion), the order of quantifiers is essential.  The constants $\eta,C_1,c_\Gamma$ are chosen before the scale $m$.  The set $\Gamma_m$ may depend upon $m$ and upon the slope parameters (since our final usage of the witness function $h$ occurs at a fixed scale and along a fixed direction). However, these constants must remain uniform over every signed chart and every slope, including slopes equal or arbitrarily close to zero.

The two most critical parameters are $\eta$ and $C_1$. Here, the parameter $\eta$ is the \textit{measure gain}: it records how much larger $\Gamma_m$ is than the baseline interval scale $L^{-m}$.  The parameter $C_1$ is the \textit{pointwise loss}: it records the cost of keeping the difference set away from the accumulating zeros through all $m$ generations.  These two quantities play different roles in the final optimization, and we therefore retain them separately.

\begin{proposition}\label{prop:mainslvestimate}
Assume that every coordinate digit set satisfies at least one of the three hypotheses in Theorem \ref{thm:totalFavbound}.  Then the accumulating product
\[
\phi_t''(\xi)=\phi_{A_1}''(\xi)\phi_{A_2}''(t_1\xi)\cdots\phi_{A_d}''(t_{d-1}\xi)
\]
is SLV-structured uniformly over the finite signed chart family and all slope parameters.
\end{proposition}

\subsection{Single Scale SLV Sets}

We begin with a simple model.  Let $A=\{0,2,3,4,6\}$.  Its mask polynomial factors as $A(X) : = \Phi_5(X)\Phi_6(X)$, and, since $\#A=5$, its accumulating factor is $A''=\Phi_6$.  The zero set of $\Phi_6(e^{2\pi i\xi})$ consists of two points modulo one.  If $\Gamma_A$ is a sufficiently short periodic interval centred at the origin, then $\Gamma_A-\Gamma_A$ remains a positive distance from both zeros.  By continuity and periodicity, there is consequently a constant $c_A>0$ such that
\[
 |\phi_A''(\xi)|\geq c_A,
 \qquad \xi\in\Gamma_A-\Gamma_A.
\]
This is the basic one-scale mechanism: we choose a set whose difference set avoids the zeros of the accumulating factor, and compactness converts that avoidance into a uniform pointwise lower bound. To prove Proposition \ref{prop:mainslvestimate}, we must carry out this construction for every coordinate digit set and then combine the resulting one-scale sets through all generations. 

It is useful to separate the qualitative arithmetic problem from the quantitative multiscale problem.  At one scale, we ask whether the difference set of a reasonably large periodic set can avoid the finite zero set of $\phi_A''$ (and it is not at all obvious that such a set of reasonably large measure \textit{should even exist}).  Once such a set has been constructed for each coordinate, we intersect appropriately translated and dilated copies through the first $m$ generations.  The one-scale construction supplies the arithmetic avoidance, while the multiscale intersection determines the precise values of the exponents $\eta$ and $C_1$.

We thus proceed with the first portion of the construction: proving the existence of a relatively large and $1$-periodic set avoiding the accumulating zeroes of $A(X)$. We call such sets \textit{single scale SLV sets}.

Let $A(X)\in\mathbb Z[X]$ have the cyclotomic factorization $A=A'A''$ from Definition \ref{def:cycfactorization}.  Denote the periodic zero set of the accumulating factor by
\[
 \Sigma_A+\mathbb Z
 :=\{\xi\in\mathbb R:A''(e^{2\pi i\xi})=0\}.
\]
We call a $1$-periodic Borel set $\Gamma_A\subset\mathbb R$ a \textit{single-scale SLV set} if
\begin{equation}\label{ch5:eq:singlescaledifference}
\operatorname{dist}\bigl(\Gamma_A-\Gamma_A,\Sigma_A+\mathbb Z\bigr)=c_A>0
\end{equation}
and
\begin{equation}\label{ch5:eq:singleSLVsize}
\mathcal H^1([0,1]\cap\Gamma_A)>\frac1{\#A}.
\end{equation}
Condition \eqref{ch5:eq:singlescaledifference} supplies a positive lower bound for $|\phi_A''|$ on the difference set by compactness.  Condition \eqref{ch5:eq:singleSLVsize} says that the one-scale set has density \textit{strictly larger} than the reciprocal of the set $A$ and its relative branching number.  This strict improvement is exactly what becomes the multiscale measure gain $\eta_A>0$.

Our proof of Proposition \ref{prop:mainslvestimate} now has two stages.  First, Lemmas \ref{lma:singlescaleSLV} and \ref{lma:singlescaleSLVfibered} construct a single-scale set for every nontrivial coordinate factor allowed by Theorem \ref{thm:totalFavbound}.  Second, Lemma \ref{lma:singletomulti} intersects translated and rescaled copies of these sets to obtain one multiscale set for the full function $\phi_t''$.  All hypotheses concerning cyclotomic divisors are utilized in the first stage.

We first recall the notation attached to a nonempty digit set $A\subset\mathbb N_0$:
$$
S_{A} : = S_A^{(2)} =  \{s \in \mathbb{N} : \Phi_s (X) \mid A (X) \textrm{ and } \gcd(s,\#A) = 1 \}.
$$
and further let $s_A : = \textrm{lcm} (S_A)$. We then define
$$
\Sigma_A : = \{\xi \in [0,1] : \Phi_s (e^{2 \pi i \xi }) = 0 \textrm{ for some } s \in S_A \}.
$$
If $S_A=\varnothing$, then $A''\equiv1$ and the coordinate contributes no accumulating factor.  This includes every singleton digit set.  Such a coordinate is omitted from the SLV construction and may formally be assigned $\eta_A=1$ with zero pointwise loss.  We therefore need only consider the case $\#A\geq2$ and $S_A\neq\varnothing$.

The first single-scale result was proved by the author and I. {\L}aba in \cite{LM}; we use it as a black box.

\medskip

\begin{lemma}\label{lma:singlescaleSLV}
Suppose that $A \subset \mathbb{N}_0$ either satisfies:
\begin{enumerate}[(1)]
    \item $2 \leq \# A \leq 10$; or else,
    \medskip
    \item $s_A=p^\alpha$, or $s_A = p^{\alpha} q^{\beta}$, for distinct prime numbers $p, q$ and positive exponents $\alpha, \beta$.
\end{enumerate}
Then, there exists some $\eta_A \in (0,1)$ and a $1$-periodic set $\Gamma_A \subset \mathbb{R}$, which is a finite union of closed intervals, such that
\begin{equation}
\dist(\Gamma_A -\Gamma_A ,\Sigma_A + \mathbb{Z}) > 0
\end{equation}
\begin{equation}
\mathcal{H}^1 \big([0, 1] \cap \Gamma_A \big) > \big( \# A\big)^{-(1 - \eta_A)}
\end{equation}
\end{lemma}
The fibered alternative requires a separate single-scale statement, which we prove in Subsection \ref{subsec:SLVfiberconst}.

\medskip

\begin{lemma}\label{lma:singlescaleSLVfibered}
Let $A \subset \mathbb{N}_0$ with
$$
S = S_A^{(2)} : = \big\{s \geq 2 : \Phi_s (X) \mid A(X), \, \gcd (s,\#A) = 1\big\}
$$
and with
$$
M := \lcm (S) = \prod_{k=1}^{K} p_k^{n_k}
$$
for some distinct prime numbers $p_k$ and exponents $n_k \geq 1$. 

Suppose that $S\neq\varnothing$ and that $A$ admits a fibered subset (see Definition \ref{def:persistentfiber}). Then, there exists some $\eta_A \in (0,1)$ as well as a $1$-periodic set $\Gamma_{A} \subset \mathbb{R}$, which is a finite union of closed intervals, and such that
\begin{equation}\label{eq:singlegamma1}
\dist(\Gamma_A -\Gamma_A ,\Sigma_A + \mathbb{Z}) > 0
\end{equation}
\begin{equation}\label{eq:singlegamma2}
\mathcal{H}^1 \big([0, 1] \cap \Gamma_A \big) > \big( \#A\big)^{-(1-\eta_A)}
\end{equation}
\end{lemma}

The next lemma performs the transition from these one-scale sets to the multiscale SLV set required by Definition \ref{def:SLVset}.  A related result appeared in \cite{LM}; we include the proof because uniformity as a slope tends to zero is important in the present higher-dimensional setting.

\medskip

\begin{lemma}\label{lma:singletomulti}
Suppose that $A_1,...,A_d \subset \mathbb{N}_{0}$ are non-empty.  For every
$i$ with $A_i''\not\equiv1$, assume that there are $\eta_{A_i}\in(0,1)$ and
a $1$-periodic finite union of closed intervals $\Gamma_{A_i}$ satisfying
\[
 \dist(\Gamma_{A_i}-\Gamma_{A_i},\Sigma_{A_i}+\mathbb Z)>0
\]
and
\[
 \mathcal H^1([0,1]\cap\Gamma_{A_i})
 >(\#A_i)^{-(1-\eta_{A_i})}.
\]
Coordinates with
$A_i''\equiv1$ are omitted.  Then, uniformly for every
$t \in [0,1]^{d-1}$, the function $\phi_t''$ is SLV-structured.
\end{lemma}

Combining either single-scale lemma with Lemma \ref{lma:singletomulti} proves Proposition \ref{prop:mainslvestimate}.  Proposition \ref{prop:keyproposition} then converts the resulting SLV structure into the required Riesz-product lower bound.

\medskip

\begin{remark}
\rm{
The argument separates the remaining algebraic obstruction from the Fourier analysis.  For a more general rational product Cantor set, it is enough to construct a single-scale set satisfying the separation and density conclusions of Lemma \ref{lma:singlescaleSLV}; the multiscale and analytic portions of the proof would then apply without further change. The author is currently unaware of how to handle the general case, without an additional cyclotomic hypothesis as presented in this work.
}
\end{remark}

\subsection{The multiscale construction}\label{subsec:singletomultiSLV}
We now reintroduce the scale parameter $m$.  Starting from one single-scale set $\Gamma_{A_i}$ for each active coordinate, we construct a common set which works simultaneously for every coordinate and every generation $0\leq k<m$.  The proof has two parts.  We first average independent translations in order to retain the required measure.  We then show that the difference set of the resulting intersection remains inside the appropriate one-scale difference set after each dilation by $L^k$.

$$
S_{i} : = S_{A_i}, \quad s_i : = s_{A_i}, \quad \Sigma_i = \Sigma_{A_i}.
$$
We remind the reader that the SLV set property is defined in Definition \ref{def:SLVset}.

\begin{proof}[Proof of Lemma \ref{lma:singletomulti}]
We begin by preparing the sets at each generation.  For every coordinate with a nontrivial accumulating factor, let $\Gamma_i=\Gamma_{A_i}$ be supplied by the hypothesis.  If $A_i''\equiv1$, set $\Gamma_i=\mathbb R$,
$\eta_{A_i}=1$, and $c_i=1$; this convention imposes no restriction and
contributes no loss.  For each $k = 0,...,m-1$, let
$
\Gamma_{i,k} : = L^{-k} \,\Gamma_i.
$
The scaling is chosen so that
$$
\xi\in\Gamma_{i,k}-\Gamma_{i,k}
\quad\Longrightarrow\quad
L^k\xi\in\Gamma_i-\Gamma_i.
$$
Thus the same one-scale separation from $\Sigma_i$ controls the $k$th factor of the Riesz product.

For every nontrivial factor, the single-scale separation hypothesis produces
a constant $c_i>0$ which, uniformly in $k$, satisfies
\begin{equation}\label{eq:SLVsize}
\big\lvert \phi_{A_i}'' (L^k \xi) \big\rvert \geq c_i, \textrm{ for every } \xi \in \Gamma_{i,k} - \Gamma_{i,k}.
\end{equation}
We next choose translations of these sets so that their common intersection remains large.  Put, as a matter of notational convention,
\[
 q_1:=1,\qquad q_i:=t_{i-1}\quad(2\leq i\leq d),
 \qquad \mathcal I(t):=\{i:q_i>0\}.
\]
If $q_i=0$, then the corresponding factor is
$\phi_{A_i}''(0)=1$, and so it is omitted from the construction and considered inactive. For any active
coordinate $i\in\mathcal I(t)$, the set $q_i^{-1}\Gamma_{i,k}$ is non-empty and periodic with
period $q_i^{-1}L^{-k}$.

For each $i \in \mathcal{I}(t)$, we set
\[
 \nu_i:=\mathcal H^1([0,1]\cap\Gamma_i)>0,
\]
so that $\nu_i$ denotes the density of the set $\Gamma_i$ along a single period. Fixing a large and positive integer $R$, we observe that for each active pair $(i,k) \in \mathcal{I}(t) \times \{0,...,m-1\}$, that the following identity must hold
\begin{equation}\label{eq:slopeuniformavg}
 \frac{q_i}{R}\int_0^{R/q_i}
 \mathbf 1_{q_i^{-1}\Gamma_{i,k}}(x+\tau)\,d\tau =\frac1R\int_{q_ix}^{q_ix+R}\mathbf 1_{\Gamma_{i,k}}(u)\,du = \nu_i.      
\end{equation}
Thus, the identity \eqref{eq:slopeuniformavg} demonstrates that the average density of the sets $q_i^{-1} \Gamma_{i,k}$ over intervals of length $R_i/q_i$ is equal to the $1$-periodic density of the set $\Gamma_i$. As we shall see, then, this has the effect of removing (in the average) any error introduced slope parameter $q_i$ to the density estimate along scales $\Gamma_{i,0},...,\Gamma_{i,m-1}$.

To this end, let $I_*:=[0,c_*/2]$.  We fix the slope parameter $t$ throughout the following argument.  Thus the active index set
\[
\mathcal I(t):=\{i:q_i>0\}
\]
and the corresponding numbers $q_i$ are now fixed.  The averaging below is performed only over the translation parameters.

For every $i\in\mathcal I(t)$ and $0\leq k<m$, let
\[
T_{i,k}:=[0,R/q_i]
\]
and equip $T_{i,k}$ with the probability measure
\[
d\mu_{i,k}(\tau)
:=
\frac{q_i}{R}\mathbf 1_{[0,R/q_i]}(\tau)\,d\tau.
\]

We take each of these translation parameters as independently random variables, and hence work on the product probability space
\[
\Omega_t
:=
\prod_{i\in\mathcal I(t)}
\prod_{k=0}^{m-1}T_{i,k},
\qquad
\mu_t
:=
\bigotimes_{i\in\mathcal I(t)}
\bigotimes_{k=0}^{m-1}\mu_{i,k}.
\]
For
\[
\boldsymbol{\tau}
=
(\tau_{i,k})_{\substack{i\in\mathcal I(t)\\0\leq k<m}}
\in\Omega_t,
\]
define the translated intersection sets
\[
\Gamma(\boldsymbol{\tau})
:=
I_*\cap
\bigcap_{i\in\mathcal I(t)}
\bigcap_{k=0}^{m-1}
\bigl(q_i^{-1}\Gamma_{i,k}-\tau_{i,k}\bigr).
\]
and observe that, since the single-scale sets are finite unions of closed intervals, all translated sets and associated indicator functions appearing in this construction are measurable.  

For every $x\in I_*$, we have
\[
x\in\Gamma(\boldsymbol{\tau})
\quad\Longleftrightarrow\quad
x+\tau_{i,k}\in q_i^{-1}\Gamma_{i,k}
\]
for every $i\in\mathcal I(t)$ and every $0\leq k<m$.  Consequently,
\[
\mathbf 1_{\Gamma(\boldsymbol{\tau})}(x)
=
\mathbf 1_{I_*}(x)
\prod_{i\in\mathcal I(t)}
\prod_{k=0}^{m-1}
\mathbf 1_{q_i^{-1}\Gamma_{i,k}}
      (x+\tau_{i,k}),
\]
and it follows from this simple pointwise identity that
\[
\mathcal H^1\bigl(\Gamma(\boldsymbol{\tau})\bigr)
=
\int_{I_*}
\prod_{i\in\mathcal I(t)}
\prod_{k=0}^{m-1}
\mathbf 1_{q_i^{-1}\Gamma_{i,k}}
      (x+\tau_{i,k})\,dx.
\]

Averaging and applying Tonelli's theorem thus gives
\begin{eqnarray*}
\int_{\Omega_t}
\mathcal H^1\bigl(\Gamma(\boldsymbol{\tau})\bigr)
\,d\mu_t(\boldsymbol{\tau})
& = &
\int_{\Omega_t}
\int_{I_*}
\prod_{i\in\mathcal I(t)}
\prod_{k=0}^{m-1}
\mathbf 1_{q_i^{-1}\Gamma_{i,k}}
      (x+\tau_{i,k})
\,dx\,d\mu_t(\boldsymbol{\tau})\\
& = &
\int_{I_*}
\int_{\Omega_t}
\prod_{i\in\mathcal I(t)}
\prod_{k=0}^{m-1}
\mathbf 1_{q_i^{-1}\Gamma_{i,k}}
      (x+\tau_{i,k})
\,d\mu_t(\boldsymbol{\tau})\,dx.
\end{eqnarray*}
Because $\mu_t$ is the product measure of the measures $\mu_{i,k}$, the translation parameters are independent, and we may thus evaluate each of the inner integrations as
\begin{align*}
&\int_{\Omega_t}
\prod_{i\in\mathcal I(t)}
\prod_{k=0}^{m-1}
\mathbf 1_{q_i^{-1}\Gamma_{i,k}}
      (x+\tau_{i,k})
\,d\mu_t(\boldsymbol{\tau})\\
&\qquad=
\prod_{i\in\mathcal I(t)}
\prod_{k=0}^{m-1}
\left(
\int_{T_{i,k}}
\mathbf 1_{q_i^{-1}\Gamma_{i,k}}(x+\tau)
\,d\mu_{i,k}(\tau)
\right).
\end{align*}
For each active pair $(i,k)$, the corresponding factor is
\begin{align*}
\int_{T_{i,k}}
\mathbf 1_{q_i^{-1}\Gamma_{i,k}}(x+\tau)
\,d\mu_{i,k}(\tau)
&=
\frac{q_i}{R}
\int_0^{R/q_i}
\mathbf 1_{q_i^{-1}\Gamma_{i,k}}(x+\tau)\,d\tau\\
&=
\frac{q_i}{R}
\int_0^{R/q_i}
\mathbf 1_{\Gamma_{i,k}}
       \bigl(q_i(x+\tau)\bigr)\,d\tau.
\end{align*}
Making the suitable change of variables and applying the coordinate-wise estimate \eqref{eq:slopeuniformavg}, we thus obtain that
\begin{equation}\label{eq:averagetranslationpigeonhole}
\int_{\Omega_t}
\mathcal H^1\bigl(\Gamma(\boldsymbol{\tau})\bigr)
\,d\mu_t(\boldsymbol{\tau})
=
\int_{I_*}
\left(
\prod_{i\in\mathcal I(t)}\nu_i
\right)^m dx =
|I_*|
\left(
\prod_{i\in\mathcal I(t)}\nu_i
\right)^m.
\end{equation}
Thus, the equality \eqref{eq:averagetranslationpigeonhole} shows that the average value of the measure of our translated sets $\Gamma (\bm{\tau})$ is exactly equal to the product of all coordinate $1$-periodic densities, raised to the appropriate scale exponent $m$.

Hence, by pigeonholing \eqref{eq:averagetranslationpigeonhole}, there necessarily exists a deterministic choice of translation vector
\[
\boldsymbol{\rho}
=
(\rho_{i,k})_{\substack{i\in\mathcal I(t)\\0\leq k<m}}
\in\Omega_t
\]
such that
\[
\mathcal H^1\bigl(\Gamma(\boldsymbol{\rho})\bigr)
\geq
|I_*|
\left(
\prod_{i\in\mathcal I(t)}\nu_i
\right)^m.
\]
We hence view this choice of $\bm{\rho} : = (\rho_{i,k})$ as fixed, and then define the (localized) Set of Large Values as
\begin{equation}\label{eq:nonperiodicmultiscalegamma}
\Gamma
:=
I_*\cap
\bigcap_{i\in\mathcal I(t)}
\bigcap_{k=0}^{m-1}
\bigl(q_i^{-1}\Gamma_{i,k}-\rho_{i,k}\bigr).
\end{equation}
The chosen translations $\rho_{i,k}$, and hence the set $\Gamma$, may depend upon the fixed slope parameter $t$.  The quantitative constants obtained below, however, will be independent of $t$.

We now estimate the measure of $\Gamma$.  Whether we apply Lemma
\ref{lma:singlescaleSLV} or Lemma
\ref{lma:singlescaleSLVfibered}, the single-scale sets can be chosen with parameters $\eta_{A_i}\in(0,1)$ such that
\[
\nu_i
>
(\#A_i)^{-(1-\eta_{A_i})}
\]
whenever $A_i''\not\equiv1$; and, for any coordinate with $A_i''\equiv1$, recall that we have adopted the conventions that $\nu_i = \eta_i = 1$.
Set
\[
\eta:=\min_{1\leq i\leq d}\eta_{A_i}>0.
\]
Since $0<\nu_i\leq1$, omitting the coordinates for which $q_i=0$ can only increase the product.  It follows that
\begin{align*}
\prod_{i\in\mathcal I(t)}\nu_i
&\geq
\prod_{i=1}^d\nu_i\\
&\geq
\prod_{i=1}^d
(\#A_i)^{-(1-\eta_{A_i})}\\
&\geq
\prod_{i=1}^d
(\#A_i)^{-(1-\eta)}\\
&=
\left(
\prod_{i=1}^d\#A_i
\right)^{-(1-\eta)}\\
&=
L^{-(1-\eta)}.
\end{align*}
Combining this estimate with \eqref{eq:averagetranslationpigeonhole} and the identity $|I_*|=\frac{c_*}{2}$ gives
\begin{equation}\label{eq:slopeuniformgammasize}
\mathcal H^1(\Gamma)
\geq
\frac{c_*}{2}L^{-(1-\eta)m}.
\end{equation}
Since the constants $c_*/2$ and $\eta$ depend only upon the fixed digit sets and are independent of the slope parameter $t$, this verifies the measure theoretic properties of the set $\Gamma$.

To verify the pointwise lower bound, take any
$\xi\in\Gamma-\Gamma$ and write $\xi : = x - y$ for some $x,y \in \Gamma$. Fix an active coordinate $i\in\mathcal I(t)$ and an integer $0\leq k \leq m -1$, and observe then that
\[
x+\rho_{i,k},
\ y+\rho_{i,k}
\in q_i^{-1}\Gamma_{i,k}.
\]
Therefore,
\[
q_i(x+\rho_{i,k}),
\ q_i(y+\rho_{i,k})
\in\Gamma_{i,k}.
\]
Subtracting these two elements gives
\[
q_i\xi
=
q_i(x-y)
=
q_i(x+\rho_{i,k})
-
q_i(y+\rho_{i,k})
\in
\Gamma_{i,k}-\Gamma_{i,k}.
\]
It follows from \eqref{eq:SLVsize} that
\[
\bigl|\phi_{A_i}''(L^kq_i\xi)\bigr|
\geq c_i.
\]
If $q_i=0$, then the corresponding factor instead satisfies
\[
\phi_{A_i}''(L^kq_i\xi)
=
\phi_{A_i}''(0)
=
1.
\]
Similarly, if $A_i''\equiv1$, then the corresponding factor is identically equal to one.

For every coordinate, put
\[
\widetilde c_i:=\min\{c_i,1\},
\]
where $c_i=1$ whenever $A_i''\equiv1$.  We have therefore shown that
\[
\bigl|\phi_{A_i}''(L^kq_i\xi)\bigr|
\geq\widetilde c_i
\]
for every $1\leq i\leq d$, every $0\leq k<m$, and every
$\xi\in\Gamma-\Gamma$.  Recalling that
\[
\phi_t''(L^k\xi)
=
\prod_{i=1}^d
\phi_{A_i}''(L^kq_i\xi),
\]
we conclude that
\begin{align*}
\left|
\prod_{k=0}^{m-1}\phi_t''(L^k\xi)
\right|
=
\prod_{k=0}^{m-1}
\prod_{i=1}^d
\bigl|\phi_{A_i}''(L^kq_i\xi)\bigr|
\geq
\left(
\prod_{i=1}^d\widetilde c_i
\right)^m.
\end{align*}
Thus defining our final SLV parameter
\[
C_1
:=
-\log_L
\left(
\prod_{i=1}^d\widetilde c_i
\right)
\geq0,
\]
the preceding estimate then becomes
\[
\left|
\prod_{k=0}^{m-1}\phi_t''(L^k\xi)
\right|
\geq
L^{-C_1m},
\qquad
\xi\in\Gamma-\Gamma.
\]
The constant $C_1$ depends only upon the fixed single-scale sets and is therefore independent of $m$ and $t$.

To finish the argument, simply observe that \eqref{eq:nonperiodicmultiscalegamma} gives
\[
\Gamma\subset I_* : =[0,c_*/2], \Longrightarrow \Gamma - \Gamma \subset [-c_*, c_*].
\]
Thus the entire Salem-kernel argument remains inside the safe-frequency interval from Section \ref{subsec:countingandtrig}, and thus periodic extension of the final set is not required.  Together with \eqref{eq:slopeuniformgammasize}, the pointwise estimate verifies both conditions in Definition \ref{def:SLVset}, uniformly over all slopes, including slopes equal or arbitrarily close to zero.  This completes the proof of Lemma \ref{lma:singletomulti}.
\end{proof}

\subsection{Large-value sets arising from fibered digit sets}\label{subsec:SLVfiberconst}
We now prove the fibered single-scale lemma.  The proof uses the \textit{cluster method} developed by the author and I. {\L}aba in \cite{LM}.\label{subsec:clustermethod}  Its purpose is to divide a potentially complicated collection of cyclotomic zeros into simpler arithmetic families which can be avoided one at a time.

For one cluster $\mathcal C$, we construct a periodic neighbourhood of a lattice $Q^{-1}\mathbb Z$.  The integer $Q$ is chosen so that no zero associated to an index $s\in\mathcal C$ lies in the difference set of this neighbourhood.  We then translate and intersect the sets associated to the different clusters.  Translation averaging preserves the product of their densities, while intersection preserves every previously established zero-avoidance property.

Fix a digit set $A\subset\mathbb N_0$ throughout this subsection and enumerate its accumulating cyclotomic indices as
\[
 S_A=\{s_1,\ldots,s_J\}.
\]
We use $(m,n)$ and $\gcd(m,n)$ interchangeably for the greatest common divisor of two positive integers.

We begin with the following definition from \cite{LM}.

\medskip

\begin{definition}
A subset $\mathcal{C} \subset S_A$ is called a \textbf{cyclotomic divisor cluster}, or (for short) a \textbf{cluster}, of $A$.
\end{definition}

Fix a cluster $\mathcal C=\{s_1,\ldots,s_I\}\subset S_A$ and write
\[
 N_{\mathcal C}:=\operatorname{lcm}(s_1,\ldots,s_I)=QU.
\]
We choose the factor $Q$ so that
\begin{equation}\label{q-e1}
s_j \hbox{ does not divide } Q \hbox{ for any } j\in \{1,\dots,I\},
\end{equation}
The condition $s_j\nmid Q$ ensures that the complementary factor $U$ retains a nontrivial common divisor with every $s_j\in\mathcal C$.  Write
$$
r_j:=(s_j,Q), \quad t_j:= s_j/(s_j,Q).
$$
and put
$$
T := \max(t_1,\dots,t_I)  = \max \big(s_1/(s_1,Q),...,s_I/(s_I,Q)\big),
$$
For an auxiliary parameter $0<\rho<(QT)^{-1}$, define
\begin{equation}\label{eq:singlegammadefn}
\Gamma (\mathcal{C},\rho) =\left\{\xi\in \RR :\ \dist(\xi,\frac{1}{Q}\mathbb{Z})\leq\frac{\rho}{2}\right\}.
\end{equation}

The following lemma verifies the two properties of this lattice neighbourhood which we will need.

\medskip

\begin{lemma}\label{q-lemma1}\label{q-lemma2}
The set $\Gamma := \Gamma (\mathcal{C},\rho)$ defined above satisfies
$$
\text{dist} \big( \Gamma -  \Gamma, \Sigma(\mathcal{C}) \big) > 0
$$
where
\begin{equation}\label{q-cluster-sigma}
\Sigma(\mathcal{C}) : =  \{\xi \in \RR: \ \Phi_{s_j} (e^{2\pi i\xi})=0 \textrm{ for some } s_j\in \mathcal{C} \}.
\end{equation}
Recalling the definition of $T$, we then claim that for each $0<\lambda<T^{-1}$, there exists a choice of $0<\rho<(QT)^{-1}$ such that
\begin{equation}\label{e-lambda}
\mathcal{H}^1 \big( [0,1] \cap \Gamma  \big) > \lambda
\end{equation}
and the resulting set $\Gamma$ is 1-periodic.
\end{lemma}

\begin{proof}
Let $\xi\in\Sigma(\calC)$.  Then
\[
 \xi=\frac{b}{s_j}=\frac{b}{r_jt_j}
 \quad\text{for some }s_j\in\calC,
 \qquad (b,s_j)=1.
\]
For every $a\in\mathbb Z$,
\begin{equation}\label{eq:cluster-rational-separation}
 \left|\frac aQ-\frac b{s_j}\right|
 =\frac1{Qt_j}
 \left|at_j-\frac{bQ}{r_j}\right|.
\end{equation}
The quantity in the last absolute value is an integer because $r_j\mid Q$.
It is nonzero: equality would imply
\[
 t_j\mid\frac{bQ}{r_j},
 \qquad (b,t_j)=1,
 \qquad t_j\mid\frac Q{r_j},
\]
and hence $s_j=r_jt_j\mid Q$, contrary to \eqref{q-e1}.  Thus
\[
 \left|\frac aQ-\xi\right|
 \geq\frac1{Qt_j}\geq\frac1{QT}.
\]
Since every element of $\Gamma-\Gamma$ lies within distance $\rho$ of
$Q^{-1}\mathbb Z$,
\[
 \operatorname{dist}(\Gamma-\Gamma,\Sigma(\calC))
 \geq\frac1{QT}-\rho>0.
\]

Finally, $\rho<Q^{-1}$, so the $Q$ component intervals of $\Gamma$ in one
period are disjoint and
\begin{equation}\label{eq:cluster-density}
 \mathcal H^1([0,1]\cap\Gamma)=Q\rho.
\end{equation}
If $0<\lambda<T^{-1}$, choose
\[
 \frac{\lambda}{Q}<\rho<\frac1{QT}.
\]
Then \eqref{eq:cluster-density} gives \eqref{e-lambda}, and the definition
shows that $\Gamma$ is $1$-periodic.
\end{proof}

We also use the following result from \cite{LM}.

\medskip

\begin{lemma}\label{q-lemma3}
Suppose that $\mathcal{C}^1, ..., \mathcal{C}^k$ are clusters associated to some $A \subset \mathbb{N}_{0}$. For each $l\in\{1,\dots,k\}$, let $\Gamma^l:=\Gamma(\mathcal{C}^l,\rho^l)$ be the set defined in \eqref{eq:singlegammadefn} and satisfying \eqref{e-lambda}, with corresponding parameters $Q_l, T_l, \rho^l, \lambda_l$. Then there exist parameters $\tau_1, ..., \tau_k \in [0,1]$ such that
\begin{equation}\label{q-e4}
\mathcal{H}^1\bigg( [0,1] \cap  \bigcap_{l=1}^{k} \big( \Gamma^l+ \tau_l \big)  \bigg) > \prod_{l = 1}^{k} \lambda_l,
\end{equation}
Furthermore, if we define $\Gamma^{1,\dots,k}:=  \bigcap_{l=1}^{k} \big( \Gamma^l+ \tau_l \big)$ with this choice of $\tau_l$, then
\begin{equation}\label{q-e4a}
\text{dist} \big( \Gamma^{1,\dots,k} -  \Gamma^{1,\dots,k}, \bigcup_{l=1}^k \Sigma(\mathcal{C}^l) \big) > 0
\end{equation}
where $\Sigma(\mathcal{C}^l)$ is defined as in \eqref{q-cluster-sigma} with $
\mathcal{C}=\mathcal{C}^l$.
\end{lemma}

\begin{proof}
For arbitrary translations,
\[
 \Gamma^{1,\dots,k}-\Gamma^{1,\dots,k}
 \subset\Gamma^l-\Gamma^l
 \qquad(1\leq l\leq k).
\]
Lemma \ref{q-lemma1} therefore gives
\[
 \operatorname{dist}\left(
 \Gamma^{1,\dots,k}-\Gamma^{1,\dots,k},
 \bigcup_{l=1}^k\Sigma(\mathcal C^l)\right)
 \geq\min_{1\leq l\leq k}
 \operatorname{dist}(\Gamma^l-\Gamma^l,\Sigma(\mathcal C^l))>0.
\]

For $x\in[0,1)$, periodicity and Fubini give
\begin{align*}
 \Psi(x)
 &:=\int_{[0,1]^k}\prod_{l=1}^k
 \mathbf1_{\Gamma^l}(x-\tau_l)\,d\tau_1\cdots d\tau_k\\
 &=\prod_{l=1}^k\int_0^1
 \mathbf1_{\Gamma^l}(x-\tau_l)\,d\tau_l\\
 &=\prod_{l=1}^k\mathcal H^1([0,1]\cap\Gamma^l)
 >\prod_{l=1}^k\lambda_l.
\end{align*}
Integrating in $x$ and reversing the order of integration yields
\[
 \int_{[0,1]^k}
 \mathcal H^1\left([0,1]\cap\bigcap_{l=1}^k
 (\Gamma^l+\tau_l)\right)d\tau
 >\prod_{l=1}^k\lambda_l.
\]
At least one translation vector satisfies \eqref{q-e4}.
\end{proof}

We now apply the cluster construction to a fibered subset and prove Lemma \ref{lma:singlescaleSLVfibered}.  The only external cardinality input is Proposition \ref{prop:multifibered}, proved in \cite{KLMS}; in particular, the argument below is unconditional.

\medskip

\begin{proposition}\label{prop:singlegammafibered}
Let $A \subset \mathbb{N}_0$ with $S_A^{(2)} = S$ and
$M = \lcm (S) = \prod_{k=1}^K p_k^{n_k}$.  Fix a nonempty subset
$A'\subseteq A$ and an assignment function
$\sigma:S\rightarrow\{1,\ldots,K\}$ for which $A'$ is
$(S,\sigma)$-fibered.  Write
$$
\mathsf{FIB}(S,\sigma) := p_1^{E_1(S,\sigma)} \cdots p_K^{E_K(S,\sigma)},
$$
where recall that
$$
\capexp_k(S, \sigma) : = \{ \alpha\in \mathbb{N} : \exists \, s \in S \textrm{ with } (s, p_k^{n_k}) = p_k^{\alpha} \textrm{ and } \sigma (s) = k \}.
$$
Then for any list of parameters $(\lambda_1,...,\lambda_K)$ with $0 < \lambda_k < p_k^{-1}$, there exists a 1-periodic set $\Gamma_A \subset \mathbb{R}$ satisfying
\begin{equation}\label{eq:distancegammafiber}
\dist(\Gamma_A -\Gamma_A ,\Sigma_A)>0,
\end{equation}
\begin{equation}\label{eq:sizegammafiber}
\mathcal{H}^1 \big([0, 1] \cap \Gamma_A \big) > \bigg( \lambda_{1}^{E_1(S,\sigma)} \cdot \lambda_{2}^{E_2(S,\sigma)} \cdots \lambda_{K}^{E_K(S,\sigma)} \bigg)
\end{equation}
In particular, since $s \in S_A^{(2)}$, we know that $\gcd (\# A, p_k) = 1$ for each $k = 1,...,K$. By the above calculation, there exists a choice of $(\lambda_1,...,\lambda_K)$ and some $\eta_A \in (0,1)$ such that
\begin{equation}\label{ch5:eq:singlescaleta-2}
\mathcal{H}^1 \big([0, 1] \cap \Gamma_A \big) > \lvert A \rvert^{-(1 - \eta_A)}.
\end{equation}
\end{proposition}

To simplify notation, we let
\[
 F:=\mathsf{FIB}(S,\sigma)=\prod_{k=1}^Kp_k^{E_k(S,\sigma)}.
\]
Before proving the proposition, we verify that its density bound is strictly larger than $(\#A)^{-1}$.  Proposition \ref{prop:multifibered} gives
\[
 \#A\geq\#A'\geq F.
\]
The first and last quantities cannot be equal.  Indeed, since $S\neq\varnothing$,
there is a pair $(k,\alpha)$ with $E_k(S,\sigma)>0$ and an $s\in S$ assigned
to $k$.  Then
\[
 p_k\mid s,
 \qquad p_k\mid F.
\]
If $\#A=F$, this would give $p_k\mid(s,\#A)$, which is a contradiction of the definition of
$S=S_A^{(2)}$.  Hence $\# A > F$, and by taking $\lambda_k \rightarrow (p_k^{-1})^-$, we thus obtain that
\[
\Lambda = \Lambda (\lambda_1,...,\lambda_K) : =\prod_{k=1}^K\lambda_k^{E_k(S,\sigma)}\longrightarrow F^{-1}>(\#A)^{-1}.
\]
Hence, we may choose $\lambda_k$ so that the product, denoted by $\Lambda (\lambda_1,...,\lambda_K)$, satisfies
$(\#A)^{-1}<\Lambda(\lambda_1,...,\lambda_K)<F^{-1}$. Under this assumption, we thus obtain that
\[
 \eta_A:= \big(1+\log_{\#A} \Lambda \big)\in(0,1),
 \textrm{ where, by definition, }
 \Lambda=(\#A)^{-(1-\eta_A)}.
\]
This proves that a choice of density parameters sufficiently close to their endpoints yields a genuine exponent $\eta_A>0$, thus producing the measure gain required by Lemma \ref{lma:singlescaleSLVfibered}.

\begin{proof}[Proof of Proposition \ref{prop:singlegammafibered}]
We partition the cyclotomic indices according to the prime direction assigned by $\sigma$ and the corresponding prime-power valuation.  For $1\leq k\leq K$ and $\alpha\in\capexp_k(S,\sigma)$, put
\[
 \mathcal C^{k,\alpha}
 :=\{s\in S:\sigma(s)=k,\ v_{p_k}(s)=\alpha\}.
\]
These nonempty clusters partition $S$, and therefore
\begin{equation}\label{eq:fiber-cluster-partition}
 S=\bigcup_{k=1}^K\bigcup_{\alpha\in\capexp_k(S,\sigma)}
 \mathcal C^{k,\alpha},
 \qquad
 \Sigma_A=\bigcup_{k=1}^K\bigcup_{\alpha\in\capexp_k(S,\sigma)}
 \Sigma(\mathcal C^{k,\alpha}).
\end{equation}

Fix one pair $(k,\alpha)$ and write
\[
 \mathcal C^{k,\alpha}=\{p_k^\alpha q_1,\ldots,p_k^\alpha q_I\},
 \qquad (p_k,q_j)=1.
\]
Because $\sigma(s)=k$ implies $p_k\mid s$, one has $\alpha\geq1$.  Set
\[
 Q_{k,\alpha}:=p_k^{\alpha-1}\operatorname{lcm}(q_1,\ldots,q_I).
\]
Then $p_k^\alpha q_j\nmid Q_{k,\alpha}$ and
\[
 r_j=(p_k^\alpha q_j,Q_{k,\alpha})=p_k^{\alpha-1}q_j,
 \qquad
 t_j=\frac{p_k^\alpha q_j}{r_j}=p_k,
 \qquad T=p_k.
\]
The calculation $T=p_k$ is the reason for this particular clustering: Lemma \ref{q-lemma1} therefore supplies, for each
$0<\lambda_k<p_k^{-1}$, a $1$-periodic set $\Gamma_{k,\alpha}$ such that
\begin{equation}\label{eq:fiber-cluster-gamma}
 |[0,1]\cap\Gamma_{k,\alpha}|>\lambda_k,
 \qquad
 \operatorname{dist}(\Gamma_{k,\alpha}-\Gamma_{k,\alpha},
 \Sigma(\mathcal C^{k,\alpha}))>0.
\end{equation}

We now apply Lemma \ref{q-lemma3} simultaneously to all clusters and define
\[
 \Gamma_A:=\bigcap_{k=1}^K
 \bigcap_{\alpha\in\capexp_k(S,\sigma)}
 (\Gamma_{k,\alpha}+\tau_{k,\alpha}).
\]
Since $\#\capexp_k(S,\sigma)=E_k(S,\sigma)$,
\[
 |[0,1]\cap\Gamma_A|
 >\prod_{k=1}^K\lambda_k^{E_k(S,\sigma)}.
\]
Equations \eqref{eq:fiber-cluster-partition} and
\eqref{eq:fiber-cluster-gamma} also give
\[
 \operatorname{dist}(\Gamma_A-\Gamma_A,\Sigma_A)>0.
\]
These are exactly \eqref{eq:sizegammafiber} and \eqref{eq:distancegammafiber}.  Together with the preceding density calculation, they complete the proof of Lemma \ref{lma:singlescaleSLVfibered}.
\end{proof}

This concludes our verification that the Set of Large Values property is satisfied by the function $\phi_t''$, under the suitable cyclotomic hypotheses identified by our main results.

\section{The Improved Combinatorial Argument}\label{sec:comblemmata}
\raggedbottom

We now prove the two combinatorial statements used earlier in the Fourier argument: the low-multiplicity $L^2$ estimate (which is Lemma \ref{lma:L2exceptionalset}) and the propagation inequality (which is Lemma \ref{lma:reverseholder}).  Both arguments are most naturally formulated on the symbolic cylinder tree, which was introduced in Section \ref{subsec:cylindertree}.  This language removes the ambiguity created by set intersections in Euclidean space, since every cylinder has a unique symbolic parent, and any pair of incomparable cylinders have disjoint descendant families as collections of words, even though their projected intervals may overlap.

\subsection{An outline of the two combinatorial arguments}\label{subsec:combinatorialoutline} We continue to use the cylinder-tree notation introduced in Section
\ref{subsec:cylindertree}.  Fix a direction $\theta\in\mathbb S^{d-1}$ and
recall that
\[
    f_N^*(x):=\max_{1\leq n\leq N}f_{n,\theta}(x)
\]
records the largest projected-cylinder multiplicity attained at $x$ during
the first $N$ generations.  For every $A\geq1$, put
\[
    G_A=G_{N,A,\theta}
    :=\{x\in\mathbb R:f_N^*(x)\geq A\},
    \qquad
    \mu_\theta(A):=\mathcal H^1(G_{N,A,\theta}).
\]
When the direction is fixed, we suppress the subscript $\theta$ and write
simply $\mu(A)$.

Recall also that
\[
    E_{N,K}
    =
    \left\{
    \theta\in\mathbb S^{d-1}:
    \mu_\theta(K)\leq K^{-\rho}
    \right\},
    \qquad \rho>2.
\]
Thus the proof naturally separates into the alternatives
\[
    \theta\in E_{N,K}
    \quad\Longleftrightarrow\quad
    \mu(K)\leq K^{-\rho},
\]
and
\[
    \theta\notin E_{N,K}
    \quad\Longleftrightarrow\quad
    \mu(K)>K^{-\rho}.
\]
We retain all definitions, terminology and notation as in Section \ref{subsec:cylindertree}. For us, it is especially important to recall that two words are \textit{incomparable} if neither is a prefix of the other; and, that a family of cylinders is called
\textit{prefix-free} when its corresponding words are pairwise incomparable.  We call the depth-$N$ cylinders, or equivalently the
words of length $N$, the \textit{terminal cylinders} or \textit{terminal
words}.

As this section has two logically distinct outputs, associated to distinct but co-existing combinatorial arguments, which provide the
foundational inputs to the proofs of Lemma \ref{lma:L2exceptionalset} and
Lemma \ref{lma:reverseholder}. We provide a sketch of the arguments below.

\begin{enumerate}[leftmargin=*,labelindent=0pt]
    \item \underline{\textit{Selection of Maximal Prefix-Free Cylinder
    Families}}. Suppose first that the underlying direction belongs to
    $E_{N,K}$.  Equivalently,
    \[
        \mu(K)=\mathcal H^1(G_K)\leq K^{-\rho}.
    \]
    For every height $A\geq2$, we mark the cylinders appearing at the first
    depth at which the projected multiplicity reaches $2A$, and then retain
    only the maximal marked cylinders in the prefix ordering.  The resulting
    family $\mathscr M_A$ is prefix-free and satisfies
    \[
        \sum_{R\in\mathscr M_A}\ell(R)\lesssim A\mu(A).
    \]
    Extending every retained cylinder to the common terminal depth $N$
    produces a sub-counting function $g_{A,N}$.  Since the terminal descendant
    families of distinct members of $\mathscr M_A$ are disjoint as collections
    of words,
    \[
        \|g_{A,N}\|_1
        \approx_d \sum_{R\in\mathscr M_A}\ell(R)
        \lesssim A\mu(A).
    \]
    Moreover, we show that, if one obtains a large contribution of size $T$ at some intermediate scale (so that the maximal counting function satisfies $g_A^* \geq T$), then, in fact, $g_A^*$ is well-controlled by the maximal average value of $g_{A,N}$, which means that
    $g_A^*\lesssim_d\cM g_{A,N}$,
    where $\cM$ is the one-dimensional Hardy--Littlewood maximal operator.
    The weak $(1,1)$-boundedness of $\mathcal{M}$ then gives that
    \[
        \mu(T)\lesssim\frac{A}{T}\mu(A),
        \qquad 1\leq A\leq T.
    \]

    To control the complete dyadic energy, we retain the contribution of every
    maximal marked subtree and separate these contributions into dyadic
    classes.  This gives the refined estimate
    \[
        \mu(A_{\mathrm{dp}}LKM)
        \lesssim
        \mu(K)\sum_{r=0}^{R_K}
        (r+1)^2 2^r\mu(2^rM),
        \qquad
        R_K=\left\lceil\log_2(A_{\mathrm{dp}}LK)\right\rceil,
    \]
    and this refinement is the main combinatorial improvement over previous iterations of this argument, such as \cite{BV,BLV,NPV}.
    Substituting this estimate directly into
    \[
        \mathcal E_2:=\sum_{j\geq0}2^{2j}\mu(2^j)
    \]
    yields
    \[
        \mathcal E_2
        \lesssim K+K^2\mu(K)\mathcal E_2.
    \]
    The passage to the prefix-free families therefore creates no
    depth-dependent loss.  Instead, the entire high-multiplicity contribution
    returns to the original energy with coefficient
    \[
        K^2\mu(K)\leq K^{2-\rho}.
    \]
    If $\rho>2$ and $K$ is sufficiently large, this coefficient is strictly
    smaller than the constant required for absorption.  Consequently,
    \[
        \mathcal E_2\lesssim K \Longrightarrow 
        \|f_N^*\|_2^2\lesssim K.
    \]
    This is the low-multiplicity estimate required in Lemma
    \ref{lma:L2exceptionalset}.

    \medskip

    \item \underline{\textit{Terminal Marking and Block Propagation}}.
    Suppose instead that the underlying direction parameter lies outside $E_{N,K}$, so
    that
    \[
        \mu(K)=\mathcal H^1(G_K)>K^{-\rho}.
    \]
    By the definition of $G_K$, for every $x\in G_K$, there is some depth $1 \leq n(x)\leq N$ at which at least $\lceil K\rceil$-many projected cylinders contain the point $x$.  However, since this depth necessarily
    depends upon $x$, these witness depths $n(x)$ do not initially occur at some unified scale.  The purpose of the \textit{terminal marking lemma} is to resolve this difficulty, by passing to a unified and minimal scale where one observes $\lceil K \rceil$ stacking.

    We then concatenate cylinders $\mathcal{C}(x) : = \{Q_{x,1},...,Q_{x,\lceil K \rceil}\}$ with each $Q_{x,k}$ being a cylinder at depth $n(x)$ whose projections contain the point $x$, and further let $J_x$ denote the projected image of the union over $\mathcal{C}(x)$. Since the intervals $\{J_x\}_{x \in G_K}$ form a cover of $G_K$, Vitali's lemma extracts a pairwise disjoint sub-family $\{J_{\lambda}\}$ whose dilates cover $G_K$, and whose associated concatenated cylinder subfamilies $\mathcal{C}(\lambda)$ consist of cylinders which are pairwise incomparable across subfamilies, and each concatenated cylinder has exactly $L^{N-n(\lambda)}$ descendants at depth $N$. 
    
    We then retain only those cylinders provided by Vitali's lemma and associated to our concatenated families $\mathcal{C}(\lambda)$. By extending each such retained cylinder to depth $N$, we thus produce a subfamily whose cardinality satisfies
    $\mathscr W_K\subset\mathscr Q_N$ satisfying
    \[
    \begin{aligned}
        \#\mathscr W_K
        &\gtrsim
        K\sum_\lambda L^{N-n(\lambda)}\\
        &\asymp
        KL^N\sum_\lambda |J_\lambda|\\
        &\gtrsim
        K\mu(K)L^N.
    \end{aligned}
    \]
    Moreover, the pairwise incomparability of the selected witness cylinders ensures that no two terminal words in this subfamily $\mathscr{W}_K$ are counted twice across the pairwise disjoint subfamilies $\mathcal{C}(\lambda)$. Consequently, the marked terminal words occupy a proportion
    \[
        p:=\min\left\{
        \frac{\#\mathscr W_K}{L^N},\frac12
        \right\}
        \gtrsim K\mu(K)
        >K^{-(\rho-1)}
    \]
    of the truncated tree $\mathcal{A}^N$ of cardinality $\# \mathcal{A}^N : = L^N$. 
    
    We now use self-similarity, beginning at the common terminal scale $N$, to propagate this proportionality estimate forward across scales $N, 2N,..., J_{K,\rho}N$.  Inside every cylinder at depth $(j-1)N$, the next $N$ generations form a translated and rescaled copy of the original depth-$N$ tree.  We identify in each such cylinder, up to rescaling, the same $N$-letter words determined by
    $\mathscr W_K$, declaring the corresponding retained descendants green,
    and the remaining discarded  descendants to be black.  Every black block parent therefore has, at most,
    \[
        (1-p)L^N
    \]
    black descendants after the next block.  If $\mathscr B_j$ denotes the
    black family at depth $jN$, then
    \[
        \#\mathscr B_j\leq L^{jN}(1-p)^j,
    \]
    and hence
    \[
        \mathcal H^1\left(
        \pi_\theta\left(\bigcup_{Q\in\mathscr B_j}Q\right)
        \right)
        \lesssim(1-p)^j
        \leq e^{-pj}.
    \]
    The green contribution is controlled by the multiplicity retained inside each concatenated cluster.  If
    \[
        g_j:=\sum_{Q\in\mathscr G_j}\mathbf 1_{I_Q},
    \]
    then
    \[
        \|g_j\|_1
        \approx_d \#\mathscr G_j L^{-jN}
        \lesssim_d 1
    \]
    Every point in the green projection belongs to a cluster arising from at least a constant multiple of $K$-many concatenated cylinders with a common projected point, and we thus have the following pointwise lower bound
    \[
        \cM g_j(x)\gtrsim_d K
    \]
    where $g_j$ is the counting function associated to the green family $\mathscr{G}_j$ and $\cM$ again denotes the Hardy-Littlewood maximal function.  The weak $(1,1)$ inequality for the Hardy-Littlewood maximal function therefore furnishes the estimate
    \[
    \begin{aligned}
        \mathcal H^1\left(
        \pi_\theta\left(\bigcup_{Q\in\mathscr G_j}Q\right)
        \right)
        &\leq
        \mathcal H^1\{\cM g_j\gtrsim_d K\}\\[1ex]
        &\lesssim K^{-1}\|g_j\|_1\\[1ex]
        &\lesssim_d K^{-1}.
    \end{aligned}
    \]

    Finally, take
    \[
        j=J_{K,\rho}
        \asymp K^{\rho-1}\log(2+K).
    \]
    Since $p\gtrsim K^{-(\rho-1)}$, we have
    \[
        pj\gtrsim\log(2+K),
        \qquad
        e^{-pj}\lesssim K^{-1}.
    \]
    Combining the green and black estimates gives
    \[
        \mathcal H^1\bigl(
        \pi_\theta(\mathcal S^{NJ_{K,\rho}})
        \bigr)\lesssim K^{-1}.
    \]
    This is the propagation estimate required in Lemma
    \ref{lma:reverseholder}.
\end{enumerate}

\begin{figure}[H]
  \centering
  \begin{subfigure}[t]{.56\linewidth}
    \centering
    \includegraphics[width=\linewidth]{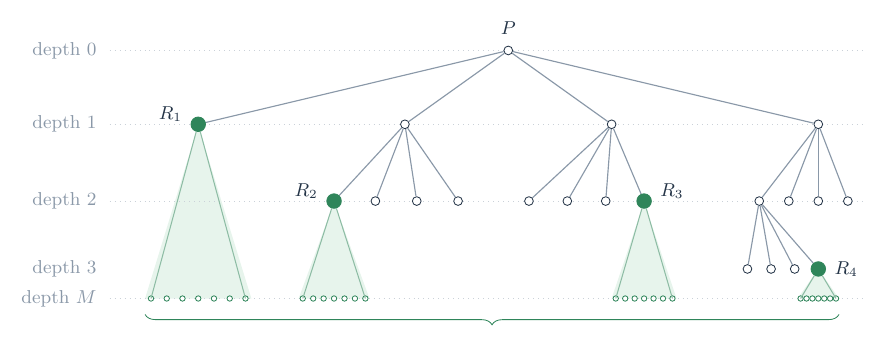}
    \caption{A prefix-free family in the symbolic tree.  The green vertices
      are pairwise incomparable, so their terminal descendant families are
      disjoint as collections of words.}
    \label{fig:symbolic-prefix-free-family}
  \end{subfigure}\hfill
  \begin{subfigure}[t]{.39\linewidth}
    \centering
    \includegraphics[width=\linewidth]{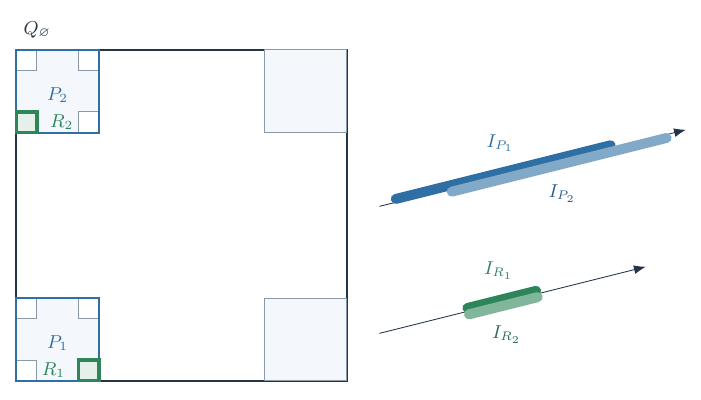}
    \caption{Two incomparable cylinders may nevertheless have coincident
      projections.  Symbolic disjointness is therefore distinct from
      geometric disjointness on the projection line.}
    \label{fig:planar-coincident-projections}
  \end{subfigure}
  \caption{The symbolic and projected viewpoints used throughout the
    combinatorial argument.}
  \label{fig:symbolic-and-projected-viewpoints}
\end{figure}

In both arguments, we retain the cylinder-tree notation introduced in Section \ref{subsec:cylindertree} and recalled in Section \ref{subsec:combinatorialoutline}. We also include Figure
\ref{fig:symbolic-prefix-free-family}, which illustrates a prefix-free family at distinct depths; whereas, the complementary Figure \ref{fig:planar-coincident-projections} gives a visual depiction as to why symbolic incomparability does not imply disjointness after projection.

Throughout both arguments, we fix a direction $\theta \in \mathbb{S}^{d-1}$, and again utilize the notation
\[
 \sigma_\theta:=\|\theta\|_{\ell^1}
 =\sum_{i=1}^d|\theta_i|, \textrm{ with } 1 \leq \sigma_{\theta} \leq \sqrt{d}.
\]
For each intermediate scale $1 \leq n \leq N$, we continue to write
$$
\mathscr{Q}_n : = \big\{Q_\omega : \omega \in \mathcal{A}^n \big\}
$$
and, since our projection parameter is fixed, we simply write $I_{\omega, \theta} = I_\omega : = \pi_{\theta} (Q_\omega)$. For each integer $A \geq 0$ and intermediate scale $1 \leq n \leq N$, we further use the notation
$$
\mu_n (A) : = \mathcal{H}^1 \big(\big\{x \in \mathbb{R} : \max_{1 \leq k \leq n} f_{k,\theta} (x) \geq A \big\} \big),
$$
so that,
$$
\mu_n (A) \leq \mu_{n+1} (A) \Rightarrow \max_{1 \leq n \leq N} \mu_n (A) = \mu (A), \quad \forall A \geq 0,
$$
where we remind the reader that $\mu (A) : = \mu_N (A)$.

All constants below may depend upon $L$, $d$, and the digit sets $A_1,...,A_d$, but not upon $N$, $K$, or $\theta$, and we regularly absorb such constants into $\lesssim$ inequalities with underlying dependencies suppressed.

\subsection{The low-multiplicity energy estimate}\label{subsec:L2exceptional}
We now give a proof of Lemma \ref{lma:L2exceptionalset}. Since this result is considered pointwise for directions $\theta \in E_{N,K}$, we therefore perform all calculations along some fixed direction $\theta_0$ for which
$$
\mathcal{H}^1 \big( \big\{x \in \ell_{\theta_0} :  f_{N,\theta_0}^* (x) \geq K  \big\}\big) \leq K^{- \rho }.
$$
This allows us to identify $\ell_{\theta_0} \cong \mathbb{R}$, and thus to perform all calculations relative to this single, fixed low-multiplicity direction.
 
We organize the argument (which is rather involved) into three general steps, which are included as the following sections.
\begin{itemize}
    \item In Section \ref{subsubsec:stoppingtime}, we construct, for some given stacking threshold $A \gg 1$, a stopping time forest $\mathscr{M}_A$. This stopping time forest identifies those prefix-free symbolic subtrees where overlaps of size $\geq A$ first appear, while simultaneously controlling their growth. As a consequence of the structure of $\mathscr{M}_A$ (in particular, its minimality among scales where large stacking first appeared, as well as the prefix-free hypothesis) we obtain a proportionality inequality between the measures $\mu(A)$ and $\mu(T)$ of the superlevel sets at depths $1 \leq A \leq T$.
    \medskip
    \item In Section \ref{subsubsec:selfsimtails}, we examine what happens inside of one such stopping time subtree $R \in \mathscr{M}_A$. In particular, we are able to include a renormalization in scale for all stopping time subtrees $R \in \mathscr{M}_A$, which proceeds essentially by rescaling and self-similarity. We also prove an elementary dyadic pigeonholing inequality, which eventually will allow us to recombine this local subtree renormalization across all such stopping time subtrees $R \in \mathscr{M}_A$.
    \medskip
    \item In Section \ref{subsubsec:twoparametertail}, we leverage all previous intermediate results to prove a powerful two-parameter tail estimate. In essence, given depth parameters $1 \ll K \ll M$, this result allows us to bound the measure $\mu(KM)$ of the superlevel set at the product scale $KM$ by the product of the measure $\mu (K)$ and the geometric sum of the measures $\mu(2^r M)$ along scales $1 \leq r \lesssim \log K$. This two-parameter inequality, together with the assumed bound $\mu(K) \lesssim K^{-\rho}$, then directly implies the $L^2$ estimate for the maximal counting function $f_{N,\theta_0}^*$ via a direct bound for the dyadic distribution function associated to $f_{N,\theta_0}^*$, which follows from the two-parameter tail estimate using the estimate on $\mu(K)$ as a base hypothesis.
\end{itemize}

\subsubsection{Stopping-time packing and relative weak-type control}\label{subsubsec:stoppingtime}
\begin{figure}[H]
  \centering
  \includegraphics[width=.98\linewidth]{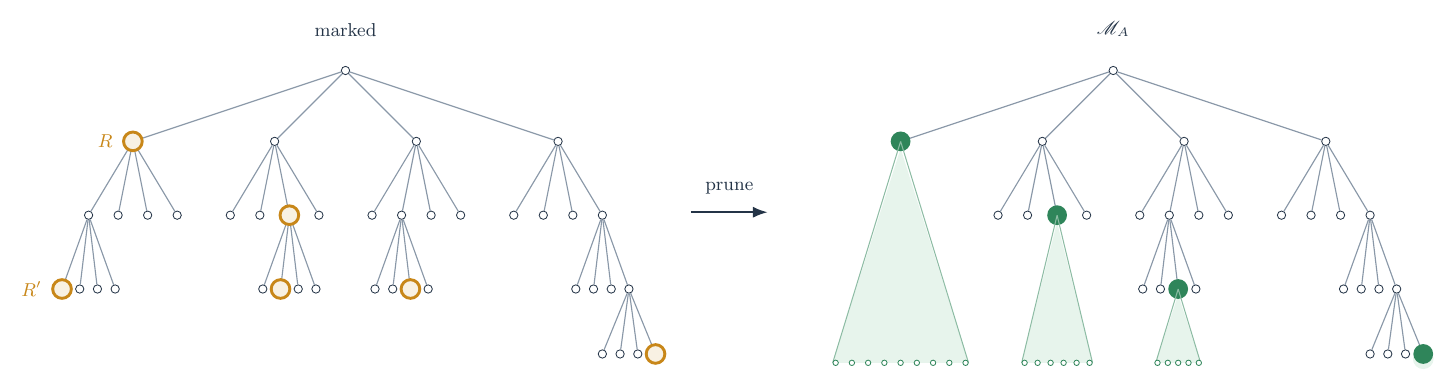}
  \caption{The first-crossing stopping-time construction.  Marks may occur
    at different depths.  Removing every marked cylinder strictly contained
    in another marked cylinder leaves the prefix-free family $\mathscr M_A$;
    the green cones represent its disjoint terminal descendant families.}
  \label{fig:first-crossing-pruning-tree}
\end{figure}

We begin with the following lemma, which allows us to identify a prefix-free family of cylinders which exhibits large stacking, while simultaneously having the corrected sum of geometric lengths. This procedure is illustrated in Figure \ref{fig:first-crossing-pruning-tree}.

\begin{lemma}[Maximal marked-cylinder packing]\label{lem:marked-packing}
For every $A\geq2$, there is a prefix-free family $\mathscr M_A$ of cylinder cubes, at possibly different depths, such that the following two conditions hold simultaneously.
\begin{enumerate}
\item The sum of geometric lengths along $\mathscr M_A$ is bounded by the corresponding measure of the good set at threshold $A$; which is to say,
\begin{equation}\label{eq:marked-packing}
\sum_{R\in\mathscr M_A}\ell(R)\lesssim A\mu(A).
\end{equation}
\item If $T\geq2LA$, $x\in G_T$, $f_n(x)\geq T$, then one has the following implication
\begin{equation}\label{eq:marked-packing-2}
x \in I_{\omega} \Longrightarrow Q_{\omega}\subseteq R
\text{ for some }R\in\mathscr{M}_A, \qquad \forall \omega \in \mathcal{A}^n.
\end{equation}
\end{enumerate}

\end{lemma}

\begin{proof}
We first construct the prefix-free family $\mathscr{M}_A$; we then estimate the total geometric length as in \eqref{eq:marked-packing} and establish the implication \eqref{eq:marked-packing-2}.

Set $f_0=\mathbf1_{I_{Q_\varnothing}}$.  For each $x\in G_{2A}$, define
$$
n_A (x) : = \min \big\{1 \leq n \leq N : f_{n} (x) \geq 2 A \big\},
$$
so that $n_A (x)$ denotes the minimal depth at which a stack of height $\geq 2 A$ is observed at $x$. We then define, relative to this parameter, the family of marked cylinders by setting
$$
 \mathscr{M}^{\text{marked}}_A
 :=
 \bigcup_{x\in G_{2A}}
 \left\{
 Q\in\mathscr Q_{n_A(x)}:x\in I_Q
 \right\}.
$$
Thus, a cylinder $Q$ at depth $m$ is a \textbf{marked cylinder} (or, simply, \textbf{marked}) if there exists some $x \in G_{2A}$ such that $m = n_A (x)$ and $x$ is contained in the projection $I_Q$ of $Q$. Since every cylinder has exactly $L$ children, $f_n(x)\leq Lf_{n-1}(x)$ for every $1 \leq n \leq N$ and for every $x \in \mathbb{R}$; hence,
\begin{equation}\label{eq:first-crossing}
 2A\leq f_{n_A(x)}(x)
 \leq Lf_{n_A(x)-1}(x)<2LA, \qquad \forall x \in \mathbb{R}.
\end{equation}
We record this inequality, as it is ubiquitous in later calculations.
 
To construct the family $\mathscr{M}_A$, we proceed as follows. For each $x \in G_{2A}$, consider all words $\omega \in \mathcal{A}^{n_A(x)}$ such that $x \in I_{\omega}$; which, in our notation, is the collection
$$
\mathscr{M}_{A,x}^* : = \big\{\omega \in \mathcal{A}^{n_A(x)} : x \in I_{\omega} \big\}.
$$
We then take 
$$
\widetilde{\mathscr{M}}^*_{A} : = \bigcup_{x \in G_{2A}} \mathscr{M}_{A,x}^*
$$
and remark that each word $\omega \in \widetilde{\mathscr{M}}_{A}$ need not live at some common depth. Finally, we generate a maximal prefix-free subfamily from $\widetilde{\mathscr{M}}_{A}^*$ by letting
$$
\mathscr{M}_{A}^* := \big\{\omega \in \widetilde{\mathscr{M}}_{A}^* : \omega' \not\preceq \omega, \textrm{ for any } \omega' \in \widetilde{\mathscr{M}}^*_{A} \setminus \{\omega \}\big\}.
$$
Hence, if $\tau \in \widetilde{\mathscr{M}}_{A}^* \setminus  \mathscr{M}_{A}^*$, then there necessarily exists some $\omega = \omega(\tau) \in \mathscr{M}_{A}^*$ such that $\omega \preceq \tau$ (and, in fact, this choice of prefix is unique).

We thus take
$$
\mathscr{M}_A : = \big\{Q_{\omega} : \omega \in \mathscr{M}_A^* \big\}
$$
so that each cylinder cube $Q_{\omega} \in \mathscr{M}_A$ corresponds to a unique prefix-free word $\omega \in \mathscr{M}_A^*$. Notice that, in particular, this implies that if $Q_{\omega}, Q_{\omega'} \in \mathscr{M}_A$ are distinct cylinder cubes, then necessarily $Q_{\omega} \not\subseteq Q_{\omega'}$ and $Q_{\omega'} \not\subseteq Q_{\omega}$.

We first prove \eqref{eq:marked-packing}.  Let $\mathcal D$ be the standard closed $L$-adic grid on $\mathbb R$; so that,
$$
\mathcal{D} : = \bigcup_{k \in \mathbb{Z}} \mathcal{D}_k, \textrm{ with } \mathcal{D}_k : = \big\{[m L^{-k}, (m+1) L^{-k}] : m \in \mathbb{Z} \big\}.
$$
We will use this complete $L$-adic grid $\mathcal{D}$ to organize marked intervals of different lengths.  To this end, we first define a sub-collection
\begin{equation}\label{eq:adic-candidate}
\mathcal{C}_A : = \big\{D\in\mathcal D:D\cap I_R\neq\varnothing
 \text{ for some }R\in\mathscr M_A\text{ with }\ell(R)\geq|D| \big\},
\end{equation}
before letting
\begin{equation}\label{eq:maximal-adic}
\mathcal D_{\max} : = \big\{D \in \mathcal{C}_A : \not\exists D' \in \mathcal{C}_A \textrm{ with }  D \subsetneq D' \big\}.
\end{equation}
Thus, by the construction of this family $\mathcal{D}_{\max}$, we observe that the interiors of all such $D \in \mathcal{D}_{\max}$ are pairwise disjoint; and, moreover, their union covers all possible intervals $I_R$ with $R \in \mathscr{M}_A$.

We thus fix $D\in\mathcal D_{\max}$ and write $|D|=L^{-k}$ for some $k \in \mathbb{N}$ (so that $k$, itself, is also fixed).  Maximality and the discreteness of the geometric lengths imply
\begin{equation}\label{eq:no-larger-marked}
 \ell(R)\leq|D|,
 \text{ whenever }R\in\mathscr M_A\text{ and }I_R\cap D\neq\varnothing.
\end{equation}
Indeed, this holds since, if $\ell(R)>|D|$, then the $L$-adic parent of $D$ would satisfy \eqref{eq:adic-candidate}. The same argument shows that some $R_D\in\mathscr M_A$ meeting $D$ satisfies $\ell(R_D)=|D|$. For every $R\in\mathscr M_A$ satisfying $I_R \cap D \neq \emptyset$, let $P_k(R)\in\mathscr Q_k$ denote its depth-$k$ ancestor. That such an ancestor exists is implied by \eqref{eq:no-larger-marked}. 

Suppose first, for contradiction, that there are strictly greater than $4LA$-many distinct cylinders $P_k(R)$. Then, each of their associated intervals, which have length $\sigma_\theta|D|\geq|D|$, must meet $D$; and in particular, one such interval must contain at least one of the two endpoints $y_-, y_+$ of $D$. In particular, 
\[
\max \{f_k(y_-), f_k(y_+)\} \geq\frac12\#\{P_k(R):I_{P_k(R)}\cap D\neq\varnothing\}>2LA,
\]
and we let $y$ denote one such maximal-stacking endpoint.
In particular, this shows that if
$$
\#\{P_k(R):I_{P_k(R)}\cap D\neq\varnothing\} > 4LA,
$$
then
\begin{equation}\label{eq:contradiction-index}
m :=\min \big\{1 \leq n \leq N : f_n (y) \geq 2 A\big\} \leq k - 1.
\end{equation}
Now, let $Q \in \mathscr{Q}_m$ satisfy $y \in I_Q$. Since we have established that $f_{m-1} (y) < 2 A \leq f_m (y)$, there necessarily exists some $R' \in \mathscr{M}_A$ such that $Q \subseteq R'$, and so, in particular, $\ell (R') \geq \ell (Q) > |D|$, since $m < k$. Moreover, we know that $y \in I_Q \subseteq I_{R'}$, and so $I_{R'} \cap D \neq \emptyset$. Hence, we have contradicted \eqref{eq:no-larger-marked}, and thus proven that
\begin{equation}\label{eq:minimal-depth-ancestor-bound}
\# \{P_k (R) : I_{P_k(R)} \cap D \neq \emptyset \} \leq 4LA.
\end{equation}

Now, continuing to keep the $L$-adic interval $D \in \mathcal{D}_{\max}$ fixed, we thus define a sub-family
$$
\mathscr{P}_D : = \big\{P_k (R) : R \in \mathscr{M}_A, I_R \cap D \neq \emptyset \} \subseteq \mathscr{Q}_k.
$$
For each $P\in \mathscr{P}_D$, the members $R \in \mathscr M_A$ contained in $P$ are prefix-free descendants; and, hence, the $L$-ary Kraft inequality (which we have stated in the form of Lemma \ref{lma:prefixfreepacking}) gives
$$
 \sum_{\substack{R\in\mathscr M_A\\R\subset P}}\ell(R)\leq\ell(P)=|D|.
$$
By combining this estimate with \eqref{eq:minimal-depth-ancestor-bound}, we thus obtain the inequality
\begin{equation}\label{eq:local-packing}
 \sum_{\substack{R\in\mathscr M_A\\I_R\cap D\neq\varnothing}}\ell(R)
 \leq4LA|D|.
\end{equation}

Now, choose a cylinder $R_D \in \mathscr{M}_A$ satisfying $I_{R_D} \cap D \neq \emptyset$ and $\ell (R_D) = |D|$. That such a cylinder must exist follows from the definitions of sub-families $\mathcal{C}_A$ and $\mathcal{D}_{\max}$ given in \eqref{eq:adic-candidate} and \eqref{eq:maximal-adic}, respectively. Then, there necessarily exists a point $x_D\in I_{R_D}$ at which at least $\lceil2A\rceil$ depth-$k$ projected cylinders overlap.  At least half of these projections extend from $x_D$ by $\sigma_\theta|D|/2$ on the same side.  Since
$\lceil2A\rceil/2\geq A$, there is a one-sided interval $J_D$ simultaneously satisfying
\begin{equation}\label{eq:charge-interval}
 |J_D|= \frac{|I_{R_D}|}{2}=\frac{\sigma_\theta}{2} \lvert D \rvert, \textrm{ and also }
 J_D\subset G_A.
\end{equation}
We thus let $C_D \subseteq G_A$ be the connected component containing $J_D$.  Since $I_{R_D}$ meets $D$, there exists a dimensional constant $c_d \geq 1$ such that
\begin{equation}\label{eq:component-enlargement}
 D\subset c_d C_D,
\end{equation}
where the right-hand side denotes the concentric enlargement by the factor $c_d$.  Indeed, the interval containment in \eqref{eq:component-enlargement} follows since
$$
 |C_D|\geq|J_D|=\frac{\sigma_\theta}{2}|D|\geq\frac{|D|}{2},
$$
while also one has
$$
\operatorname{dist}(D,C_D)\leq |I_{R_D}|\leq\sqrt d\,|D|.
$$

For each connected component $C$ of $G_A$, we thus define
$$
 \mathcal D(C):=\{D\in\mathcal D_{\max}:C_D=C\}.
$$
The intervals in $\mathcal D_{\max}$ have disjoint interiors, and
\eqref{eq:component-enlargement} places every member of $\mathcal D(C)$ inside
the same fixed enlargement $c_dC$.  Hence
$$
 \sum_{D\in\mathcal D(C)}|D|
 =\left|\bigcup_{D\in\mathcal D(C)}D\right|
 \leq c_d |C|.
$$
The connected components of $G_A$ are pairwise disjoint and their lengths sum to $\mu(A)$. Hence, summing over these connected components thus gives
\begin{equation}\label{eq:adic-packing}
 \sum_{D\in\mathcal D_{\max}}|D|
 \leq c_d 
 \sum_{C\in\operatorname{Comp}(G_A)}|C|
 =c_d \, \mu(A).
\end{equation}
As every $R\in\mathscr M_A$ meets at least one $L$-adic interval in $\mathcal D_{\max}$, we may thus sum over \eqref{eq:local-packing} while also summing over $D\in\mathcal D_{\max}$ and use \eqref{eq:adic-packing} to obtain
$$
 \sum_{R\in\mathscr M_A}\ell(R)
 \lesssim A\sum_{D\in\mathcal D_{\max}}|D|
 \lesssim A\mu(A),
$$
which is \eqref{eq:marked-packing}.

Finally, suppose that $T\geq 2LA$, $x\in G_T$, and
$f_n(x)\geq T$.  Set $m:=n_A(x)$.  Since
$$
2A\leq f_m(x)<2LA\leq T\leq f_n(x),
$$
we have $m<n$.  Let $Q=Q_\omega\in\mathscr Q_n$ be any
cylinder over $x$ at depth $n$ (so that $x\in I_Q$) and let
$
P \in \mathscr Q_m
$
denote its ancestor at depth $m$.  Since $Q\subseteq P$, we have
$I_Q\subseteq I_P$, and therefore $x\in I_P$.  By the definition
of $m=n_A(x)$, every cylinder at depth $m$ whose projection contains
$x$ is marked; and so, in particular, $P$ is a marked cylinder.  Every marked cylinder is contained in a member of $\mathscr M_A$, so there is
some $R\in\mathscr M_A$ such that
\[
Q\subseteq P\subseteq R.
\]
Thus every cylinder at depth $n$ over $x$ is contained in a member of
$\mathscr M_A$, proving the final assertion.
\end{proof}

\begin{figure}[H]
  \centering
  \begin{subfigure}[t]{.66\linewidth}
    \centering
    \includegraphics[width=\linewidth]{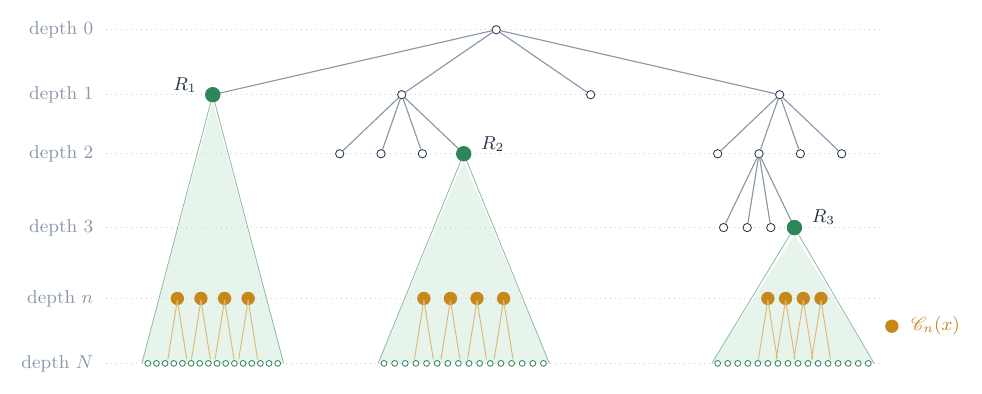}
    \caption{Each retained cylinder is extended to the common terminal depth
      $N$.  The gold vertices are the intermediate cylinders counted by
      $g_{A,n}$, and the green vertices are their terminal descendants.}
    \label{fig:terminal-extension-tree}
  \end{subfigure}\hfill
  \begin{subfigure}[t]{.28\linewidth}
    \centering
    \includegraphics[width=\linewidth]{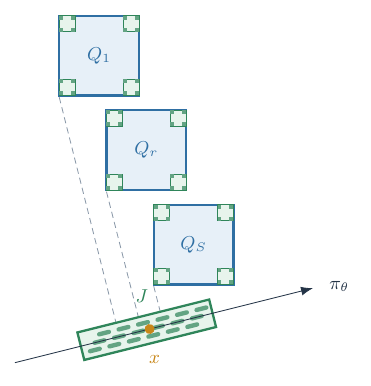}
    \caption{The averaging interval used to dominate an intermediate stack
      by terminal descendant mass.}
    \label{fig:planar-terminal-domination}
  \end{subfigure}
  \caption{Terminal extension and the associated maximal-function estimate.}
  \label{fig:terminal-extension-and-domination}
\end{figure}

We next extend each cylinder $R \in \mathscr{M}_A$ to a common terminal depth $N$.  More precisely, if
$R=Q_\tau \in \mathscr{M}_A$ and $m=\operatorname{depth}(R)$, let
\[
 \operatorname{Desc}_N(R)
 :=
 \{Q_{\tau\nu}:\nu\in\mathcal A^{N-m}\}
\]
denote the family of all depth-$N$ descendants of $R$.  This family has
cardinality $L^{N-m}$.  Moreover, the families
$\operatorname{Desc}_N(R)$ are pairwise disjoint as $R$ ranges over the
retained cylinders!  Indeed, if a terminal word $\omega$ descended from two
retained cylinders $Q_\tau$ and $Q_{\tau'}$, then $\tau$ and $\tau'$ would
both be prefixes of $\omega$.  Two prefixes of the same word are comparable,
so one of $\tau$ and $\tau'$ would be a prefix of the other. This is a contradiction of the prefix-freeness of the family $\mathscr{M}_A$, unless $\tau=\tau'$.  Thus every terminal word belongs to at most one retained subtree, and the terminal descendants may be counted
injectively. This cylinder extension process is illustrated in Figure \ref{fig:terminal-extension-tree}

For a cylinder $R$, we again let $\operatorname{depth}(R)$ denote the length of its defining word, which is distinct from its geometric length $\ell(R)$ and from the length $|I_R|$ of its projection.  For $A\geq2$ and $0\leq n\leq N$, define for each $x \in \mathbb{R}$, the modified counting and maximal counting functions
\begin{equation}\label{eq:green-counting}
 g_{A,n}(x)
 :=\sum_{\substack{R\in\mathscr M_A\\\operatorname{depth}(R)\leq n}}
   \sum_{\substack{Q\in\mathscr Q_n\\Q\subset R}}\mathbf1_{I_Q}(x),
 \qquad
 g_A^*(x):=\max_{0\leq n\leq N}g_{A,n}(x).
\end{equation}
Since $\mathscr M_A$ is prefix-free, its terminal descendant families are disjoint as collections of words, and this has the immediate impact of establishing the following inequality.
\begin{equation}\label{eq:green-L1}
 \|g_{A,N}\|_1
 =\sigma_\theta\sum_{R\in\mathscr M_A}\ell(R)
 \lesssim A\mu(A),
\end{equation}
where, in the final inequality, we applied \eqref{eq:marked-packing}. This $L^1$ estimate will allow us to bound the maximal stacking function $g_A^*$, once we have the following Lemma in hand.

\begin{lemma}[Terminal domination]\label{lem:terminal-domination}
One has the inequality
\begin{equation}\label{eq:terminal-domination}
 g_A^*\lesssim_d\cM g_{A,N}, \qquad \forall x \in \mathbb{R}
\end{equation}
where $\cM$ is the one-dimensional Hardy--Littlewood maximal operator.
\end{lemma}

We remark that the averaging interval used in the proof below is illustrated in Figure
\ref{fig:planar-terminal-domination}.

\begin{proof}
Suppose $g_{A,n}(x)\geq S$, and let $\mathscr C_n(x)$ be the family of
depth-$n$ cylinders counted by $g_{A,n}(x)$. Then $\# \mathscr{C}_n(x) \geq S$, while $x \in I_Q$ for every $Q \in \mathscr{C}_n (x)$. Every interval $I_Q$ has length $\sigma_\theta L^{-n}$; and so, consequently, if
\[
 J:=[x-\sigma_\theta L^{-n},x+\sigma_\theta L^{-n}],
\]
then $I_Q\subset J$ for every $Q\in\mathscr C_n(x)$, while
$|J|=2\sigma_\theta L^{-n}$.  

Now, the projections of the terminal descendants $\operatorname{Desc}_N(Q)$ of all such cylinders $Q \in \mathscr{C}_n (x)$ remain inside $I_Q$, while the terminal descendant families of distinct
depth-$n$ cylinders are disjoint as collections of words.  Therefore
\begin{align*}
 \int_J g_{A,N}(y)\,dy
 &\geq
 \sum_{Q\in\mathscr C_n(x)}
 \sum_{\substack{Q'\in\mathscr Q_N\\Q'\subset Q}}|I_{Q'}|\\
 &=\#\mathscr C_n(x)L^{N-n}\sigma_\theta L^{-N}\\
 &\geq S\sigma_\theta L^{-n}
 =\frac S2|J|.
\end{align*}
Dividing by $|J|$ and taking a supremum thus proves that $\cM g_{A,N}(x)\geq S/2$; and then, taking a maximum over scales
$0\leq n\leq N$ proves \eqref{eq:terminal-domination}.
\end{proof}

This establishes the sub-family $\mathscr{M}_A$ of good cylinders, a restricted maximal counting function $g_A^*$ associated to these good cylinders, as well as a corresponding boundedness estimate relating this restricted maximal counting function $g_A^*$ to the Hardy-Littlewood maximal function of the regular restricted stacking function $g_{A,N}$ at the terminal depth $N$.

We can easily convert this pointwise estimate on $g_A^*$ into a statement comparing the super-level sets $\mu(A)$ and $\mu(T)$ for $1 \leq A \leq T$. This is the content of the following result.

\begin{proposition}[Relative weak type estimate for $\mu$]\label{prop:relative-weak}
For all $1\leq A\leq T$,
\begin{equation}\label{eq:relative-weak}
 \mu(T)\lesssim\frac{A}{T}\mu(A).
\end{equation}
Equivalently, for every integer $r\geq0$,
\begin{equation}\label{eq:relative-dyadic}
 \mu(2^rA)\lesssim2^{-r}\mu(A).
\end{equation}
\end{proposition}

\begin{proof}
First suppose $A\geq2$ and $T\geq2LA$.  If $x\in G_T$, choose $n$ with $f_n(x)\geq T$.  The last assertion of Lemma \ref{lem:marked-packing} gives $g_{A,n}(x)\geq T$; and so, we have the containments
\[
 G_T\subset\{g_A^*\geq T\}
 \subset\{\cM g_{A,N}\gtrsim_d T\}.
\]
Combining the weak $(1,1)$ boundedness of the maximal function, the $L^1$ estimate \eqref{eq:green-L1}, and Lemma \ref{lem:terminal-domination} now yield
\[
 \mu(T)\lesssim\frac1T\|g_{A,N}\|_1
 \lesssim\frac{A}{T}\mu(A).
\]
If $A\leq T<2LA$, then
\[
 \mu(T)\leq\mu(A)\leq2L\frac AT\mu(A).
\]
If $1\leq A<2$ and $T\geq4L$, apply the established estimate at height $2$, and thus obtain
\[
 \mu(T)\lesssim\frac2T\mu(2)
 \leq2\frac AT\mu(A).
\]
The remaining range $T<4L$ follows from monotonicity.  This proves
\eqref{eq:relative-weak}; set $T=2^rA$ to obtain \eqref{eq:relative-dyadic} (which will be the main estimate we use in later calculations).
\end{proof}

\subsubsection{Self-similar subtree tails and dyadic recombination}\label{subsubsec:selfsimtails}

For $R\in\mathscr M_K$, define its descendant maximal contribution by
\begin{equation}\label{eq:local-contribution}
 h_R(x)
 :=
 \max_{\operatorname{depth}(R)\leq n\leq N}
 \#\bigl\{
 Q\in\mathscr Q_n:Q\subseteq R,\ x\in I_Q
 \bigr\}.
\end{equation}
Thus $h_R$ measures only the stacking contributed by the symbolic subtree
rooted at $R$.  The following local tail estimate compares the distributional
tail of this single-subtree contribution with the corresponding tail of the
full construction.

\begin{lemma}[Local tail estimate]\label{lem:local-tail}
For every $R\in\mathscr M_K$ and every $S\geq2$,
\begin{equation}\label{eq:local-tail}
 \mathcal H^1\bigl(\{x:h_R(x)\geq S\}\bigr)
 \leq \ell(R)\mu(S).
\end{equation}
\end{lemma}

\begin{proof}
Let
\[
 r:=\operatorname{depth}(R).
\]
If $R=Q_\omega$, write $z_R:=z_\omega$ for the minimal vertex of $R$,
and set
\[
 a_R:=\pi_\theta(z_R).
\]
For $M\geq0$, define the root-inclusive maximal counting function by
\[
 f_{M,\theta}^{*,0}(y)
 :=
 \max_{0\leq q\leq M}f_{q,\theta}(y).
\]

The symbolic subtree rooted at $R$ is a translated,
$\ell(R)$-dilated copy of the original construction.  More precisely, the
depth-$(r+q)$ descendants of $R$ correspond bijectively to the depth-$q$
cylinders of the original construction.  After projection, this
correspondence gives the exact identity
\begin{equation}\label{eq:local-self-similarity}
 h_R\bigl(a_R+\ell(R)y\bigr)
 =
 f_{N-r,\theta}^{*,0}(y).
\end{equation}

If $r=N$, then $R$ has no proper descendants at depths at most $N$, and
\[
 h_R=\mathbf 1_{I_R}.
\]
Hence $\{h_R\geq S\}=\varnothing$ for $S\geq2$, and
\eqref{eq:local-tail} is immediate.  We may therefore assume that $r<N$.

The depth-zero term is
\[
 f_{0,\theta}=\mathbf 1_{I_{Q_\varnothing}},
\]
so it takes only the values zero and one.  Since $S\geq2$, this term cannot
cause the root-inclusive maximum to reach height $S$.  Consequently,
\[
 \left\{
 y:f_{N-r,\theta}^{*,0}(y)\geq S
 \right\}
 =
 \left\{
 y:\max_{1\leq q\leq N-r}f_{q,\theta}(y)\geq S
 \right\}.
\]
The set on the right is precisely the depth-$(N-r)$ maximal level set
$G_{N-r,S,\theta}$, whose measure is $\mu_{N-r}(S)$.

Applying the affine change of variables
\[
 x=a_R+\ell(R)y
\]
to \eqref{eq:local-self-similarity}, and then using depth monotonicity, gives
\[
 \begin{aligned}
 \mathcal H^1\{x:h_R(x)\geq S\}
 &=
 \ell(R)\,
 \mathcal H^1
 \big(\left\{
 y:f_{N-r,\theta}^{*,0}(y)\geq S
 \right\} \big)\\[1ex]
 &=
 \ell(R)\mu_{N-r}(S)\\[1ex]
 &\leq
 \ell(R)\mu(S).
 \end{aligned}
\]
Note that, in the final inequality, we use the fact that $\mu (S) : = \mu_N (S)$, and $f_{N}^{*,0} \geq f_{N-r}^{*,0}$ pointwise. This proves \eqref{eq:local-tail}.
\end{proof}

We also have the following standard elementary pigeonholing argument, which will allow us to distinguish two distinct ways in which a large total multiplicity can occur.

\begin{lemma}[Dyadic pigeonholing]\label{lem:dyadic-pigeonhole}
There are constants $A_{\mathrm{dp}}>1$ and $c_{\mathrm{dp}}>0$, depending only on $L$, with the following property.

If $\{a_1,...,a_V\}$ is a collection of non-negative numbers of any cardinality $1 \leq V \leq 2LK$ which satisfies the estimate
\[
 \sum_{j=1}^{V} a_j\geq A_{\mathrm{dp}}LKM,
\]
then, letting
\[
 R_K:=\left\lceil\log_2(A_{\mathrm{dp}}LK)\right\rceil.
\]
there necessarily exists some $0\leq r\leq R_K$ such that
\begin{equation}\label{eq:dyadic-pigeonhole}
 \#\{j:a_j\geq2^rM\}
 \geq c_{\mathrm{dp}}\frac{K}{(r+1)^2 2^r}.
\end{equation}
\end{lemma}

\begin{proof}
For $0\leq r\leq R_K$, put
$
N_r:=\#\{j:a_j\geq2^rM\},
$
and suppose that \eqref{eq:dyadic-pigeonhole} fails for every $0 \leq r \leq R_K$.  Since
$2^{R_K}\geq A_{\mathrm{dp}}LK$, one has
\[
 c_{\mathrm{dp}}\frac{K}{(R_K+1)^22^{R_K}}
 \leq\frac{c_{\mathrm{dp}}}{A_{\mathrm{dp}}L(R_K+1)^2}<1.
\]
Hence, the integrality of $N_{R_K}$ implies that
$N_{R_K}=0$, which is equivalent to the uniform upper bound $a_j<2^{R_K}M$.  The terms below $M$ contribute at
most $2LKM$, and the remaining terms lie in the disjoint ranges
$[2^rM,2^{r+1}M)$, $0\leq r<R_K$.  Consequently,
\begin{align*}
 \sum_j a_j
 &\leq2LKM
   +\sum_{r=0}^{R_K-1}2^{r+1}M
     N_r\\
 &<2LKM+2c_{\mathrm{dp}}KM\sum_{r\geq0}\frac1{(r+1)^2}.
\end{align*}
Choosing $c_{\mathrm{dp}}>0$ so that
$2c_{\mathrm{dp}}\sum_{r}(r+1)^{-2}\leq L$ and
$A_{\mathrm{dp}}>3$, the above then implies that
$$
\sum_{j=1}^{V} a_j<3LKM<A_{\mathrm{dp}}LKM
$$
which is a contradiction.
\end{proof}

These results now allow us to pass between scales, and also to pigeonhole a good scale 

\subsubsection{A two-parameter tail estimate for Lemma \ref{lma:L2exceptionalset}}\label{subsubsec:twoparametertail}
We now combine the previous intermediate results to obtain the main combinatorial vehicle for proving Lemma \ref{lma:L2exceptionalset}. Since the author suspects that such a result is likely useful in other applications, we state the two-parameter tail estimate as the following theorem.

\begin{theorem}[A two-parameter tail estimate for $f_{N,\theta_0}^*$]\label{thm:refined-tail}
Fix a digit set $\mathcal{A} : = A_1 \times \cdots \times A_d$ and a direction $\theta_0 \in \mathbb{S}^{d-1}$. Choosing $N \gg 1$, we then let $K = K(N) \gg 1$ be chosen as in Definition \ref{def:stackingparameter} and let $\rho \in (2, \infty)$ be chosen independent of $\mathcal{A}$, $\theta_0$, $N$ and $K$. Define $f_{N,\theta_0}^* : = \max_{1 \leq n \leq N} f_{n,\theta_0}$, where
$$
f_{n,\theta_0} (x) : = \sum\limits_{\omega \in \mathcal{A}^n} \mathbf{1}_{I_{Q_{\omega}}} (x), \qquad \forall x \in \ell_{\theta_0},
$$
and then denote, for all $T \geq 1$, the associated superlevel set measure function
$$
\mu (T) = \mu_{N,K,\rho} (T) : = \mathcal{H}^1 \big( \big\{ x \in \ell_{\theta_0} : f_{N,\theta_0}^* (x) \geq T \big\} \big).
$$

\noindent
Suppose that $\theta_0 \in \mathbb{S}^{d-1}$ is a low-multiplicity direction; which means that
$$
\mu (K) : = \mathcal{H}^1 \big(\big\{x \in \ell_{\theta_0} : f_{N,\theta_0}^* (x) \geq K\big\} \big) \leq K^{-\rho}.
$$
Then, there are constants $A_{\mathrm{dp}} > 0,$ and $C>0$, which depend only upon $L : = \# \mathcal{A}$ and $d$, such that, for all auxiliary scales $M \geq 2$, one has the estimate
\begin{equation}\label{eq:refined-tail}
 \mu(A_{\mathrm{dp}}LKM)
 \leq C\mu(K)
 \sum_{r=0}^{R_K}(r+1)^2 2^r\mu(2^rM),
 \qquad
 R_K=\left\lceil\log_2(A_{\mathrm{dp}}LK)\right\rceil.
\end{equation}
\end{theorem}

We refer to the estimate \eqref{eq:refined-tail} as a \textit{two-parameter tail estimate} because we are able to bound the measure of the (product) depth $KM$ superlevel set in terms of the product of the measure of the depth $K$ superlevel set times a \textit{weighted tail estimate} along a dyadic partition of those depths between $M$ and $KM$.

\begin{proof}
Construct the stopping time forest $\mathscr M_K$ using Lemma \ref{lem:marked-packing}.  Let $x\in G_{A_{\mathrm{dp}}LKM}$ and choose a depth $n$ with
\[
 f_n(x)\geq A_{\mathrm{dp}}LKM.
\]
Let
\[
 m = n_K (x) := \min \big\{1 \leq k \leq N : f_{k} (x) \geq 2 K \big\}
\]
be the first depth at which the multiplicity over $x$ reaches $2K$.  The first-crossing estimate then gives
$$
 2K\leq f_m(x)<2LK.
$$
In particular, the family
$$
 \mathscr P_m(x):=\{P\in\mathscr Q_m:x\in I_P\}
$$
contains fewer than $2LK$-many cylinders.  Moreover, since
$f_n(x)\geq A_{\mathrm{dp}}LKM\geq2LK$, we necessarily have $m<n$.

Let $Q\in\mathscr Q_n$ be any depth-$n$ cylinder over $x$, and denote its unique depth-$m$ ancestor by $P_m(Q)$.  The condition $x\in I_Q$ implies that $x\in I_{P_m(Q)}$; and so,
$P_m(Q)$ was marked at the first crossing over $x$.

Every marked cylinder is contained in a member of $\mathscr M_K$.  This member is unique: if two members of the prefix-free family $\mathscr M_K$ both contained $P_m(Q)$, then they would be comparable in the symbolic prefix order and hence would have to coincide.  Let $R(Q)\in\mathscr M_K$
denote this unique retained cylinder.  We then have
\[
 Q\subseteq P_m(Q)\subseteq R(Q).
\]

For each $R\in\mathscr M_K$, define
\[
 a_R=a_R(x,n)
 :=
 \#\bigl\{
 Q\in\mathscr Q_n:
 x\in I_Q,\ R(Q)=R
 \bigr\}.
\]
Thus $a_R$ is the portion of the depth-$n$ multiplicity over $x$ contributed
by cylinders assigned to the symbolic subtree rooted at $R$.  Every depth-$n$ cylinder over $x$ is assigned to exactly one such subtree, and therefore
\[
 \sum_{R\in\mathscr M_K}a_R
 =
 \#\{Q\in\mathscr Q_n:x\in I_Q\}
 =
 f_n(x)
 \geq A_{\mathrm{dp}}LKM.
\]
Furthermore, every $R$ for which $a_R>0$ contains at least one cylinder in
$\mathscr P_m(x)$.  Since $\#\mathscr P_m(x)<2LK$, and several members of $\mathscr P_m(x)$ may be contained in the same retained cylinder, at most $2LK$ of the numbers $a_R$ are nonzero. Finally, all cylinders counted by $a_R$ are depth-$n$ descendants of $R$. The depth-$n$ contribution appearing in the definition of $h_R(x)$ is therefore at least $a_R$, so
\[
 h_R(x)\geq a_R.
\]
We have thus decomposed the large multiplicity at $x$ into at most $2LK$ nonzero single-subtree contributions satisfying
\[
 \sum_{R\in\mathscr M_K}a_R=f_n(x)
 \geq A_{\mathrm{dp}}LKM,
 \qquad
 0\leq a_R\leq h_R(x).
\]

For each intermediate dyadic depth $0\leq r\leq R_K$, we thus let
$$
C_r(x):=\#\{R\in\mathscr M_K:h_R(x)\geq2^rM\},
$$
and define an associated depth-sensitive threshold by setting
$$
q_r:=c_{\mathrm{dp}}\frac{K}{(r+1)^2 2^r}.
$$
Lemma \ref{lem:dyadic-pigeonhole} gives
\begin{equation}\label{eq:high-cover}
 G_{A_{\mathrm{dp}}LKM}
 \subset\bigcup_{r=0}^{R_K}\{x:C_r(x)\geq q_r\}.
\end{equation}
An application of Markov's inequality, together with the local-tail Lemma \ref{lem:local-tail} and the geometric length bound \eqref{eq:marked-packing} along $\mathscr{M}_A$ (which was a part of Lemma \ref{lem:marked-packing}), thus gives
\begin{align*}
 \cH^1 \big(\big\{C_r\geq q_r\big\}\big)
 &\leq q_r^{-1}\int C_r(x)\,dx\\
 &=q_r^{-1}\sum_{R\in\mathscr M_K}
   \cH^1 \big(\{x:h_R(x)\geq2^rM\} \big)\\
 &\lesssim\frac{(r+1)^2 2^r}{K}\mu(2^rM)
   \sum_{R\in\mathscr M_K}\ell(R)\\
 &\lesssim(r+1)^2 2^r\mu(K)\mu(2^rM).
\end{align*}
Using this estimate for all $0 \leq r \leq R_K$ on the right-hand side of \eqref{eq:high-cover} and summing in $r$ thus proves \eqref{eq:refined-tail}.
\end{proof}

We state a consequence of Theorem \ref{thm:refined-tail} which will be utilized in the upcoming calculation.

\begin{corollary}[Polylogarithmic tail]\label{cor:polylog-tail}
For all $K,M\geq2$,
\begin{equation}\label{eq:polylog-tail}
 \mu(A_{\mathrm{dp}}LKM)
 \lesssim(1+\log K)^3\mu(K)\mu(M).
\end{equation}
\end{corollary}

\begin{proof}
By Proposition \ref{prop:relative-weak},
$
 2^r\mu(2^rM)\lesssim\mu(M).
$
Inserting this estimate into \eqref{eq:refined-tail} and using the simple upper bound
$$
 \sum_{r=0}^{R_K}(r+1)^2\lesssim(1+\log K)^3
$$
thus proves the corollary.
\end{proof}

We are now finally in a position to prove the $L^2$ estimate for the maximal counting function along low-multiplicity directions.

\begin{theorem}[Low-multiplicity $L^2$ estimate for $\rho>2$]\label{thm:improved-L2}
Let $\rho>2$.  Suppose that
\begin{equation}\label{eq:low-multiplicity}
 \mu(K) : = \mathcal{H}^1 \big(\big\{x \in \mathbb{R} : f_{N}^* (x) \geq K \big\} \big)\leq K^{-\rho}.
\end{equation}
Then, for all sufficiently large $K=K(N;L,d,\rho)$, one has
\begin{equation}\label{eq:improved-L2}
 \|f_N^*\|_2^2\lesssim_{L,d,\rho}K.
\end{equation}
Consequently, we have that
\[
 \max_{1\leq n\leq N}\|f_n\|_2^2\lesssim_{L,d,\rho}K.
\]
\end{theorem}

To prove this result, we convert the single level-set hypothesis of \eqref{eq:low-multiplicity} at height $K$ into a bound for the complete dyadic energy.

\begin{proof}
We organize the proof through the dyadic distribution function of $f_N^*$.  Define
\begin{equation}\label{eq:dyadic-energy}
 \mathcal E_2:=\sum_{j\geq0}2^{2j}\mu(2^j).
\end{equation}
Since $N \gg 1$ is fixed, the previous sum is finite because $f_N^*\leq L^N$.  Using a layer-cake estimate for the $L^2$ norm, we easily obtain that
\begin{equation}\label{eq:energy-controls-L2}
 \|f_N^*\|_2^2\lesssim\mathcal E_2.
\end{equation}
To verify this, observe that on the annulus $A_j:=\{2^j\leq f_N^*<2^{j+1}\}$, one has
\[
 (f_N^*)^2\mathbf1_{A_j}\leq2^{2j+2}\mathbf1_{G_{2^j}}.
\]
As the sets $A_j$ are disjoint, we thus have
\[
 \|f_N^*\|_2^2
 =\sum_{j\geq0}\int_{A_j}(f_N^*)^2
 \leq4\sum_{j\geq0}2^{2j}\mu(2^j)
 =4\mathcal E_2,
\]
which gives \eqref{eq:energy-controls-L2}.

We now split the energy estimate for $\mathcal{E}_2$ into its corresponding low and high parts, with cutoff heights comparable to $K$.  To this end, put
\[
 D_K:=A_{\mathrm{dp}}LK
 \quad \textrm{and} \quad 
 R_K:=\lceil\log_2D_K\rceil.
\]
For $0\leq j\leq R_K$, Proposition \ref{prop:relative-weak} with $A=1$ gives
\[
 \mu(2^j)\lesssim2^{-j}\mu(1).
\]
Directly from the definition, we know that $G_1$ is contained in the projection of the cylinders at depth one; and so, we easily have that
\[
 \mu(1)\leq\sum_{Q\in\mathscr Q_1}|I_Q|
 =L\sigma_\theta L^{-1}=\sigma_\theta\leq\sqrt d.
\]
Therefore the low part of the energy satisfies
\begin{equation}\label{eq:low-energy}
 \sum_{j=0}^{R_K}2^{2j}\mu(2^j)
 \lesssim\sum_{j=0}^{R_K}2^j
 \lesssim K.
\end{equation}

For the high-energy tail, fix $j\geq R_K+1$ and set $M_j:=2^{j-R_K}$.  Since $D_K\leq2^{R_K}$, we know that
$ D_KM_j\leq2^j$ and so obtain that $\mu(2^j)\leq\mu(D_KM_j)$. We then apply the refined two-parameter tail estimate of Theorem \ref{thm:refined-tail}, which gives
\begin{equation}\label{eq:high-energy-pointwise}
 \mu(2^j)
 \lesssim\mu(K)\sum_{r=0}^{R_K}(r+1)^2 2^r
 \mu(2^{j-R_K+r}).
\end{equation}
We multiply by the energy weight $2^{2j}$ and write $q=j-R_K+r$.  Since
\[
 2^{2j}2^r
 =2^{2R_K}2^{-r}2^{2q},
\]
we obtain
\begin{align*}
 2^{2j}\mu(2^j)
 &\lesssim2^{2R_K}\mu(K)
 \sum_{r=0}^{R_K}(r+1)^2 2^{-r}
 \bigl[2^{2q}\mu(2^q)\bigr].
\end{align*}
Summing first over $j\geq R_K+1$ and then over $r$ gives
\begin{align*}
 \sum_{j\geq R_K+1}2^{2j}\mu(2^j)
 &\lesssim2^{2R_K}\mu(K)
 \sum_{r=0}^{R_K}(r+1)^22^{-r}
 \sum_{q\geq1+r}2^{2q}\mu(2^q)\\[1ex]
 &\leq C2^{2R_K}\mu(K)\mathcal E_2 \\[1ex]
 & \lesssim K^2\mu(K)\mathcal E_2,
\end{align*}
because $\sum_{r\geq0}(r+1)^22^{-r}<\infty$ and
$2^{R_K}\lesssim_LK$. 

Putting together this high-energy estimate together with the low-energy estimate \eqref{eq:low-energy}, we obtain
\begin{equation}\label{eq:energy-absorption}
 \mathcal E_2
 \lesssim K+K^2\mu(K)\mathcal E_2.
\end{equation}
Using the low-multiplicity direction hypothesis \eqref{eq:low-multiplicity}, we can bound the high-energy part so that
\[
 K^2\mu(K)\leq K^{2-\rho}=o(1).
\]
We then choose $K$ sufficiently large so that the implicit constant multiplying $K^2\mu(K)$ in \eqref{eq:energy-absorption} is at most $1/2$.  We may then move that term to the left-hand side and obtain
\[
 \mathcal E_2\lesssim K.
\]
Equation \eqref{eq:energy-controls-L2} now proves \eqref{eq:improved-L2}.
\end{proof}

\begin{remark}
The restriction $\rho>2$ comes only from the absorption coefficient
$K^2\mu(K)\leq K^{2-\rho}$.  No additional power is lost through the depth
of the tree or the ambient dimension; and so, this handling of the high-energy part of the energy estimate is exactly the moment where we need the $\rho > 2$ hypothesis.
\end{remark}

\subsection{Terminal marking and block propagation}\label{subsec:reverseholder}

We now turn to the complementary high-multiplicity alternative
$\theta\notin E_{N,K}$ and prove the propagation inequality of Lemma
\ref{lma:reverseholder}.  When comparing this argument to the previous low-multiplicity argument, the associated level set $G_K$ is now too large in measure
for the energy argument from Section \ref{subsec:L2exceptional} to be leveraged. Instead, we take a different approach,
using the self-similar tree structure to convert a proportion of the
terminal words in every block into a family whose projection has high local
multiplicity. In fact, showing that this proportion is an absolute and fixed positive proportion of the entire family of terminal words will be the exact place where we utilize the high-multiplicity hypothesis for the direction $\theta \not\in E_{N,K}$.

There are two steps.  
\begin{enumerate}
\item The terminal marking lemma replaces the witness
depth $n(x)$, which may vary with $x\in G_K$, by a single family
$\mathscr W_K\subset\mathscr Q_N$ at the common terminal depth $N$. One can think of $\mathscr{W}_K$ as a sub-family of cylinders, for which the corresponding projection inequality holds at the corresponding earlier depth and with the appropriate stacking threshold at this depth.
\medskip
\item We then copy (using the self-similar structure of our word tree) the same terminal pattern $\mathscr{W}_K$ inside every block at the earlier depth which did not belong to this good family.  We show that the
unconverted (i.e. \textit{black}) family loses a fixed proportion of its words in
each block, while the converted (i.e. \textit{green}) family is controlled by the weak
$(1,1)$ inequality for the Hardy--Littlewood maximal operator. This weak $(1,1)$ control then allows us to conclude that the projection of those converted blocks accounts for a fixed and constant positive proportion of the terminal words.
\end{enumerate}

\subsubsection{Terminal marking at a common depth}

We begin by making the terminal selection precise.  For a cylinder
$Q\in\mathscr Q_m$, recall that $\operatorname{Desc}_N(Q)$ denotes its family
of descendants at depth $N$.

\begin{lemma}[Terminal marking lemma]\label{lma:terminalmarking}
Let $A\geq2$.  There exist a pairwise disjoint family of intervals
$\{J_\lambda\}_{\lambda\in\Lambda}$ such that, for every $\lambda\in\Lambda$, there exists a
depth $n(\lambda)\leq N$, a point $x_\lambda\in J_\lambda$, as well as a family of distinct witness cylinders
\[
 \mathcal C(\lambda)
 =\{Q_{\lambda,1},\ldots,Q_{\lambda,\lceil A\rceil}\}
 \subset\mathscr Q_{n(\lambda)}
\]
such that the following properties are all satisfied simultaneously.
\begin{enumerate}
\item The witness projections all contain the common point
\[
 x_\lambda\in
 \bigcap_{s=1}^{\lceil A\rceil}I_{Q_{\lambda,s}},
 \textrm{ where one also has }
 J_\lambda=\bigcup_{s=1}^{\lceil A\rceil}I_{Q_{\lambda,s}}.
\]
\item The intervals $J_{\lambda}$ form a Vitali-type cover for the associated superlevel set $G_A$, in the sense that
\[
 G_A\subset\bigcup_{\lambda\in\Lambda}5J_\lambda.
\]
\item The associated witness cylinders contained in
$\bigcup_{\lambda\in\Lambda}\mathcal C(\lambda)$ are all pairwise incomparable
in the symbolic prefix order.
\end{enumerate}
Consequently, the family
\begin{equation}\label{eq:terminal-marked-family}
 \mathscr W_A
 :=\bigcup_{\lambda\in\Lambda}
   \bigcup_{Q\in\mathcal C(\lambda)}
   \operatorname{Desc}_N(Q)
 \subset\mathscr Q_N
\end{equation}
is a disjoint union as a family of terminal words and satisfies
\begin{equation}\label{eq:markedterminalcount}
 \#\mathscr W_A\gtrsim_d A\mu(A)L^N.
\end{equation}
\end{lemma}

\begin{proof}
Fix $x\in G_A$.  By the definition of $G_A$, there is a witness depth
$n(x)\leq N$ for which $f_{n(x)}(x)\geq A$.  Since the counting function is
integer-valued, we may choose $\lceil A\rceil$ distinct cylinders
\[
 Q_{x,1},\ldots,Q_{x,\lceil A\rceil}\in\mathscr Q_{n(x)}
 \quad\text{such that}\quad
 x\in\bigcap_{s=1}^{\lceil A\rceil}I_{Q_{x,s}}.
\]
Define
\[
 J_x:=\bigcup_{s=1}^{\lceil A\rceil}I_{Q_{x,s}}.
\]
Every interval in this union contains $x$ and has the common length
$\sigma_\theta L^{-n(x)}$.  Their union is therefore an interval, and its
length satisfies
\begin{equation}\label{eq:vitali-block-length}
 \sigma_\theta L^{-n(x)}
 \leq |J_x|
 \leq2\sigma_\theta L^{-n(x)}.
\end{equation}
The family $\{J_x:x\in G_A\}$ covers $G_A$.  The one-dimensional Vitali
covering lemma furnishes a pairwise disjoint subfamily
$\{J_\lambda\}_{\lambda\in\Lambda}$ such that
\begin{equation}\label{eq:vitali-selected-cover}
 G_A\subset\bigcup_{\lambda\in\Lambda}5J_\lambda,
 \qquad
 \mu(A)\leq5\sum_{\lambda\in\Lambda}|J_\lambda|.
\end{equation}
Here $5J_\lambda$ denotes the interval with the same centre as $J_\lambda$
and five times its length.  We call the selected intervals the
\emph{Vitali blocks}.  For each $\lambda$, retain the corresponding point
$x_\lambda$, depth $n(\lambda)$, and witness family
\[
 \mathcal C(\lambda)
 :=\{Q_{\lambda,1},\ldots,Q_{\lambda,\lceil A\rceil}\}.
\]

We next verify the symbolic incomparability asserted in the lemma.  Distinct
witness cylinders associated with the same $\lambda$ have the same depth and
are therefore incomparable.  Suppose that witness cylinders associated with
different Vitali blocks were comparable.  After interchanging the two
blocks if necessary, we would have
\[
 Q_{\lambda',s'}\subseteq Q_{\lambda,s}.
\]
Projection preserves containment, and hence
\[
 x_{\lambda'}
 \in I_{Q_{\lambda',s'}}
 \subseteq I_{Q_{\lambda,s}}
 \subseteq J_\lambda.
\]
But $x_{\lambda'}\in J_{\lambda'}$ as well, contradicting the disjointness
of $J_\lambda$ and $J_{\lambda'}$.  Thus all selected witness cylinders are
pairwise incomparable.

It follows that their depth-$N$ descendant families are pairwise disjoint as
collections of words.  Therefore the union in
\eqref{eq:terminal-marked-family} is disjoint and
\begin{align*}
 \#\mathscr W_A
 &=\sum_{\lambda\in\Lambda}
   \sum_{Q\in\mathcal C(\lambda)}
   \#\operatorname{Desc}_N(Q)\\
 &=\lceil A\rceil
   \sum_{\lambda\in\Lambda}L^{N-n(\lambda)}.
\end{align*}
By the upper bound in \eqref{eq:vitali-block-length} and the inequality
$\sigma_\theta\leq\sqrt d$,
\[
 L^{N-n(\lambda)}
 \geq\frac{1}{2\sqrt d}L^N|J_\lambda|.
\]
Combining this with \eqref{eq:vitali-selected-cover} gives
\[
 \#\mathscr W_A
 \geq\frac{A}{2\sqrt d}L^N
      \sum_{\lambda\in\Lambda}|J_\lambda|
 \gtrsim_d A\mu(A)L^N,
\]
which proves \eqref{eq:markedterminalcount}.
\end{proof}

\subsubsection{Iteration of the terminal pattern}

\begin{figure}[H]
  \centering
  \includegraphics[width=.98\linewidth]{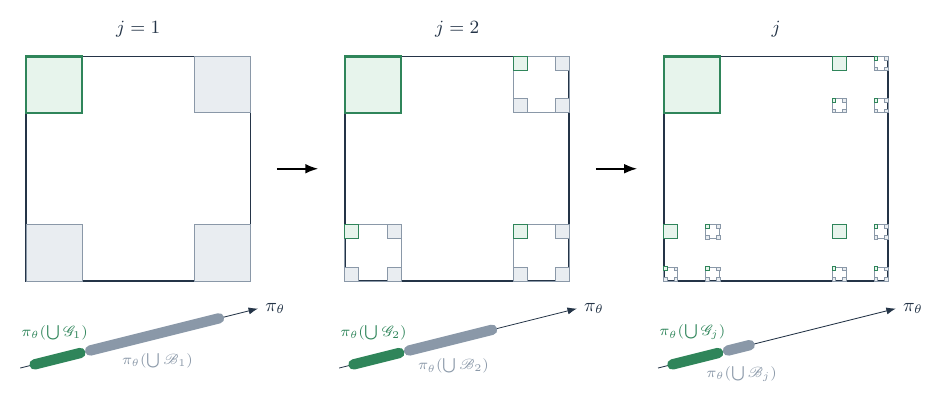}
  \caption{Block propagation across the depths $N,2N,\ldots,jN$.  Inside
    every surviving black block, the same terminal pattern $\mathscr W_K$ is
    copied and converted to green.  Green descendants persist, while the
    black family loses at least a $p$-proportion in every block.  The lower
    target lines depict the corresponding green and black projections.}
  \label{fig:block-propagation}
\end{figure}

We now apply the terminal marking lemma at height $K$ and iterate the
resulting pattern $\mathscr{W}_K$ in blocks with $N$-gap generations.  Figure
\ref{fig:block-propagation} illustrates the green--black decomposition used
below. The final result of this application is a complete proof of Lemma \ref{lma:reverseholder}.

\begin{proof}[Proof of Lemma \ref{lma:reverseholder}]
Throughout the proof, we view $\theta\notin E_{N,K}$, and so suppress the direction
notation in all of our calculations.

By assumption, we know that
\begin{equation}\label{eq:outsideexceptionalmu}
 \mu(K)>K^{-\rho}.
\end{equation}
Apply Lemma \ref{lma:terminalmarking} with $A=K$, and define two auxiliary parameters
\[
 \delta:=\frac{\#\mathscr W_K}{L^N},
 \qquad
 p:=\min\left\{\delta,\frac12\right\}.
\]
In words: $\delta$ is the actual proportion of marked terminal words, while $p$
is an appropriately truncated conversion of this proportion that is always bounded away from one.  The
terminal marking estimate already guarantees that
\begin{equation}\label{eq:terminal-density-lower-bound}
 \delta\gtrsim_d K\mu(K),
\end{equation}
and so the truncation $\delta \mapsto p$ does not change this lower bound by more than a fixed
constant.  Indeed, Proposition \ref{prop:relative-weak}, applied with $A=1$
and $T=K$, gives
\[
 K\mu(K)\lesssim\mu(1)\lesssim_d1.
\]
If $\delta\leq1/2$, then $p=\delta$ and
\eqref{eq:terminal-density-lower-bound} applies directly.  If
$\delta>1/2$, then $p=1/2$, which is still bounded below by a fixed
multiple of $K\mu(K)$.  Consequently,
\begin{equation}\label{eq:conversionproportion}
 p\gtrsim_d K\mu(K)\gtrsim K^{-(\rho-1)},
\end{equation}
where the final inequality follows from \eqref{eq:outsideexceptionalmu}.

We now make the self-similar block iteration explicit.  Let
\[
 \mathscr W_K^*
 :=\{\omega\in\mathcal A^N:Q_\omega\in\mathscr W_K\}
\]
be the set of marked words of length $N$.  Begin with
\[
 \mathscr B_0:=\{Q_\varnothing\},
 \qquad
 \mathscr G_0:=\varnothing.
\]
Suppose that the black family $\mathscr B_{j-1}$ has been defined at depth
$(j-1)N$.  For every black block parent $Q_\tau\in\mathscr B_{j-1}$, use
the canonical identification of its $N$-block descendants with
$\mathcal A^N$ and declare
\[
 Q_{\tau\omega}\ \text{GREEN if }\omega\in\mathscr W_K^*,
 \qquad
 Q_{\tau\omega}\ \text{BLACK if }\omega\notin\mathscr W_K^*.
\]
All block descendants of a green parent remain green.  Equivalently, we may
define
\begin{equation}\label{eq:black-family-recursion}
 \mathscr B_j
 :=\left\{
 Q_{\tau\omega}:
 Q_\tau\in\mathscr B_{j-1},\ 
 \omega\in\mathcal A^N\setminus\mathscr W_K^*
 \right\},
 \qquad
 \mathscr G_j:=\mathscr Q_{jN}\setminus\mathscr B_j.
\end{equation}
Thus $\mathscr G_j$ and $\mathscr B_j$ form a disjoint partition of
$\mathscr Q_{jN}$.

We estimate the black and green projections separately.

\smallskip
\noindent\emph{The black projection.}
Each black block parent has exactly $L^N-\#\mathscr W_K$ black descendants
after the next block.  Since $p\leq\delta$, this number is bounded by
\[
 L^N-\#\mathscr W_K
 =(1-\delta)L^N
 \leq(1-p)L^N.
\]
Starting from $\#\mathscr B_0=1$ and iterating the recurrence gives
\begin{equation}\label{eq:blackcylindercount}
 \#\mathscr B_j
 \leq L^{jN}(1-p)^j.
\end{equation}
Every cylinder in $\mathscr Q_{jN}$ has projected interval of length
$\sigma_\theta L^{-jN}$.  Therefore
\begin{align}
 \mathcal H^1\left(\bigcup_{Q\in\mathscr B_j}I_Q\right)
 &\leq\sum_{Q\in\mathscr B_j}|I_Q|\notag\\
 &=\sigma_\theta L^{-jN}\#\mathscr B_j\notag\\
 &\leq\sigma_\theta(1-p)^j
 \leq\sqrt d\,e^{-pj}.
 \label{eq:blackprojectionbound}
\end{align}
This is the exponential decay furnished by the repeated removal of a
$p$-proportion of every surviving black block.

\smallskip
\noindent\emph{The green projection.}
Define the depth-$jN$ green counting function by
\begin{equation}\label{eq:green-block-counting}
 g_j:=\sum_{Q\in\mathscr G_j}\mathbf1_{I_Q}.
\end{equation}
Since $\#\mathscr G_j\leq L^{jN}$,
\begin{equation}\label{eq:green-block-L1}
 \|g_j\|_1
 =\sigma_\theta L^{-jN}\#\mathscr G_j
 \leq\sigma_\theta
 \leq\sqrt d.
\end{equation}
We claim that every point in the green projection belongs to a region on
which $g_j$ has average value comparable to $K$.

To this end, fix a point
\[
 x\in\bigcup_{Q\in\mathscr G_j}I_Q
\]
and choose $Q\in\mathscr G_j$ with $x\in I_Q$.  Let $b\in\{1,\ldots,j\}$
be the first block for which the depth-$bN$ ancestor of $Q$ is green.  Its
block parent $Q_\tau$ at depth $(b-1)N$ is black, and the length-$N$ word
that converts this parent belongs to $\mathscr W_K^*$.  By the decomposition
\eqref{eq:terminal-marked-family}, that terminal word descends from a unique
witness family $\mathcal C(\lambda)$ selected by the terminal marking lemma.

Write $n=n(\lambda)$ and
\[
 \mathcal C(\lambda)
 =\{Q_{\lambda,1},\ldots,Q_{\lambda,\lceil K\rceil}\}.
\]
Copying this family into the block parent $Q_\tau$ produces
$\lceil K\rceil$ incomparable cylinders
\[
 \widetilde Q_{\lambda,1},\ldots,
 \widetilde Q_{\lambda,\lceil K\rceil}
\]
at the common global depth $(b-1)N+n$.  Their projections have a common
point, which we denote by $y_{\lambda,\tau}$, and each copied witness
cylinder has geometric length
\begin{equation}\label{eq:copied-witness-length}
 \ell_*:=L^{-((b-1)N+n)}.
\end{equation}
The cylinder $Q$ descends from one of these copied witnesses.  Hence $x$ and
$y_{\lambda,\tau}$ both belong to the projection of that witness, and
\begin{equation}\label{eq:green-point-distance}
 |x-y_{\lambda,\tau}|\leq\sigma_\theta\ell_*.
\end{equation}

When the block parent converts, all of the depth-$bN$ descendants of every
copied witness cylinder are declared green.  Since green status persists,
all of their descendants at depth $jN$ belong to $\mathscr G_j$.  Each copied
witness has exactly $L^{jN-((b-1)N+n)}$ such descendants, whose total
projected $L^1$ mass is
\begin{equation}\label{eq:one-witness-green-mass}
 L^{jN-((b-1)N+n)}\sigma_\theta L^{-jN}
 =\sigma_\theta\ell_*.
\end{equation}
The copied witnesses are incomparable, so these descendant families are
disjoint as collections of words and their contributions to $g_j$ may be
summed.

To obtain a centred maximal-function estimate, set
\[
 J_x:=[x-2\sigma_\theta\ell_*,x+2\sigma_\theta\ell_*].
\]
Every copied witness projection contains $y_{\lambda,\tau}$ and has length
$\sigma_\theta\ell_*$.  Together with
\eqref{eq:green-point-distance}, this shows that every copied witness
projection, and hence every projection of one of its depth-$jN$ descendants,
is contained in $J_x$.  Using \eqref{eq:one-witness-green-mass} and the fact
that there are at least $K$ copied witnesses, we obtain
\[
 \int_{J_x}g_j(t)\,dt
 \geq K\sigma_\theta\ell_*,
 \qquad
 |J_x|=4\sigma_\theta\ell_*.
\]
Therefore the centred Hardy--Littlewood maximal function satisfies
\[
 \cM g_j(x)\geq\frac{1}{|J_x|}\int_{J_x}g_j(t)\,dt
 \geq\frac K4.
\]
Since $x$ was arbitrary,
\[
 \bigcup_{Q\in\mathscr G_j}I_Q
 \subseteq\{x:\cM g_j(x)\geq K/4\}.
\]
The weak $(1,1)$ inequality and \eqref{eq:green-block-L1} now give
\begin{align}
 \mathcal H^1\left(\bigcup_{Q\in\mathscr G_j}I_Q\right)
 &\leq\mathcal H^1\{x:\cM g_j(x)\geq K/4\}\notag\\
 &\lesssim K^{-1}\|g_j\|_1
 \lesssim_d K^{-1}.
 \label{eq:greenprojectionbound}
\end{align}

We finally choose the number of blocks.  Put $j=J_{K,\rho}$, where
$J_{K,\rho}$ is defined in \eqref{eq:propagationblocks}.  From
\eqref{eq:conversionproportion},
\[
 pj
 \geq cK^{-(\rho-1)}
 \left\lceil
 C_{\mathrm{prop}}K^{\rho-1}\log(2+K)
 \right\rceil
 \geq cC_{\mathrm{prop}}\log(2+K).
\]
Choose the fixed constant $C_{\mathrm{prop}}$ sufficiently large that
$cC_{\mathrm{prop}}\geq1$.  Then
\begin{equation}\label{eq:black-final-decay}
 e^{-pJ_{K,\rho}}\leq(2+K)^{-1}\lesssim K^{-1}.
\end{equation}
Because $\mathscr G_j$ and $\mathscr B_j$ partition $\mathscr Q_{jN}$,
their projected unions cover $\pi_\theta(\mathcal S^{jN})$.  Combining
\eqref{eq:blackprojectionbound}, \eqref{eq:greenprojectionbound}, and
\eqref{eq:black-final-decay}, we conclude that
\begin{align*}
 \mathcal H^1\bigl(\pi_\theta(\mathcal S^{NJ_{K,\rho}})\bigr)
 &\leq
 \mathcal H^1\left(\bigcup_{Q\in\mathscr G_{J_{K,\rho}}}I_Q\right)
 +
 \mathcal H^1\left(\bigcup_{Q\in\mathscr B_{J_{K,\rho}}}I_Q\right)\\
 &\lesssim_d K^{-1}.
\end{align*}
This proves Lemma \ref{lma:reverseholder}.
\end{proof}

\begin{remark}[Propagation cost and the role of $\rho$]
The block iteration uses the exponent $\rho$ only through the conversion
proportion
\[
 p\gtrsim K^{-(\rho-1)}.
\]
Reducing the black projection to size $O(K^{-1})$ requires
$pJ_{K,\rho}\gtrsim\log K$, which explains the propagation cost
\[
 J_{K,\rho}\asymp K^{\rho-1}\log(2+K).
\]
The stricter assumption $\rho>2$ does not arise from the propagation
mechanism itself; it enters earlier, through the energy-absorption
coefficient $K^2\mu(K)\leq K^{2-\rho}$, which accounts for the high-part of the energy estimate necessary to prove our $L^2$ estimate for the low multiplicity directions.  Letting $\rho\downarrow2$ makes the
block cost $K^{1+o(1)}$, whereas the endpoint $\rho=2$ is not reached because
the present $L^2$ argument no longer yields a strictly absorbable
coefficient.  This is a limitation of the proof rather than a claimed
obstruction for the Favard-length problem.
\end{remark}

\clearpage
\subsection*{Figures for the combinatorial argument}

The figures are collected here in the order in which their constructions
occur in Sections \ref{subsec:combinatorialoutline}--\ref{subsec:reverseholder}.

\clearpage

\section{The Exponent Calculation and Principal Examples}\label{sec:exponentexamples}

All of the analytic and combinatorial ingredients are now in place; and so, the
purpose of this final section is to make their quantitative consequences
explicit.  There are two distinct scales in the proof.  At the word depth
$N$, the Set of Small Values and Set of Large Values estimates determine how large the stacking threshold $K$
may be while still keeping the exceptional set of directions small.  The
propagation lemma then carries the resulting high-multiplicity information
from depth $N$ to the terminal construction depth
\[
 M\asymp NJ_{K,\rho}.
\]
The second passage \textit{is not cost-free}. However, Section \ref{sec:comblemmata} precisely quantifies this cost, showing that it is only $K^{\rho-1}\log(2+K)$ for every $\rho>2$.  Keeping the
analytic and combinatorial stages of our argument separate, as well as the careful book-keeping which we have engaged in so far, allows us both to optimize the general power law theorem, while also isolating the special Fourier input for the four-corner Cantor set, as well as other rational product Cantor sets whose digit sets have mask polynomials with special cyclotomic structure.

Four consequences will emerge from our previous arguments and the explicit calculations in this section.
\begin{enumerate}
\item The SSV--SLV mechanism gives the most general form of
the power and sub-power alternatives stated in Theorem \ref{thm:totalFavbound}.
\medskip
\item The Nazarov--Peres--Volberg Fourier estimate, combined with the
improved propagation, gives the exponent $1/5-u$ for the four-corner Cantor
set.
\medskip
\item When all coordinate mask polynomials are zero-free on the unit
circle, both the SSV deletion and the SLV construction become trivial.  In
that case the exponent is governed only by the minimum moduli of the mask
polynomials. We prove that the explicit planar digit set $\{0,1,2,6,9\}$ gives the much improved exponent of
$1/3-u$.
\medskip
\item We then provide a six-digit example in base $36$ which has minimum modulus strictly larger than one. Therefore, its associated power law exponent (which we again calculate explicitly) crosses the one-third threshold.
\end{enumerate}

The purpose of these four consequences is to demonstrate that, while our current methods prove general power laws for large classes of self-similar sets, in order to better refine the exponents (and, thus, the asymptotics of this problem) one must examine the underlying cyclotomic structure of the involved digit sets. Further, while the four-corner Cantor set is a geometrically natural canonical example of such a set (since, in the plane, base $L = 4$ is the smallest base in which one can construct an unrectifiable rational product Cantor set), it is not exemplary from the perspective of improved asymptotics. In fact, our observations thus complement the work of Bateman and Volberg \cite{BV} as well as the very recent work of Chang, Shmerkin and Suomala \cite{CSS}, which demonstrates that the four-corner Cantor set has $4$-adic neighbourhoods whose Favard length decays strictly slower than the general lower bound.

\subsection{A guide to the parameters}

We begin by collecting the parameters which enter the final calculation.  It
is important that they be chosen in the order recorded in the parameter table at the end of this section: the digit data are fixed first, the analytic
parameters are then chosen at depth $N$, and only afterwards is the
propagation length selected. One can see Section \ref{subsec:parametertable} for the included parameter table, which lists all parameters, their role in the argument, and the main calculative restrictions placed by the web of the argument.

The dependency chain may be summarized as follows.  The digit sets determine
$L$ and all SSV and SLV constants.  We then choose an analytic exponent
$\kappa$ and set $K=N^\kappa$.  The scale parameter $s_0$, and hence the
low-frequency length $m$, must satisfy two competing requirements.  The SLV
mass must be large enough to dominate the error term of size $KL^{-n}$, while
the pointwise loss $L^{-2Dm}$ must remain small enough to contradict the
Fourier upper bound.  Only after these choices have been fixed do we apply the
propagation lemma and pass from $N$ to $M$.

We shall call $\delta>0$ \emph{admissible} if, for the rational product Cantor
set under consideration,
\[
 \Fav\bigl(\mathcal N_{L^{-M}}(\mathcal S^\infty)\bigr)
 \lesssim M^{-\delta}
\]
for every sufficiently large $M$, with an implicit constant allowed to depend
upon the fixed digit data and the fixed auxiliary exponents.

\subsection{The general exponent mechanism}

We first treat the ordinary power-law branch.  Write
\[
 K=N^\kappa,
 \qquad
 m=\lceil s_0\log_LN\rceil.
\]
The SLV witness has measure at least
$c_\Gamma L^{-(1-\eta)m}$.  Salem's argument therefore gains the factor
$L^{\eta m}$ over the baseline scale $L^{-n}$.  In order for this gain to
dominate both the Poisson-localization and SSV-deletion errors, we require
\begin{equation}\label{eq:SLVscalecondition}
 \kappa<\eta s_0.
\end{equation}
On the other hand, Proposition \ref{prop:keyproposition} contributes the
squared pointwise loss $L^{-2Dm}$.  The final comparison
\eqref{eq:finalparametercontradiction} fails, as required, provided that
\begin{equation}\label{eq:Fourierscalecondition}
 2Ds_0<1.
\end{equation}
If $D>0$, the two restrictions are equivalent to the existence of a choice
\[
 \frac{\kappa}{\eta}<s_0<\frac1{2D}.
\]
Such a choice is possible precisely when
\begin{equation}\label{eq:optimizedkappa}
 \kappa<\frac{\eta}{2D}.
\end{equation}
This elementary interval condition is the entire analytic optimization in the
power branch.  It identifies the largest stacking scale supported by the SSV
and SLV inputs before the combinatorial propagation cost is imposed.

Fix $\beta\in(0,1)$.  Once \eqref{eq:SLVscalecondition},
\eqref{eq:Fourierscalecondition}, and the strict SSV margins are satisfied,
the argument of Sections \ref{sec:rieszproductbounds}--\ref{sec:multiSLV}
gives
\[
 |E_{N,K}|\lesssim K^{-\beta}.
\]
At the terminal depth,
\[
 M\asymp NJ_{K,\rho}
 \asymp NK^{\rho-1}\log(2+K).
\]
Putting
\[
 a:=1+(\rho-1)\kappa,
\]
we have $M\asymp N^a\log N$.  Thus, for every $\varepsilon>0$ and all
sufficiently large $N$,
\[
 cN^a\leq M\leq C_\varepsilon N^{a+\varepsilon},
 \qquad
 N\geq c_\varepsilon M^{1/(a+\varepsilon)}.
\]
It follows that
\begin{equation}\label{eq:propagated-power-calculation}
 K^{-\beta}=N^{-\beta\kappa}
 \leq C_\varepsilon
 M^{-\beta\kappa/(a+\varepsilon)}.
\end{equation}
Consequently, every exponent
\begin{equation}\label{eq:generalpropagatedexponent}
 \delta<\frac{\beta\kappa}{1+(\rho-1)\kappa}
\end{equation}
is admissible.  

We record the optimized form of this calculation.

\begin{theorem}[Abstract exponent conversion]\label{thm:abstractexponentconversion}
Fix $\beta\in(0,1)$ and $\rho>2$.  Assume the following inputs uniformly on
every signed projective chart:
\begin{enumerate}[(1)]
\item The $i$th primed coordinate product has threshold $L^{-c_{1,i}m}$ and
an SSV cover with constants $c_{2,i}$ and $c_{3,i}$ satisfying the strict
margin \eqref{eq:exactSSVmargin}.
\medskip
\item For every fixed slope parameter $t$, the double-primed product has an
SLV witness $\Gamma_{m,t}$ of measure at least
$c_\Gamma L^{-(1-\eta)m}$ and pointwise loss at most $L^{-C_1m}$ on
$\Gamma_{m,t}-\Gamma_{m,t}$, where $c_\Gamma$, $\eta$, and $C_1$ are
independent of $t$.
\end{enumerate}
Put
\[
 D:=C_1+\sum_{i=1}^d c_{1,i}.
\]
Then every exponent
\begin{equation}\label{eq:optimizedgeneralexponent}
 \delta<
 \frac{\beta\eta}{2D+(\rho-1)\eta}
\end{equation}
is admissible.
\end{theorem}

\begin{proof}
Suppose first that $D>0$.  Fix
$0<\kappa<\eta/(2D)$.  The interval
\[
 \frac\kappa\eta<s_0<\frac1{2D}
\]
is nonempty, so $s_0$ may be chosen to satisfy both
\eqref{eq:SLVscalecondition} and \eqref{eq:Fourierscalecondition}.  We now
briefly indicate how the preceding sections enter the exceptional-set
estimate.  If $|E_{N,K}|$ were larger than a sufficiently large multiple of
$K^{-\beta}$, the Fourier pigeonholing in Subsection
\ref{subsec:countingandtrig}, followed by Corollary
\ref{cor:ssvgooddirection}, would leave a direction for which the origin and
SSV errors are both $O(KL^{-n})$.  Proposition
\ref{prop:keyproposition} would then give the lower bound
\[
 KL^{-n}L^{-2Dm}.
\]
The upper bound obtained from the Fourier pigeonholing is $O(Km/(NL^n))$.
Their comparison is exactly \eqref{eq:finalparametercontradiction}, which is
impossible because $2Ds_0<1$.  Hence
\[
 |E_{N,K}|\lesssim K^{-\beta}.
\]
By \eqref{eq:propagated-power-calculation}, every
\[
 \delta<F(\kappa)
 :=\frac{\beta\kappa}{1+(\rho-1)\kappa}
\]
is admissible.  Since
\[
 F'(\kappa)
 =\frac{\beta}{[1+(\rho-1)\kappa]^2}>0,
\]
letting $\kappa\uparrow\eta/(2D)$ gives
\[
 \sup_{\kappa<\eta/(2D)}F(\kappa)
 =\frac{\beta\eta}{2D+(\rho-1)\eta}.
\]

If $D=0$, the Fourier restriction \eqref{eq:Fourierscalecondition} is
automatic.  For any fixed $\kappa>0$, choose $s_0>\kappa/\eta$ and repeat the
same argument.  Letting $\kappa\to\infty$ in
\eqref{eq:generalpropagatedexponent} gives every
$\delta<\beta/(\rho-1)$, which is again
\eqref{eq:optimizedgeneralexponent} with $D=0$.
\end{proof}

When $\beta$ approaches one and $\rho$ approaches two, Theorem
\ref{thm:abstractexponentconversion} then explicitly gives every exponent strictly below
\begin{equation}\label{eq:improvedpowerceiling}
 \frac{\eta}{2D+\eta}.
\end{equation}
The constants in this expression must come from proved estimates which are
uniform on every signed chart.  In particular, this means that any statement about the typical size of a Riesz product cannot be substituted for the interval covering in
the SSV argument, since the exponent is especially sensitive to the structure of this cover. This distinction will be essential in the four-corner
example.

Every endpoint in this subsection is a supremum rather than an attained
value.  The SSV covering argument requires a strict margin,
$\beta$ must remain below one, and the combinatorial theorem requires
$\rho>2$.  The logarithm in $J_{K,\rho}$ creates a further arbitrarily small
loss when the estimate is rewritten in terms of $M$.  These are precisely the
losses absorbed by the parameter $u>0$ in the final theorem statements.

\subsection{The sub-power regime}

We now suppose that $A_i^{(3)}$ is nontrivial for at least one coordinate.
For such a coordinate, Theorem \ref{thm:onedimensionalSSV} supplies only a
logarithmic threshold
\[
 L^{-c_{1,i}m\log_Lm}.
\]
If one kept $m\asymp\log_LN$, the squared pointwise loss would be smaller
than every fixed power of $N^{-1}$, and the comparison
\eqref{eq:finalparametercontradiction} could not close.  We therefore shorten
the low-frequency product by a factor of $\log_L\log_LN$ and take a
sub-power stacking threshold.

Let
\[
 \mathcal I_{\log}:=\{i:A_i^{(3)}\not\equiv1\},
 \qquad
 C_{\log}:=\sum_{i\in\mathcal I_{\log}}c_{1,i},
\]
where $c_{1,i}$ is the logarithmic threshold exponent for the $i$th active
coordinate.  For the remaining coordinates we retain ordinary SSV
thresholds.  Put
\[
 H:=\log_L\log_LN,
 \qquad
 K=N^{\kappa/H},
 \qquad
 m=\left\lceil s_0\frac{\log_LN}{H}\right\rceil.
\]
The SLV comparison is unchanged and requires
\begin{equation}\label{eq:log-SLV-condition}
 \kappa<\eta s_0.
\end{equation}

To compute the logarithmic SSV loss, put $B:=\log_LN$, so $H=\log_LB$.
Then
\[
 m=\frac{s_0B}{H}+O(1),
 \qquad
 \log_Lm=H-\log_LH+O(1),
\]
and consequently
\[
 m\log_Lm
 =s_0B+O\left(\frac{B\log H}{H}\right).
\]
It follows that
\begin{align}
 \prod_{i\in\mathcal I_{\log}}
 L^{-2c_{1,i}m\log_Lm}
 &=N^{-2s_0C_{\log}+o(1)},
 \label{eq:log-ssv-asymptotic}\\
 L^{-2m(C_1+\sum_{i\notin\mathcal I_{\log}}c_{1,i})}
 &=N^{-o(1)}.
 \label{eq:log-slv-asymptotic}
\end{align}
Thus the logarithmic SSV coordinates contain the leading power of $N$,
whereas the accumulating loss and all ordinary SSV losses contribute only to
the lower-order term.  The Fourier lower bound contradicts the upper bound
provided that
\begin{equation}\label{eq:log-Fourier-condition}
 2s_0C_{\log}<1.
\end{equation}
Combining \eqref{eq:log-SLV-condition} and
\eqref{eq:log-Fourier-condition}, we obtain a nonempty interval of choices for
$s_0$ exactly when
\begin{equation}\label{eq:optimizedlogcoefficient}
 \kappa<\frac{\eta}{2C_{\log}}.
\end{equation}

It remains to verify that propagation does not alter the leading
$1/\log\log N$ coefficient.  Since
\[
 K^{\rho-1}=N^{(\rho-1)\kappa/H}=N^{o(1)}
 \quad\text{and}\quad
 \log(2+K)=N^{o(1)},
\]
we have
\[
 NJ_{K,\rho}=N^{1+o(1)}.
\]
Hence $\log_LM=(1+o(1))\log_LN$ and
$\log_L\log_LM=H+o(1)$.  If
$|E_{N,K}|\lesssim K^{-\beta}$, the propagated estimate therefore has the
form
\[
 \Fav\bigl(\mathcal N_{L^{-M}}(\mathcal S^\infty)\bigr)
 \leq
 M^{-c/\log_L\log_LM}
\]
for every
\[
 0<c<\frac{\beta\eta}{2C_{\log}}.
\]
In particular, this proves a quantitative sub-power estimate of exactly the
form asserted in Theorem \ref{thm:totalFavbound}.  Notice that the parameter
$\rho$ affects only the $o(1)$ term in this branch.

\subsection{Completion of the general theorem}

We now verify that the hypotheses used above are precisely the hypotheses
available for Theorem \ref{thm:totalFavbound}.  For every coordinate,
Theorem \ref{thm:onedimensionalSSV} supplies log-SSV.  If
$A_i^{(3)}\equiv1$, it supplies ordinary SSV as well.  Its final assertion
allows the threshold to be weakened until the strict covering margin
\eqref{eq:exactSSVmargin} holds on every member of the finite signed chart
family.

The algebraic alternatives in Theorem \ref{thm:totalFavbound} are used only
to control the accumulating factors.  Lemma \ref{lma:singlescaleSLV}, which
is the single-scale result from \cite{LM}, treats the small-cardinality and
at-most-two-prime cases.  Lemma \ref{lma:singlescaleSLVfibered} treats the
fibered case; its cardinality input is Proposition
\ref{prop:multifibered}, quoted from \cite[Proposition 3.2]{KLMS}.  Trivial
accumulating factors and singleton coordinates are omitted as explained
before Lemma \ref{lma:singlescaleSLV}.  Lemma
\ref{lma:singletomulti} then constructs, for every fixed slope $t$, a witness
$\Gamma_{m,t}$ whose measure and pointwise constants are uniform in all slope
parameters, including slopes equal or arbitrarily close to zero.

Corollary \ref{cor:ssvgooddirection} selects a direction in the
Fourier-good parameter set for which the SSV error is bounded.  Proposition
\ref{prop:shortintervalbound}, whose auxiliary localization estimate is
cited there from \cite{BThesis,BV3}, controls the interval about the origin
using the same symbolic multiplicities as the Fourier sum.  We choose the
implicit constant in $L^{\eta m}\gtrsim K$ large enough to absorb both error
terms.  Proposition \ref{prop:keyproposition} then supplies the Riesz-product
lower bound, and the power or sub-power choice of $m$ gives
$|E_{N,K}|\lesssim K^{-\beta}$ for every fixed $\beta<1$ after an arbitrarily
small reduction of the scale exponent.

Finally, Theorem \ref{thm:improved-L2} applies for every $\rho>2$, Lemma
\ref{lma:reverseholder} propagates the high-multiplicity alternative, and
Proposition \ref{prop:favardreducedtoexceptional} passes from the exceptional
set to Favard length at every sufficiently large construction depth.  

This
completes the proof of Theorem \ref{thm:totalFavbound}.  Theorems
\ref{thm:fav2prime} and \ref{thm:fibered} are its corresponding
specializations, which we discuss in the final few sections of this work.

\subsection{The four-corner Cantor set}

Let $\mathcal K_2^\infty$ denote the classical four-corner Cantor set.  Its
coordinate digit set is $A=\{0,3\}$, and hence
\[
 A(X)=1+X^3,
 \qquad
 |\phi_A(u)|=|\cos(3\pi u)|.
\]
At the symmetric slope $t = (1,1)$, which corresponds to an angle of $\pi /4$ radians, the two coordinate factors coincide, so the full
low-frequency product is
\[
 \prod_{k=0}^{m-1}|\cos(3\pi4^ku)|^2.
\]
The classical logarithmic cosine integral (see, for example,
\cite[Chapter V]{Zygmund}) gives
\[
 \int_0^1-\log_4|\cos(3\pi u)|^2\,du=1.
\]
Thus $4^{-m}$ is the geometric-mean scale of the full low product.  This
identity does \emph{not} imply that the sublevel set below a constant multiple
of $4^{-m}$ admits the exponentially efficient interval covering required by
ordinary SSV.  In particular, the formal substitution $D=1$ (which would yield a final exponent of $1/4 -u$) in
\eqref{eq:improvedpowerceiling} \textit{is not justified} by our current methods.

We instead retain the specialized Fourier argument of Nazarov, Peres, and
Volberg and replace only its combinatorial propagation.  Their exact
factorization gives
\begin{equation}\label{eq:NPV-low-product-bound}
 |P_{\mathrm{lo},t}(y)|
 \gtrsim
 4^{-2m}
 |\sin(4^m(y+ty))|\,
 |\sin(4^m(y-ty))|,
\end{equation}
where $4^m$ is a sufficiently large constant multiple of $K$; see
\cite[pp. 8--11]{NPV}.  Their Fourier analysis uses two independent
parameters.  The number $K$ is the multiplicity threshold, while $S^{-1}$
measures the exceptional set of directions.  The two parameters must remain
separate throughout the calculation.

The analytic portion of \cite{NPV} closes at a depth larger than
$(KS)^{4+\gamma}$ for any fixed $\gamma>0$.  Their original combinatorial
iteration then used $K^{2+o(1)}$ blocks and led to the conditions $p>6$ and
$q>4$; see \cite[pp. 13--15]{NPV}.  The improved propagation replaces the
factor $K^{2+o(1)}$ by $K^{1+o(1)}$ and yields the following two-parameter
estimate.

\begin{proposition}[Two-parameter four-corner estimate]\label{prop:fourcornertwoparameter}
Let $p>5$ and $q>4$.  For all sufficiently large $K$ and $S$, if
\[
 M\geq K^pS^q,
\]
then
\[
 \left|\left\{\theta\in(0,\pi/4):
 |\proj_\theta\mathcal K_2^M|\geq\frac{C}{K}\right\}\right|
 \lesssim_{p,q}\frac1S.
\]
\end{proposition}

\begin{proof}
Fix $\gamma>0$.  The Fourier portion of the Nazarov--Peres--Volberg
argument gives a contradiction at an analytic depth $N_0$ once
\[
 N_0>(KS)^{4+\gamma};
\]
this is the conclusion of the calculation on
\cite[pp. 10--11]{NPV}.  Choose $\rho=2+\gamma$.  By
\eqref{eq:propagationblocks},
\[
 J_{K,\rho}
 \lesssim K^{1+\gamma}\log(2+K)
 \lesssim_\gamma K^{1+2\gamma}
\]
for all sufficiently large $K$.  The terminal depth required by Lemma
\ref{lma:reverseholder} is therefore at most
\[
 N_0J_{K,\rho}
 \lesssim_\gamma
 K^{5+3\gamma}S^{4+\gamma}.
\]
Given $p>5$ and $q>4$, choose $\gamma>0$ so small that
\[
 5+3\gamma<p,
 \qquad
 4+\gamma<q.
\]
Thus $M\geq K^pS^q$ is sufficient for the improved terminal depth.  The
marking and self-similarity argument now gives a projection bound $C/K$
outside a set of directions of measure $O(S^{-1})$, as required.
\end{proof}

We now are in a position to prove the following improvement of the Nazarov, Peres and Volberg exponent for the four-corner Cantor set.

\begin{theorem}[The four-corner Cantor set]\label{thm:fourcorneronefifth}
For every $u>0$, there are constants $C_u$ and $N_u$ such that
\begin{equation}\label{eq:fourcorneronefifth}
 \Fav\bigl(\mathcal N_{4^{-N}}(\mathcal K_2^\infty)\bigr)
 \leq C_uN^{-1/5+u}
\end{equation}
for every $N\geq N_u$.
\end{theorem}

\begin{proof}
Fix $p>5$ and choose $q$ with
\[
 4<q<p.
\]
We first convert Proposition \ref{prop:fourcornertwoparameter} into a
distributional estimate for the projection length.  Put
\[
 F_M(\theta):=|\proj_\theta\mathcal K_2^M|,
 \qquad 0<\theta<\pi/4.
\]
The functions $F_M$ are bounded by an absolute constant $C_0$.  Fix a
sufficiently large constant $K_0$.  If
$
 c_0M^{-1/p}\leq t\leq C/K_0,
$
set
\[
 K=\frac Ct,
 \qquad
 S=cM^{1/q}t^{p/q},
\]
where $c>0$ is chosen so that $K^pS^q\leq M$.  Choosing $c_0$ sufficiently
large makes $S$ admissible.  Proposition
\ref{prop:fourcornertwoparameter} then gives
\begin{equation}\label{eq:four-corner-distribution}
 |\{\theta\in(0,\pi/4):F_M(\theta)\geq t\}|
 \lesssim M^{-1/q}t^{-p/q}.
\end{equation}
For $t>C/K_0$, we apply the proposition with $K=K_0$ and
$S=cM^{1/q}$.  Since $F_M\leq C_0$, the resulting bound
$O(M^{-1/q})$ is again dominated by the right-hand side of
\eqref{eq:four-corner-distribution}, after changing its implicit constant.
For $0<t<c_0M^{-1/p}$, we use the trivial bound $\pi/4$.

The layer-cake identity and \eqref{eq:four-corner-distribution} now yield
\begin{align*}
 \int_0^{\pi/4}F_M(\theta)\,d\theta
 &=\int_0^{C_0}
 |\{\theta\in(0,\pi/4):F_M(\theta)\geq t\}|\,dt\\
 &\lesssim M^{-1/p}
 +M^{-1/q}\int_{c_0M^{-1/p}}^{C_0}t^{-p/q}\,dt.
\end{align*}
Because $p/q>1$, the last integral is controlled by its lower endpoint:
\begin{align*}
 M^{-1/q}\int_{c_0M^{-1/p}}^{C_0}t^{-p/q}\,dt
 &\lesssim
 M^{-1/q}(M^{-1/p})^{1-p/q}\\
 &=M^{-1/p}.
\end{align*}
It follows that
\[
 \int_0^{\pi/4}|\proj_\theta\mathcal K_2^M|\,d\theta
 \lesssim M^{-1/p}.
\]
The symmetries of the four-corner construction give the same bound over the
remaining direction intervals.  Hence $\Fav(\mathcal K_2^M)\lesssim
M^{-1/p}$.  The comparison \eqref{eq:squaresandballs} then gives
\[
 \Fav\bigl(\mathcal N_{4^{-M}}(\mathcal K_2^\infty)\bigr)
 \lesssim M^{-1/p}.
\]
Finally, choose $p>5$ sufficiently close to $5$ that
$1/p>1/5-u$, and rename $M$ as $N$.  This proves
\eqref{eq:fourcorneronefifth}.
\end{proof}

The exponent $1/5$ is a direct numerical consequence of the two
main parts of the Nazarov-Peres-Volberg argument: the Fourier estimate contributes four powers of
$K$, while the improved propagation contributes only one (as opposed to two, as in the original work of Nazarov, Peres and Volberg in \cite{NPV}).  Indeed, the earlier
$K^{2+o(1)}$ propagation produced the denominator six; the present
$K^{1+o(1)}$ propagation reduces it to five. Anything beyond this $1/5$ barrier would require either a significant reworking of the combinatorial arguments in Section \ref{sec:comblemmata}, or else another perspective replacing the current Fourier analytic arguments developed in Sections \ref{sec:toolbox} through \ref{sec:multiSLV}.

\subsection{General zero-free digit sets}

We next consider the opposite analytic extreme to the previous four-corner set example.  Suppose that every
coordinate mask polynomial is zero-free on the unit circle, and put
\[
 \mu_{A_i}:=\min_{|z|=1}|A_i(z)|>0
 \qquad(1\leq i\leq d).
\]
Define
\begin{equation}\label{eq:zero-free-coordinate-loss}
 d_i:=\log_L\frac{\#A_i}{\mu_{A_i}}\geq0.
\end{equation}
By the definition of the normalized mask polynomial,
\[
 |\phi_{A_i}(u)|
 =\frac{|A_i(e^{2\pi iu})|}{\#A_i}
 \geq\frac{\mu_{A_i}}{\#A_i}
 =L^{-d_i}
 \qquad(u\in\mathbb R).
\]
Consequently,
\begin{equation}\label{eq:zerofree-coordinate-product}
 \prod_{k=0}^{m-1}|\phi_{A_i}(L^ku)|
 \geq L^{-d_im}
 \qquad(u\in\mathbb R).
\end{equation}

There are no unit-circle roots, so $A_i'=A_i$ and $A_i''\equiv1$.  In the
notation of \eqref{eq:zerofreeD}, the limiting total pointwise loss is
\begin{equation}\label{eq:zero-free-total-loss}
 D:=\sum_{i=1}^dd_i
 =\sum_{i=1}^d\log_L\frac{\#A_i}{\mu_{A_i}}.
\end{equation}
There is one minor endpoint issue.  Equality may occur in
\eqref{eq:zerofree-coordinate-product}, so the sublevel set at the threshold
$L^{-d_im}$ need not be empty.  Fix $\varepsilon>0$ and choose positive
$\varepsilon_i$ with $\sum_i\varepsilon_i=\varepsilon$.  At the strictly
smaller threshold
\[
 \psi_i(m):=L^{-(d_i+\varepsilon_i)m},
\]
the coordinate SSV set is empty.  It may therefore be covered with any fixed
covering exponents satisfying \eqref{eq:exactSSVmargin}, and the total loss in
Theorem \ref{thm:abstractexponentconversion} is
\[
 D_\varepsilon
 :=\sum_{i=1}^d(d_i+\varepsilon_i)
 =D+\varepsilon.
\]

Moreover, by the zero-root hypothesis, the double-primed product is identically one.  Thus we may take
$
\Gamma_{m,t}=[0,c_*/2]
$
for every $m$ and $t$, with $C_1=0$, $\eta=1$, and
$c_\Gamma=c_*/2$.  There is no deleted SSV set and no nontrivial SLV
construction.  

The abstract exponent conversion therefore gives the
following favourable result.

\begin{corollary}[Zero-free digit sets]\label{cor:zerofreedigitsets}
Let $L\geq4$, let $A_1,\ldots,A_d\subset\{0,\ldots,L-1\}$ satisfy
\[
 \prod_{i=1}^d\#A_i=L,
 \qquad
 \#A_i,\#A_j\geq2\quad\text{for some }i\ne j,
\]
and let $\mathcal S^\infty$ be the associated rational product Cantor set.
Assume that every coordinate mask polynomial is zero-free on the unit circle.
Then, for every $u>0$,
\[
 \Fav\bigl(\mathcal N_{L^{-N}}(\mathcal S^\infty)\bigr)
 \lesssim_{A_1,\ldots,A_d,d,u}
 N^{-1/(2D+1)+u},
\]
where $D$ is defined by \eqref{eq:zero-free-total-loss}.
\end{corollary}

\begin{proof}
Fix $\varepsilon>0$ and use the empty SSV sets constructed above.  Since the
SLV witness is the fixed interval $[0,c_*/2]$, Theorem
\ref{thm:abstractexponentconversion} gives every exponent
\[
 \delta<
 \frac{\beta}{2(D+\varepsilon)+(\rho-1)}.
\]
Let $\varepsilon\downarrow0$, $\beta\uparrow1$, and $\rho\downarrow2$.
Every $\delta<1/(2D+1)$ is therefore admissible, and the arbitrarily small
difference from the supremum is absorbed into $u$.
\end{proof}

This formula makes the gain transparent.  Once unit-circle zeros disappear,
the exponent is no longer constrained by a covering theorem for small values.
It depends only upon the pointwise minimum moduli of the finitely many mask
polynomials.

In particular, by using the abstract exponent conversion method from this section, this means that we can use a computer search to find polynomials with large modulus on the unit circle, and then immediately obtain improved power-law exponents for the associated symmetric rational product Cantor set. We give two such examples now.

\subsubsection{An explicit zero-free example}

Take
\[
 L=25,
 \qquad
 A_1=A_2=A:=\{0,1,2,6,9\},
\]
and write
\[
 P(z):=A(z)=1+z+z^2+z^6+z^9.
\]
The following elementary certificate gives the exact lower bound needed for
Corollary \ref{cor:zerofreedigitsets}.

\medskip
\begin{lemma}\label{lma:explicitzerofreedigitset}
For every $|z|=1$,
\[
 |P(z)|\geq1.
\]
Equality holds at $z=-1$.
\end{lemma}

\begin{proof}
Write $z=e^{i\theta}$ and $x=\cos\theta$.  Expanding and collecting the
trigonometric terms gives the exact factorization
\begin{equation}\label{eq:explicit-zero-free-factorization}
 |P(e^{i\theta})|^2-1=(x+1)Q(x),
\end{equation}
where
\[
\begin{aligned}
 Q(x)={}&512x^8-256x^7-768x^6+320x^5+352x^4\\
 &\quad-112x^3-48x^2+8x+4.
\end{aligned}
\]
Using this factorization, we can thus directly verify that $Q$ is positive on $[-1,1]$.  For standard background on Bernstein polynomials, including the positivity principle and basis conversion used in the exact certificates below, see \cite{Lorentz}.
For $0\leq j\leq15$, let
\[
 I_j:=\left[-1+\frac j8,-1+\frac{j+1}{8}\right]
\]
and put $q_j(t):=Q(-1+(j+t)/8)$ for $0\leq t\leq1$.  Write $q_j$ in the
degree-eight Bernstein basis,
\[
 q_j(t)=\sum_{r=0}^8b_{j,r}
 \binom8r t^r(1-t)^{8-r}.
\]
Every Bernstein basis function in this display is nonnegative on $[0,1]$,
and their sum is one.  Thus it is enough to verify that all coefficients
$b_{j,r}$ are positive.  For reproducibility, if
$q_j(t)=\sum_{k=0}^8a_{j,k}t^k$, then
\[
 b_{j,r}
 =\sum_{k=0}^r
 \frac{\binom rk}{\binom8k}a_{j,k}.
\]
Exact rational expansion gives the following minimum coefficient on each
subinterval:
\[
\begin{array}{@{}cr@{\qquad}cr@{\qquad}cr@{\qquad}cr@{}}
\toprule
 j&\min_r b_{j,r}&j&\min_r b_{j,r}&j&\min_r b_{j,r}&j&\min_r b_{j,r}\\
\midrule
0&303741/32768&4&100781/32768&8&4&12&2\\
1&41207/28672&5&211/128&9&351/128&13&443/128\\
2&311/128&6&1673/1024&10&48341/32768&14&13381/32768\\
3&6&7&83301/32768&11&79039/57344&15&15137/65536
\\
\bottomrule
\end{array}
\]
In particular,
\[
 Q(x)\geq\frac{15137}{65536}>0
 \qquad(-1\leq x\leq1).
\]
Since $x+1\geq0$ on this interval, the factorization
\eqref{eq:explicit-zero-free-factorization} proves $|P(e^{i\theta})|\geq1$.
At $z=-1$, or equivalently $x=-1$, equality holds.
\end{proof}

We now complete the proof of the explicit theorem stated in Section
\ref{subsec:RPCS}.

\begin{proof}[Proof of Theorem \ref{thm:explicitonethird}]
Lemma \ref{lma:explicitzerofreedigitset} gives $\mu_A=1$.  After normalizing
by $\#A=5$,
\[
 |\phi_A(u)|\geq\frac15=25^{-1/2}.
\]
Thus each coordinate contributes loss $1/2$, and the planar product has
\[
 D=\frac12+\frac12=1.
\]
Corollary \ref{cor:zerofreedigitsets} gives every exponent strictly below
\[
 \frac1{2D+1}=\frac13,
\]
which is exactly the conclusion of Theorem \ref{thm:explicitonethird}.
\end{proof}

This example exhibits the full benefit of the method.  The analytic deletion
is empty, the accumulating factor is absent, and the improved propagation is
the only remaining structural loss.  The resulting exponent $1/3-u$ is
therefore obtained by a direct calculation of the behaviour of minimal moduli of the associated mask polynomials on the unit circle (which is just a simple and direct computational verification).

\subsubsection{Crossing the one-third exponent via direct example}

The exponent $1/3$ in the preceding example is not a universal barrier of
the method.  Indeed, in the symmetric planar setting $A_1=A_2=A$, with
$q:=\#A$ and $L=q^2$, the total loss is
\begin{equation}\label{eq:symmetric-zero-free-loss}
 D=\log_q\frac q{\mu_A},
 \qquad
 \mu_A:=\min_{|z|=1}|A(z)|.
\end{equation}
The endpoint furnished by Corollary \ref{cor:zerofreedigitsets} exceeds
$1/3$ precisely when
\[
 \frac1{2D+1}>\frac13
 \quad\Longleftrightarrow\quad
 D<1
 \quad\Longleftrightarrow\quad
 \mu_A>1.
\]
Thus the base-$25$ example, for which $\mu_A=1$, lies exactly on the
one-third threshold.  To cross that threshold it is enough to construct a
mask polynomial whose minimum modulus is strictly larger than one.

We now give such an example.  Let
\[
 A:=\{0,7,14,15,18,24\}\subset\{0,\ldots,35\}
\]
and put
\[
 P_6(z):=1+z^7+z^{14}+z^{15}+z^{18}+z^{24}.
\]

\begin{lemma}\label{lma:explicitbeyondthird}
For every $|z|=1$,
\begin{equation}\label{eq:base36-minimum-modulus}
 |P_6(z)|\geq\frac{53}{50}.
\end{equation}
\end{lemma}

\begin{proof}
Write $z=e^{i\theta}$, $x=\cos\theta$, and let
$T_k(x):=\cos(k\arccos x)$ denote the $k$th Chebyshev polynomial.  For standard background on these polynomials and their trigonometric identities, see \cite{Rivlin}.  The
positive differences between pairs of digits in $A$, counted with
multiplicity, are
\[
 1,3,4,6,7,7,8,9,10,11,14,15,17,18,24.
\]
Consequently,
\begin{equation}\label{eq:base36-chebyshev-certificate}
\begin{aligned}
 R(x)
 &:=
 |P_6(e^{i\theta})|^2-\left(\frac{53}{50}\right)^2\\
 &=\frac{12191}{2500}
 +2\bigl(T_1+T_3+T_4+T_6+2T_7+T_8+T_9+T_{10}+T_{11}\\
 &\hspace{38mm}{}+T_{14}+T_{15}+T_{17}+T_{18}+T_{24}\bigr)(x).
\end{aligned}
\end{equation}
We certify the positivity of this degree-$24$ polynomial using exact rational
arithmetic.  If $I=[a,b]\subset[-1,1]$, write
\[
 R(a+(b-a)t)
 =\sum_{r=0}^{24}b_{I,r}\binom{24}{r}t^r(1-t)^{24-r},
 \qquad 0\leq t\leq1.
\]
If the monomial expansion on the left is
$\sum_{k=0}^{24}a_{I,k}t^k$, then
\begin{equation}\label{eq:base36-bernstein-conversion}
 b_{I,r}
 =\sum_{k=0}^r
 \frac{\binom rk}{\binom{24}{k}}a_{I,k}.
\end{equation}
The following table records a rational lower bound for the smallest
Bernstein coefficient on each member of a partition of $[-1,1]$.  Every
entry follows by direct substitution into
\eqref{eq:base36-chebyshev-certificate} and
\eqref{eq:base36-bernstein-conversion}, and we have included the table as direct verification for the exponent which is presented immediately following it.
\[
\begin{array}{@{}cr@{\qquad}cr@{}}
\toprule
 I&\min_r b_{I,r}>&I&\min_r b_{I,r}>\\
\midrule
\bigl[-1,-63/64\bigr]&1/100&\bigl[-1/2,-3/8\bigr]&1/2\\
\bigl[-63/64,-31/32\bigr]&2&\bigl[-3/8,-1/4\bigr]&1\\
\bigl[-31/32,-15/16\bigr]&2&\bigl[-1/4,0\bigr]&2\\
\bigl[-15/16,-7/8\bigr]&3/5&\bigl[0,1/4\bigr]&2/5\\
\bigl[-7/8,-13/16\bigr]&3/2&\bigl[1/4,1/2\bigr]&1/6\\
\bigl[-13/16,-25/32\bigr]&1/20&\bigl[1/2,3/4\bigr]&2/3\\
\bigl[-25/32,-3/4\bigr]&1/250&\bigl[3/4,7/8\bigr]&2/5\\
\bigl[-3/4,-1/2\bigr]&3/10&\bigl[7/8,15/16\bigr]&4/5\\
&&\bigl[15/16,31/32\bigr]&1/3\\
&&\bigl[31/32,1\bigr]&1/4
\\
\bottomrule
\end{array}
\]
The Bernstein basis functions are nonnegative and sum to one on $[0,1]$.
The table therefore proves $R(x)>0$ throughout $[-1,1]$, and
\eqref{eq:base36-minimum-modulus} follows.
\end{proof}

\begin{theorem}[An explicit exponent beyond one-third]\label{thm:explicitbeyondthird}
Let $\mathcal S_{36}^\infty\subset\mathbb R^2$ be the rational product
Cantor set in base $36$ with
\[
 A_1=A_2=\{0,7,14,15,18,24\}.
\]
Put
\begin{equation}\label{eq:base36-exponent}
 \delta_{36}:=
 \frac1{2\log_6(300/53)+1}
 =0.340720\ldots.
\end{equation}
Then, for every $u>0$,
\[
 \Fav\bigl(\mathcal N_{36^{-N}}(\mathcal S_{36}^\infty)\bigr)
 \lesssim_u N^{-\delta_{36}+u}
\]
for all sufficiently large $N$.  In particular, $\delta_{36}>1/3$.
\end{theorem}

\begin{proof}
Lemma \ref{lma:explicitbeyondthird} gives $\mu_A\geq53/50$.  By
\eqref{eq:symmetric-zero-free-loss},
\[
 D=\log_6\frac6{\mu_A}
 \leq\log_6\frac{300}{53}.
\]
Corollary \ref{cor:zerofreedigitsets} therefore gives every exponent strictly
below $\delta_{36}$.  Moreover,
\[
 \delta_{36}>\frac13
 \quad\Longleftrightarrow\quad
 \log_6\frac{300}{53}<1
 \quad\Longleftrightarrow\quad
 \frac{300}{53}<6,
\]
and the last inequality is simply $300<318$.
\end{proof}

\begin{remark}
Theorem \ref{thm:explicitbeyondthird} shows that the exponent $1/3$ is a
threshold for the particular five-digit example, rather than a barrier of
the zero-free mechanism.  The systematic construction of digit sets with
still larger minimum modulus, and the resulting optimization of the exponent,
will be pursued in a direct follow-up to the present work.
\end{remark}

\subsection{Further zero-free examples}

Several symmetric examples illustrate how the exponent changes with the
minimum modulus.  Let $A_1=A_2=A$, put $q:=\#A$ and $L=q^2$, and define
\[
 \mu_A:=\min_{|z|=1}|A(z)|.
\]
Since
\[
 2\log_{q^2}\frac q{\mu_A}
 =\log_q\frac q{\mu_A},
\]
the total loss and the corresponding exponent ceiling are
\[
 D=\log_q\frac q{\mu_A},
 \qquad
 \delta_*:=\frac1{2D+1}.
\]

\begin{center}
\renewcommand{\arraystretch}{1.15}
\begin{tabular}{@{}c c c c@{}}
\toprule
$A$ & $L$ & $\mu_A$ & Exponent ceiling $\delta_*$\\
\midrule
$\{0,1,3\}$ & $9$ & $0.607346\ldots$ & $0.2558989\ldots$\\
$\{0,1,2,4\}$ & $16$ & $0.752394\ldots$ & $0.2932173\ldots$\\
$\{0,1,2,6,9\}$ & $25$ & $1$ & $1/3$\\
$\{0,7,14,15,18,24\}$ & $36$ & $\geq1.06$ & $\geq0.3407202\ldots$\\
\bottomrule
\end{tabular}
\end{center}

We verify the two numerical entries directly, so that neither example depends upon an uncertified floating-point search.  For $A=\{0,1,3\}$ and
$x=\cos\theta$,
\[
 |1+e^{i\theta}+e^{3i\theta}|^2
 =8x^3+4x^2-4x+1.
\]
The derivative is $4(6x^2+2x-1)$, whose two zeros are
\[
 \frac{-1-\sqrt7}{6}
 \quad\text{and}\quad
 \frac{-1+\sqrt7}{6}.
\]
Evaluation at these points and at $x=\pm1$ shows that the global minimum is
attained at $x=(\sqrt7-1)/6$ and equals
\[
 \mu_A^2=\frac{47-14\sqrt7}{27}.
\]
This gives the first line of the table exactly.

For $A=\{0,1,2,4\}$, direct expansion gives
\[
 |1+e^{i\theta}+e^{2i\theta}+e^{4i\theta}|^2
 =f(x):=16x^4+8x^3-8x^2-2x+2.
\]
The interior critical points are the zeros of
\[
 g(x):=32x^3+12x^2-8x-1,
 \qquad f'(x)=2g(x).
\]
Since
\[
 g'(x)=8(12x^2+3x-1),
\]
the cubic has at most one zero in each of the three monotonicity intervals
determined by $(-3\pm\sqrt{57})/24$.  Exact substitution at the terminating
decimal endpoints below gives alternating signs, so the three zeros
$r_1<r_2<r_3$ satisfy
\[
\begin{aligned}
 -0.676326&<r_1<-0.676325,\\
 -0.111835&<r_2<-0.111834,\\
 \phantom{-}0.413160&<r_3<\phantom{-}0.413161.
\end{aligned}
\]
These sign changes, together with the monotonicity just noted, prove both
existence and uniqueness in the displayed intervals.  Rational interval
evaluation of $f$ gives
\[
\begin{aligned}
 0.566096&<f(r_1)<0.566097,\\
 2.1149&<f(r_2)<2.1150,\\
 0.8384&<f(r_3)<0.8386.
\end{aligned}
\]
Since $f(-1)=4$ and $f(1)=16$, the global minimum occurs at $r_1$.
Moreover,
\[
 0.752393^2<0.566096<f(r_1)<0.566097<0.752395^2,
\]
and hence
\[
 0.752393<\mu_A<0.752395.
\]
Substitution into $D=\log_q(q/\mu_A)$ gives the second exponent displayed in
the table.

\subsection{Curvilinear consequences of the linear estimates}

We now record the consequence of the comparison principle of {\L}aba,
McDonald, and Taylor \cite{LMT}.  Let $E_N$ be the $N$th planar approximation
of a purely unrectifiable self-similar $1$-set satisfying the strong separation
and non-collinearity assumptions in \cite[Section 2.2]{LMT}.  If
$\Phi_\alpha$ is a smooth nonlinear projection family satisfying the uniform
regularity and transversality hypotheses of \cite[Theorem 2.2]{LMT}, define
its generalized Favard length by
\[
 \Fav_\Phi(E_N):=\int_I |\Phi_\alpha(E_N)|\,d\alpha.
\]
Their global comparison theorem gives
\begin{equation}\label{eq:LMT-global-comparison}
 \Fav_\Phi(E_N)
 \lesssim_\Phi
 \Fav(E_{\lfloor N/2\rfloor});
\end{equation}
see \cite[Theorem 2.4]{LMT}.  It follows immediately that any classical
estimate
\[
 \Fav(E_N)\lesssim_u N^{-\delta+u}
\]
passes to the corresponding curvilinear problem with the same power exponent:
\begin{equation}\label{eq:LMT-power-transfer}
 \Fav_\Phi(E_N)\lesssim_{\Phi,u}N^{-\delta+u}.
\end{equation}
The replacement of $N$ by $\lfloor N/2\rfloor$ changes only the implicit
constant and produces no exponent loss.

In particular, whenever the stated hypotheses apply, Theorem
\ref{thm:fourcorneronefifth} gives the generalized exponent $1/5-u$ for the
four-corner approximants, Theorem \ref{thm:explicitonethird} gives $1/3-u$
for the base-$25$ example, and Theorem \ref{thm:explicitbeyondthird} gives
$\delta_{36}-u>1/3-u$ for the base-$36$ example.  Thus, in these regimes, the
linear estimate is the quantitatively dominant input: an improvement in the
classical projection problem automatically governs the associated smooth
curvilinear families through \eqref{eq:LMT-global-comparison}.  This is the
precise sense in which the linear bounds dominate the nonlinear estimates
here.

\section*{Future Directions}

The four-corner argument above deliberately uses only the proved Fourier
input of Nazarov, Peres and Volberg.  The geometric-mean identity suggests a stronger estimate at the
natural scale, but an unweighted small-value bound is not sufficient: the
high-frequency Riesz-product weight may concentrate on a set of very small
Lebesgue measure.  Any improvement at this point must therefore control the
low-frequency product \emph{against the high-frequency weight}, uniformly at
rational and near-rational slopes.

The numerical effect of such an estimate can be stated cleanly.  Suppose that
a four-corner Fourier argument has analytic cost $\lambda>0$ in the following
precise sense: for every $\beta<1$ and $\varepsilon>0$, one may take
\[
 K=N^{1/(\lambda+\varepsilon)}
\]
and prove
\[
 |E_{N,K}|\lesssim_{\beta,\varepsilon}K^{-\beta}.
\]
After propagation,
\[
 M\asymp NK^{\rho-1}\log(2+K)
 =K^{\lambda+\varepsilon+\rho-1+o(1)},
\]
while Proposition \ref{prop:favardreducedtoexceptional} gives a Favard bound
$K^{-\beta}$.  Hence every
\begin{equation}\label{eq:conditional-weighted-conversion}
 \delta<\frac{\beta}{\lambda+\varepsilon+\rho-1}
\end{equation}
is admissible.  Letting $\beta\uparrow1$, $\varepsilon\downarrow0$, and
$\rho\downarrow2$ gives every
\[
 \delta<\frac1{\lambda+1}.
\]

In the four-corner set bookkeeping, a weighted estimate saving one power of $K$
over the pointwise NPV lower bound would thus have analytic cost $\lambda=3$ and
would therefore give a final exponent of $1/4-u$.  A full weighted estimate at the natural geometric scale would have cost $\lambda=2$ and would give $1/3-u$.  These
conclusions are conditional applications of
\eqref{eq:conditional-weighted-conversion}; neither weighted input is used in
the proved theorems of this paper, and these are currently unavailable to us in our current methodology.  This distinction is consistent with the
original observation in \cite[p. 3]{NPV} that passing the one-quarter scale
would require new ideas.

\clearpage
\subsection{Parameter table}\label{subsec:parametertable}

\begin{center}
\small
\renewcommand{\arraystretch}{1.32}
\begin{tabularx}{\textwidth}{@{}>{\raggedright\arraybackslash}p{0.13\textwidth}
  >{\raggedright\arraybackslash}p{0.31\textwidth}X@{}}
\toprule
Parameter & Role & Main restriction or effect\\
\midrule
$L$ & Base of the construction & Fixed by the digit data; $\prod_i\#A_i=L$.\\
$N$ & Analytic construction depth & Tends to infinity.\\
$n$ & Pigeonholed frequency scale & Chosen with $N/4\leq n\leq N/2$.\\
$K$ & Stacking threshold & In the power branch, $K=N^\kappa$.\\
$m$ & Length of the low-frequency product & In the power branch, $m=\lceil s_0\log_LN\rceil$.\\
$s_0$ & Analytic scale-selection parameter & Determines $m$ and is kept distinct from every cyclotomic index $s$.\\
$\beta$ & Exceptional-set exponent & The analytic argument seeks $|E_{N,K}|\lesssim K^{-\beta}$.\\
$c_{1,i}$ & Coordinate SSV loss & The ordinary threshold is $L^{-c_{1,i}m}$; a logarithmic threshold is $L^{-c_{1,i}m\log_Lm}$.\\
$c_{2,i}$ & SSV covering-number exponent & The $i$th coordinate SSV set is covered by $O(L^{c_{2,i}m})$ intervals.\\
$c_{3,i}$ & SSV interval-length exponent & Each covering interval has length $O(L^{-c_{3,i}m})$; the pair $(c_{2,i},c_{3,i})$ must satisfy \eqref{eq:exactSSVmargin}.\\
$C_1$ & Pointwise SLV loss & In the notation of Section \ref{sec:multiSLV}, the accumulating product is bounded below by $L^{-C_1m}$.\\
$D$ & Total pointwise loss & In the power branch, $D=C_1+\sum_i c_{1,i}$.\\
$\eta$ & SLV measure gain & The witness $\Gamma_{m,t}$ has measure at least $c_\Gamma L^{-(1-\eta)m}$, uniformly in $t$.\\
$\rho$ & Low-multiplicity threshold & The improved combinatorial argument permits every $\rho>2$.\\
$J_{K,\rho}$ & Propagation length in blocks & Comparable to $K^{\rho-1}\log(2+K)$.\\
$M$ & Terminal construction depth & Comparable to $NJ_{K,\rho}$.\\
$\delta$ & Final Favard exponent & Obtained after analytic optimization and propagation.\\
\bottomrule
\end{tabularx}
\end{center}

\newpage
\bibliographystyle{amsplain}

\end{document}